\input amstex
\documentstyle{amsppt}
\magnification=\magstep1
\vsize =21 true cm
\hsize =16 true cm
\loadmsbm
\topmatter

\centerline{\bf Galois representations arising from twenty-seven
                lines on a cubic surface}
\centerline{\bf and the arithmetic associated with Hessian
                polyhedra}
\author{\smc Lei Yang}\endauthor
\endtopmatter
\document

\centerline{\bf Contents}
$$\aligned
 &\text{1. Introduction}\\
 &\text{2. Hessian polyhedra and cubic forms associated to $G_{25, 920}$}\\
 &\text{3. Hessian polyhedra and Picard modular forms}\\
 &\text{4. Hessian polyhedra and Galois representations associated with cubic surfaces}\\
 &\text{5. Hessian polyhedra and the arithmetic of rigid Calabi-Yau threefolds}
\endaligned$$

\vskip 0.5 cm

\centerline{\bf 1. Introduction}

\vskip 0.5 cm

  In his celebrated lecture: Mathematical problems (see \cite{H}),
which was delivered before the International Congress of
Mathematicians at Paris in 1900, David Hilbert gave his famous 23
problems. The 12th problem is the extension of Kronecker's theorem
on abelian fields to any algebraic realm of rationality. Hilbert
said:``The extension of Kronecker's theorem to the case that, in
place of the realm of rational numbers or of the imaginary
quadratic field, any algebraic field whatever is laid down as
realm of rationality, seems to me of the greatest importance. I
regard this problem as one of the most profound and far-reaching
in the theory of numbers and of functions. It will be seen that in
the problem the three fundamental branches of mathematics, number
theory, algebra and function theory, come into closest touch with
one another, and I am certain that the theory of analytical
functions of several variables in particular would be notably
enriched if one should succeed in finding and discussing those
functions which play the part for any algebraic number field
corresponding to that of the exponential function in the field of
rational numbers and of the elliptic modular functions in the
imaginary quadratic number field.''

  In his paper \cite{Wei}, Weil said:`` First of all, it will be
necessary to extend the theory of abelian functions to non-abelian
extensions of fields of algebraic functions. As I hope to show,
such an extension is indeed possible. In any case, we face here a
series of important and difficult problems, whose solution will
perhaps require the efforts of more than one generation.'' Later
on, Weil \cite{Wei3} pointed out:`` We take up the non-commutative
case.''

  According to Tate (\cite{Ta}, see also \cite{L1}), the biggest problem after Artin's
reciprocity law was to extend it in some way to non-abelian
extensions, or, from an analytic point of view, to identify his
non-abelian $L$-series. The breakthrough came out of Langlands's
discovery in the 1960's of the conjectural principle of
functoriality, which included a formulation of a nonabelian
reciprocity law as a special case. It turns out that in the
nonabelian case, one should describe the representations of
$\text{Gal}(\overline{K}/K)$ rather than
$\text{Gal}(\overline{K}/K)$ itself. Thus, the global reciprocity
law is formulated as a general conjectural correspondence between
Galois representations and automorphic forms.

  The Langlands program describes profound relations between motives
(arithmetic data) and automorphic representations (analytic data).
Wiles's spectacular work on the Shimura-Taniyama-Weil conjecture,
which established the proof of Fermat's Last Theorem, can be
regarded as confirmation of such a relation in the case of
elliptic curves. In general, the arithmetic information wrapped up
in motives comes from solutions of polynomial equations with
rational coefficients. The analytic information from automorphic
representations is backed up by the structure of Lie theory. Thus,
The Langlands program translates large parts of algebraic geometry
into the language of representation theory of Lie groups. The
Langlands correspondence embodies a large part of number theory,
arithmetic algebraic geometry and representation theory of Lie
groups. The Artin conjecture on $L$-functions and the
Ramanujan-Petersson conjecture would follow from the Langlands
correspondence. Over number fields, the Langlands correspondence
in its full generality seems still to be out of reach. Even its
precise formulation is very involved.

  The investigation of the relation between the arithmetic of Galois
representations and the analytic behavior of corresponding
$L$-series is one of the central topics in number theory and
arithmetic geometry. One of the outstanding problems in this field
is Artin's conjecture which predicts the holomorphy of the Artin
$L$-series of nontrivial irreducible complex representations of
the absolute Galois group of number fields. The idea of describing
extensions of ${\Bbb Q}$ via the action of the Galois group
$\text{Gal}(\overline{{\Bbb Q}}/{\Bbb Q})$ on certain groups and
other algebraic objects is very fruitful. Many examples of
constructions of abelian and non-abelian extensions of ${\Bbb Q}$
are based on this idea. A complete classification of all these
extensions in terms of Galois representations and in terms of
certain objects of analysis and algebraic geometry (automorphic
forms and motives) is the central problem of the Langlands
programme (see \cite{MP}).

  In fact, elliptic curves can be considered as a special kind
of cubic hypersurfaces. A cubic hypersurface is a projective
algebraic variety defined by a homogeneous equation $F_3(x_0,
\cdots, x_n)=0$ of degree three with coefficients in some ground
field $k$. When $n=2$, it gives a cubic curve, i.e., an elliptic
curve. Now, we will study the case of $n=3$, cubic surfaces.

  Over an algebraically closed field $k$, every irreducible cubic
surface is a rational surface. The class of a hyperplane section
$h$ of a surface $F$ is the canonical class $-K_F$. Any smooth
cubic surface can be obtained from the projective plane ${\Bbb
P}^2$ by blowing up of $6$ generic points. The appropriate
birational mapping $\phi: {\Bbb P}^2 \to F$ is determined by the
linear system of cubic curves passing through the $6$ points.
There are $27$ straight lines on $F$, each of which is
exceptional; they are the only exceptional curves on $F$. The
configuration of these $27$ lines is rich in symmetries. The
automorphism group of the corresponding graph is isomorphic to the
Weyl group of $E_6$. Cubic surfaces belong to the class of del
Pezzo surfaces, i.e. two-dimensional Fano varieties. Over an
algebraically non-closed field $k$, there are smooth cubic
surfaces $F$ which are not birationally isomorphic to ${\Bbb P}^2$
over $k$. Among these surfaces one finds surfaces possessing
$k$-points, and these are unirational over $k$. The group of
birational automorphisms of a minimal surface has been determined
and an arithmetic theory of cubic surfaces has been developed in
\cite{Ma} by the non-associative structures, such as quasi-groups
and Moufang loops.

  In his paper \cite{Wei2}, Weil showed that for nonsingular
rational surfaces over an algebraic number field $K$, the
Hasse-Weil $L$-function, except for a finite number of factors, is
the same as a suitable Artin $L$-function for a certain extension
of $K$. For instance, for a nonsingular cubic surface in the
projective $3$-space, one thus gets an $L$-function belonging to
the extension of $K$ determined by the $27$ straight lines on the
surface. The Galois group for this is $W(E_6)$, a group of order
$2^7 \cdot 3^4 \cdot 5$ and has a simple subgroup of index $2$.
Therefore, the Hasse-Weil $L$-function for a given cubic surface
is essentially an Artin $L$-function of a definitely non-abelian
type for one dimensional representations, i.e., characters. Here
is a rather unexpected connection between number theory and
algebraic geometry.

  According to \cite{Ma}, let $k$ be a perfect field and let $V$
run through the smooth cubic surfaces over $k$. we have the
realization problems.

\roster

\item Which subgroups of $W(E_6)$ are realized as an image of the
      Galois group by its representation on $\text{Pic}(V \otimes
      \overline{k})$ for some $V$?

\item Which extensions of the field $k$ correspond with kernels of
      such representations?

\endroster

  In his monograph \cite{Ma}, Manin said: ``Of course, the answer
must depend heavily on $k$. The derived group $W^{\prime}(E_6)$ is
simple, and the whole of it obviously can not be realized over a
local base field which has only solvable extensions. On the other
hand, $W(E_6)$ can be realized over some function field. The more
interesting case of a number field is completely unknown. The
search for a reasonable approach to problem (2) may involve a
pattern for a non-commutative class field theory. Here, it appears
that topological-analytic considerations and something like Hecke
operators do not suffice.''

  In the present paper, we will show that three apparently disjoint
objects: Galois representations arising from twenty-seven lines on
a cubic surface (number theory and arithmetic algebraic geometry),
Picard modular forms (automorphic forms), rigid Calabi-Yau
threefolds and their arithmetic (Diophantine geometry) are
intimately related to Hessian polyhedra and their invariants. We
investigate the Artin $L$-functions for six dimensional
representations associated to a kind of cubic surfaces. We study
the linear representations of the Galois group
$\text{Gal}(\overline{{\Bbb Q}}/{\Bbb Q})$ which constitutes a
primary object of number theory and find the correspondence
between some transcendental functions and representations.
Moreover, there exists the duality which connects them with each
other. Our linear representations come from the Picard groups of
cubic surfaces. We construct a Galois representation whose image
is a proper subgroup of $W(E_6)$, the Weyl group of the
exceptional Lie algebra $E_6$. We give a conjecture about the
identification of two different kinds of $L$-functions which can
be considered as a higher dimensional counterpart of the
Langlands-Tunnell theorem. We obtain the structure of rational
points in the rigid Calabi-Yau threefolds. More precisely, we will
investigate the relations among invariant theory (algebra), Picard
modular forms (analysis), Galois representations (arithmetic) and
cubic surfaces (geometry and topology). Just as the applications
of modular curves and elliptic modular functions to number theory,
i.e., one uses $GL(2)$ to study $GL(1)$, we will apply Picard
modular surfaces and Picard modular functions to number theory,
i.e., we will use $GL(3)$ (or $U(3)$) to study $GL(1)$. As an
application, we find an answer to Manin's problems described as
above. Moreover, we give a conjecture which involves the
noncommutative class field theory. It should be pointed out that
our Galois representations arise from twenty-seven lines on a
cubic surface, while the previous Galois representations either
arise from torsion points of an elliptic curve or an abelian
variety (e.g. the Tate modules), or arise from Mordell-Weil
lattices of an elliptic curve, i.e., both of them come from
elliptic curves (cubic curves).

\vskip 0.5 cm

{\bf Main Results}.

\vskip 0.5 cm

  We have the following dictionaries.

\roster
\item
$$\aligned
 &GL(2)\\
 &\text{elliptic curves (cubic curves)}\\
 &\text{elliptic modular forms}\\
 &\text{regular polyhedra}\\
 &\text{algebraic solutions of hypergeometric differential equations}\\
 &A_4, S_4, A_5\\
 &PSL(2, {\Cal O}_K)/PSL(2, {\Cal O}_K)(\sqrt{5})=A_5\\
 &\text{Hilbert modular surfaces over $K={\Bbb Q}(\sqrt{5})$ (Hirzebruch)}\\
 &\text{diagonal cubic surface of Clebsch and Klein, $F_1$, $27$ real lines}\\
 &\text{two dimensional real surface}
\endaligned$$
\item
$$\aligned
 &GL(3)\\
 &\text{cubic surfaces}\\
 &\text{Picard modular forms}\\
 &\text{Hessian polyhedra}\\
 &\text{algebraic solutions of Appell hypergeometric partial differential equations}\\
 &\text{Hessian groups}: G_{648}, G_{1296}\\
 &U(2, 1; {\Cal O}_K)/U(2, 1; {\Cal O}_K)(1-\omega)=S_4 \text{(see \cite{Ho1}, \cite{Ho3})}\\
 &\text{Picard modular surfaces over $K={\Bbb Q}(\sqrt{-3})$}\\
 &\text{our cubic surface, $F_2$, $15$ real lines and $12$ complex lines}\\
 &\text{two dimensional complex surface}
\endaligned$$
\endroster

  Let us consider the following Galois representation:
$$\rho: \text{Gal}(\overline{K}/K) \to \text{Aut}(X),$$
where $X$ is some geometric object.

\roster

\item $X=V_n(k)$ is an $n$-dimensional linear space over a field
      $k$, such as ${\Bbb C}$, ${\Bbb R}$, ${\Bbb Q}_p$ and ${\Bbb F}_p$,
      $\text{Aut}(V_n(k))=GL(n, k)$(Lie groups, $p$-adic Lie groups, finite groups)

\item $X=T_{\ell}(\mu)$ is the Tate module of the number field
      $K$, $\text{Aut}(T_{\ell}(\mu))=GL(1, {\Bbb Z}_{\ell})$($\ell$-adic Lie groups)

\item $X=T_{\ell}(E)$ is the Tate module of the elliptic curve
      $E$, $\text{Aut}(T_{\ell}(E))=GL(2, {\Bbb Z}_{\ell})$($\ell$-adic Lie groups)

\item $X=(E(k(C)), \langle , \rangle)$ is the Mordell-Weil lattice
      associated to a rational elliptic surface, where $k(C)$ is the function field
      of a smooth projective curve $C$ over an algebraically closed field $k$,
      $\text{Aut}(E(k(C)), \langle , \rangle)=\text{Aut}(E_r)$, $r=6, 7, 8$(finite groups)

\item $X={\Cal L}$ is the configuration of $27$ lines on the cubic
      surface $S$ preserving the intersection behavior of the lines, $\text{Aut}({\Cal L})=W(E_6)$(finite groups)

\endroster

  From the point of view of finite simple groups, the significance
of $W(E_6)$ comes from the following three pairs $(A_5, S_5)$,
$(A_6, S_6)$, $(O_5(3), W(E_6))$ and Brauer's theorem (see
\cite{Br}): If $G$ is a simple group of an order $g=5 \cdot 3^a
\cdot 2^b>5$, then $G$ is isomorphic to one of the following three
groups: the alternating group $A_5$ of order $60$, the alternating
group $A_6$ of order $360$ or the orthogonal group $O_5(3)$ of
order $25920$. The first two pairs appeared in Klein's celebrated
book \cite{Kl}. In this paper, we will study the third one.

  One of the significance of the group $W(E_6)$ is that it possesses in
terms of different Chevalley groups (see \cite{Hu}, \cite{CC}):
$$G_{25, 920} \cong PU(3, 1; {\Bbb F}_4) \cong PSp(4, {\Bbb Z}/3 {\Bbb Z}),$$
where $G_{25, 920} \subset W(E_6)$ is the simple subgroup of index
two consisting of all even elements. We have the action of the
unitary reflection group of order $25, 920$: $G_{25, 920}$ on
${\Bbb P}^4=\{ (Y_0, Y_1, Y_2, Y_3, Y_4) \}$. The following
theorems say that the invariants of $G_{25, 920}$ (see section two
for the details):
$$\Phi, u, t, \Psi_1, C_6, C_9, \Phi_3, t_3, C_{12}, C_{18}$$
can be expressed as the polynomials of the invariants $f_0$,
$f_1$, $G$, $H$, $K$ and $C_6$, where
$$\Phi=Y_0^4+8 Y_0 Y_1^3,$$
$$u=Y_0 \psi+6 Y_1 \varphi,$$
$$t=Y_0^6-20 Y_0^3 Y_1^3-8 Y_1^6,$$
$$\Psi_1=\psi(-Y_0^3+4 Y_1^3)+18 \varphi Y_0^2 Y_1,$$
$$C_6=Y_2^6+Y_3^6+Y_4^6-10 (Y_2^3 Y_3^3+Y_3^3 Y_4^3+Y_4^3 Y_2^3),$$
$$C_9=(Y_2^3-Y_3^3)(Y_3^3-Y_4^3)(Y_4^3-Y_2^3),$$
$$\Phi_3=-\psi^3 Y_0+18 \varphi \psi^2 Y_1+108 \varphi^3 Y_0,$$
$$t_3=\psi^3(Y_0^3+8 Y_1^3)-54 \varphi \psi^2 Y_0^2 Y_1+324
      \varphi^2 \psi Y_0 Y_1^2+216 \varphi^3 (Y_0^3-Y_1^3),$$
$$C_{12}=(Y_2^3+Y_3^3+Y_4^3)[(Y_2^3+Y_3^3+Y_4^3)^3+216 Y_2^3 Y_3^3 Y_4^3],$$
$$C_{18}=(Y_2^3+Y_3^3+Y_4^3)^6-540 Y_2^3 Y_3^3 Y_4^3 (Y_2^3+Y_3^3+Y_4^3)^3-5832 Y_2^6 Y_3^6 Y_4^6.$$
The other three invariants are given by
$$\Psi=Y_0^3 Y_1-Y_1^4,$$
$$\Psi_2=-\psi^2 Y_1^2-3 \varphi \psi Y_0^2+18 \varphi^2 Y_0 Y_1,$$
$${\frak C}_{12}=Y_2 Y_3 Y_4 [27 Y_2^3 Y_3^3 Y_4^3-(Y_2^3+Y_3^3+Y_4^3)^3].$$
Here, the invariants $f_0$, $f_1$, $G$, $H$, $K$ are given by
$$f_0=Y_0+2 Y_1, \quad f_1=Y_0-Y_1.\tag 1.1$$
$$G=(Y_2-Y_3)(Y_3-Y_4)(Y_4-Y_2), \quad H=\psi+6 \varphi, \quad K=\psi-3 \varphi\tag 1.2$$
with
$$\varphi=Y_2 Y_3 Y_4, \quad \psi=Y_2^3+Y_3^3+Y_4^3.\tag 1.3$$

{\smc Theorem 1.1} (see \cite{Y3}, Main Theorem 4). {\it The
invariants $G$, $H$ and $K$ satisfy the following algebraic
equations, which are the form-theoretic resolvents $($algebraic
resolvents$)$ of $G$, $H$, $K$$:$
$$\left\{\aligned
  4 G^3+H^2 G-C_6 G-4 C_9 &=0,\\
  H (H^3+8 K^3)-9 C_{12} &=0,\\
  K (K^3-H^3)-27 {\frak C}_{12} &=0.
\endaligned\right.\tag 1.4$$}
Consequently,
$$C_{18}=-\frac{1}{27}(H^6-20 H^3 K^3-8 K^6).\tag 1.5$$

{\smc Theorem 1.2} (Main Theorem 1). {\it The invariants $f_0$,
$f_1$, $H$ and $K$ satisfy the following algebraic equations,
which are the form-theoretic resolvents $($algebraic resolvents$)$
of $f_0$, $f_1$, $H$ and $K$$:$
$$u=\frac{1}{3}(f_0 H+2 f_1 K).\tag 1.6$$
$$\left\{\aligned
  \Phi &=\frac{1}{9} f_0 (f_0^3+8 f_1^3),\\
  \Psi &=\frac{1}{9} f_1 (f_0^3-f_1^3),\\
     t &=-\frac{1}{27} (f_0^6-20 f_0^3 f_1^3-8 f_1^6).
  \endaligned\right.\tag 1.7$$
$$\left\{\aligned
  \Psi_1 &=\frac{1}{9} [(f_0^3-4 f_1^3) H-6 f_0^2 f_1 K],\\
  \Psi_2 &=\frac{1}{9} (-f_0^2 H K+2 f_0 f_1 K^2-f_1^2 H^2),\\
  \Phi_3 &=\frac{1}{9} [f_0 (H^3-4 K^3)-6 f_1 H^2 K],\\
     t_3 &=\frac{1}{27} [f_0^3 (H^3+8 K^3)-18 f_0^2 f_1 H^2 K+
           36 f_0 f_1^2 H K^2+8 f_1^3 (H^3-K^3)]\\
         &=\frac{1}{27} [(f_0^3+8 f_1^3) H^3-18 f_0^2 f_1 H^2 K+
           36 f_0 f_1^2 H K^2+8 (f_0^3-f_1^3) K^3].
\endaligned\right.\tag 1.8$$}

  We find the following nonlinear duality:
$$f_0 \longleftrightarrow H, \quad f_1 \longleftrightarrow K.\tag 1.9$$
Here, $f_0$ and $f_1$ are linear, $H$ and $K$ are cubic. The
functions $u$ and $t_3$ are invariant under the above nonlinear
dual transformation. Under the nonlinear dual transformation
(1.9),
$$C_{12} \longleftrightarrow \Phi, \quad
  {\frak C}_{12} \longleftrightarrow -\frac{1}{3} \Psi, \quad
  C_{18} \longleftrightarrow t, \quad
  \Psi_1 \longleftrightarrow \Phi_3.\tag 1.10$$

  We will study the analogies between Hessian invariants and
Picard modular forms, which can be considered as the higher
dimensional counterpart of the analogies between regular
polyhedral invariants and elliptic modular forms studied by Klein.

{\smc Conjecture 1.3}. {\it
$${\Bbb C}[z_1, z_2, z_3]^{\text{Hessian groups}}
 ={\Bbb C}[C_6, C_9, C_{12}].\tag 1.11$$
The map
$$\phi: {\Bbb C}[C_6, C_9, C_{12}] \to {\Bbb C}[G_2, G_3, G_4]\tag 1.12$$
given by
$$\phi(h(C_6, C_9, C_{12}))=h(G_2, G_3, G_4)\tag 1.13$$
is an algebra isomorphism.}

  Here,
$$\left\{\aligned
  C_6(z_1, z_2, z_3) &=z_1^6+z_2^6+z_3^6-10
                     (z_1^3 z_2^3+z_2^3 z_3^3+z_3^3 z_1^3),\\
  C_9(z_1, z_2, z_3) &=(z_1^3-z_2^3)(z_2^3-z_3^3)(z_3^3-z_1^3),\\
  C_{12}(z_1, z_2, z_3) &=(z_1^3+z_2^3+z_3^3)[(z_1^3+z_2^3+z_3^3)^3
                        +216 z_1^3 z_2^3 z_3^3].
\endaligned\right.$$
$G_2$, $G_3$, $G_4$ are of weight $2$, $3$ and $4$, respectively,
which were introduced by Holzapfel in \cite{Ho1}.

  Put
$$X(z_1, z_2, z_3)=z_1^3, \quad Y(z_1, z_2, z_3)=z_2^3,
  \quad Z(z_1, z_2, z_3)=z_3^3.$$
$$\varphi=z_1 z_2 z_3, \quad Q_1(z_1, z_2, z_3)=z_1 z_2^2+z_2 z_3^2+z_3 z_1^2, \quad
  Q_2(z_1, z_2, z_3)=z_1^2 z_2+z_2^2 z_3+z_3^2 z_1.$$

{\smc Theorem 1.4} (Main Theorem 2). {\it Let
$$\left\{\aligned
  W_2 &=(X+Y+Z)^2-12(XY+YZ+ZX),\\
  W_3 &=(X-Y)(Y-Z)(Z-X),\\
  {\frak W}_3 &=XYZ-\varphi^3,\\
  W_4 &=(X+Y+Z)[(X+Y+Z)^3+216 XYZ],\\
  {\frak W}_4 &=\varphi [27 XYZ-(X+Y+Z)^3],\\
  W_6 &=(X+Y+Z)^6-540 XYZ (X+Y+Z)^3-5832 X^2 Y^2 Z^2,\\
  {\frak V}_3 &=Q_1^3+Q_2^3-(X+Y+Z+6 \varphi)(XY+YZ+ZX)-6 \varphi^2 (X+Y+Z)-9 \varphi^3,\\
  {\frak V}_2 &=Q_1 Q_2-(XY+YZ+ZX)-\varphi(X+Y+Z)-3 \varphi^2,\\
  {\frak U}_3 &=Q_1^3-Q_2^3-(X-Y)(Y-Z)(Z-X).
\endaligned\right.\tag 1.14$$
They satisfy the following four relations:
$$\left\{\aligned
  W_2^3-3 W_2 W_4+2 W_6 &=432 W_3^2,\\
  8 U^3 (W_6^2-W_4^3-1728 {\frak W}_4^3) &=27 {\frak W}_3 (W_4-9 U^4)^3,
\endaligned\right.\tag 1.15$$
where
$$27 U^8-18 W_4 U^4-8 W_6 U^2-W_4^2=0,\tag 1.16$$
and the Jacobian:
$$\frac{\partial(W_2, W_3, W_4, {\frak W}_3)}{\partial(X, Y, Z, \varphi)}
 =288 {\frak W}_4^2.\tag 1.17$$
$$\frac{\partial({\frak V}_3, {\frak V}_2)}{\partial(Q_1, Q_2)}=3({\frak U}_3+W_3).\tag 1.18$$}

 We find the following correspondence:
$$(G_2, G_3, G_4, z^2) \longleftrightarrow (W_2, W_3, W_4, {\frak W}_4).$$
Here,
$$\Delta^2=16 G_2^4 G_4-128 G_2^2 G_4^2-4 G_2^3 G_3^2+144 G_2 G_3^2 G_4-27 G_3^4+256 G_4^3, \quad z^6=\Delta^2.$$
$$6912 {\frak W}_4^3=W_2^6+9 W_2^2 W_4^2+432^2 W_3^4-4 W_4^3
  -6 W_2^4 W_4-864 W_2^3 W_3^2+2592 W_2 W_3^2 W_4.$$

{\smc Theorem 1.5} (Main Theorem 3). {\it Let
$$x=Q_1-Q_2, \quad y=X+Y+Z+6 \varphi, \quad z=X+Y+Z-3 \varphi.\tag 1.19$$
Then $x$, $y$ and $z$ satisfy the following three equations:
$$\aligned
 &2^4 \cdot 3^8 {\frak W}_3^2 (4 {\frak U}_3-4 x^3-x y^2+W_2 x+4 W_3)^3\\
=&2^{10} \cdot 3^{11} {\frak W}_3^2 {\frak V}_2^3 ({\frak
U}_3+W_3)
  -(4{\frak U}_3-4 x^3-x y^2+W_2 x+4 W_3) {\frak V}_2^2 \times\\
 &\times \{2^8 \cdot 3^{11} {\frak W}_3^2 {\frak V}_2+[(9 W_4-y
 (y^ 3+8 z^3))^2-2^6 \cdot 3^{10} W_2 {\frak W}_3^2]+\\
 &+2^5 [27 {\frak W}_4-z (z^3-y^3)][9 W_4-y (y^3+8 z^3)]+2^8
  [27 {\frak W}_4-z (z^3-y^3)]^2\}.
\endaligned$$
$$\aligned
 &[9W_4-y(y^3+8z^3)]^4+2^{12} [9W_4-y(y^3+8z^3)] \{3^{18} {\frak W}_3^4
  +[27 {\frak W}_4-z(z^3-y^3)]^3\}\\
=&2^{12} \cdot 3^{20} {\frak W}_3^4 W_4.
\endaligned$$
$$\aligned
 &[9W_4-y(y^3+8z^3)]^6-2^{11} \cdot 5 [9W_4-y(y^3+8z^3)]^3 \{3^{18} {\frak W}_3^4
  +[27 {\frak W}_4-z(z^3-y^3)]^3\}\\
 &-2^{21} \{3^{18} {\frak W}_3^4+[27 {\frak W}_4-z(z^3-y^3)]^3\}^2
  =2^{18} \cdot 3^{30} {\frak W}_3^6 W_6.
\endaligned$$}

  We have the following table which describes the relation between
Hessian polyhedra, invariant theory and Picard modular forms.
$$\text{Hessian polyhedra} \longleftrightarrow \text{Picard modular functions}$$
$$(C_6, C_9, C_{12}) \longleftrightarrow (G_2, G_3, G_4), \quad \text{weight $2$, $3$, $4$}$$
$$(G, H, K) \longleftrightarrow (\xi_1, \xi_2, \xi_3), \quad \text{weight $1$}$$
$$(\varphi, X, Y, Z, Q_1, Q_2) \longleftrightarrow (\zeta_1, \zeta_2, \zeta_3, \zeta_4, \zeta_5, \zeta_6), \quad \text{weight $1$}$$

  It is known that the Hessian group of order $216$ is generated
by the five generators (see \cite{Y3}):
$$A=\left(\matrix
  0 & 1 & 0\\
  0 & 0 & 1\\
  1 & 0 & 0
  \endmatrix\right), \quad
  B=\left(\matrix
  1 & 0 & 0\\
  0 & 0 & 1\\
  0 & 1 & 0
  \endmatrix\right),$$
$$C=\left(\matrix
  1 & 0 & 0\\
  0 & \omega & 0\\
  0 & 0 & \omega^2
  \endmatrix\right), \quad
  D=\left(\matrix
  1 & 0 & 0\\
  0 & \omega & 0\\
  0 & 0 & \omega
  \endmatrix\right),$$
$$E=\frac{1}{\sqrt{-3}} \left(\matrix
  1 & 1 & 1\\
  1 & \omega & \omega^2\\
  1 & \omega^2 & \omega
  \endmatrix\right).$$
The actions of $A$, $B$, $C$ and $E$ on $(X, Y, Z, \varphi, Q_1,
Q_2)$ are given as follows:
$$\aligned
 &A(X, Y, Z, \varphi, Q_1, Q_2)\\
=&\left(\matrix
   & 1 &   &   &   &   \\
   &   & 1 &   &   &   \\
 1 &   &   &   &   &   \\
   &   &   & 1 &   &   \\
   &   &   &   & 1 &   \\
   &   &   &   &   & 1
 \endmatrix\right)
 \left(\matrix
   X\\
   Y\\
   Z\\
   \varphi\\
   Q_1\\
   Q_2
 \endmatrix\right),
\endaligned\tag 1.20$$
$$\aligned
 &B(X, Y, Z, \varphi, Q_1, Q_2)\\
=&\left(\matrix
 1 &   &   &   &   &   \\
   &   & 1 &   &   &   \\
   & 1 &   &   &   &   \\
   &   &   & 1 &   &   \\
   &   &   &   &   & 1 \\
   &   &   &   & 1 &
 \endmatrix\right)
 \left(\matrix
   X\\
   Y\\
   Z\\
   \varphi\\
   Q_1\\
   Q_2
 \endmatrix\right),
\endaligned\tag 1.21$$
$$\aligned
 &C(X, Y, Z, \varphi, Q_1, Q_2)\\
=&\left(\matrix
 1 &   &   &   &   &   \\
   & 1 &   &   &   &   \\
   &   & 1 &   &   &   \\
   &   &   & 1 &   &   \\
   &   &   &   & \overline{\omega} &   \\
   &   &   &   &   & \omega
 \endmatrix\right)
 \left(\matrix
   X\\
   Y\\
   Z\\
   \varphi\\
   Q_1\\
   Q_2
 \endmatrix\right),
\endaligned\tag 1.22$$
$$\aligned
 &(\sqrt{-3})^3 E(X, Y, Z, \varphi, Q_1, Q_2)\\
=&\left(\matrix
 1 & 1 & 1 & 6 & 3 & 3 \\
 1 & 1 & 1 & 6 & 3 \overline{\omega} & 3 \omega\\
 1 & 1 & 1 & 6 & 3 \omega & 3 \overline{\omega}\\
 1 & 1 & 1 & -3 & 0 & 0\\
 3 & 3 \overline{\omega} & 3 \omega & 0 & 0 & 0\\
 3 & 3 \omega & 3 \overline{\omega} & 0 & 0 & 0
 \endmatrix\right)
 \left(\matrix
   X\\
   Y\\
   Z\\
   \varphi\\
   Q_1\\
   Q_2
 \endmatrix\right),
\endaligned\tag 1.23$$

{\smc Proposition 1.6}. {\it The space curves $W_2={\frak
W}_3={\frak V}_3={\frak V}_2=0$, $W_3={\frak W}_3={\frak
V}_3={\frak V}_2=0$, $W_4={\frak W}_3={\frak V}_3={\frak V}_2=0$,
$W_6={\frak W}_3={\frak V}_3={\frak V}_2=0$ and ${\frak
W}_4={\frak W}_3={\frak V}_3={\frak V}_2=0$ are invariant curves
on the invariant surface ${\frak W}_3={\frak V}_3={\frak V}_2=0$
under the action of the subgroup of Hessian group generated by
$A$, $B$, $C$ and $E$.}

  Put
$$g_1=E(x), \quad g_2=EC(x), \quad g_3=EC^2(x), \quad
  g_4=x, \quad g_5=C(x), \quad g_6=C^2(x).\tag 1.24$$
We find that
$$S: \left\{\aligned
  g_1+g_2+g_3 &=0,\\
  g_4+g_5+g_6 &=0,\\
  g_1^3+g_2^3+g_3^3 &=g_4^3+g_5^3+g_6^3,
\endaligned\right.\tag 1.25$$
which is a cubic surface.

  The study of the general cubic surface dates from 1849, in which
year the $27$ lines were found by Cayley and Salmon; the discovery
of the Sylvester pentahedron follwed two years later, and the
theory attained a remarkable degree of elegance and completeness
with the introduction of the plane representation, made
independently by Clebsch and Cremona in 1866 (see \cite{Seg}).
Neverless, even today the study of what Cayley has called `the
complicated and many-sided symmetry' of the relations between the
$27$ lines can not be said to have been exhausted, especially
Schl\"{a}fli's well-known notation for the lines readily lend
itself to the description and deeper investigation of that
symmetry.

  As far as the equation defining a cubic surface is concerned, we
have the following result, claimed by Sylvester (1851) and
Steiner, and proved by Clebsch (1861). It is the so-called
pentahedral form of the cubic.

{\smc Theorem 1.7} (see \cite{Hu}). {\it A general cubic form
$F(x_0: x_1: x_2: x_3)$ can be put, in a unique way, in the form
$$a_1 y_1^3+a_2 y_2^3+a_3 y_3^3+a_4 y_4^3+a_5 y_5^3=0,$$
where the coordinates $y_i$ are linear in $(x_0: x_1: x_2: x_3)$
and satisfy
$$y_1+y_2+y_3+y_4+y_5=0.$$}

  A second special form of the equation, the hexahedral form, is
given by the polar hexagons:

{\smc Theorem 1.8} (see \cite{Hu}). {\it A general cubic form $F$
can be put in a four-dimensional family of forms:
$$z_1^3+z_2^3+z_3^3+z_4^3+z_5^3+z_6^3=0,$$
where the $z_i$ are linear in the $x_i$.}

  In fact, the pentahedral form of the cubic is closed related to
the Clebsch-Klein diagonal cubic surface. In 1873 Klein proved
that the famous diagonal surface of Clebsch, which is the surface
in ${\Bbb P}^4({\Bbb C})$ with equations
$$x_0^3+x_1^3+x_2^3+x_3^3+x_4^3=0, \quad x_0+x_1+x_2+x_3+x_4=0,$$
can be obtained from ${\Bbb P}^2({\Bbb C})$ by blowing up $6$
points in ${\Bbb P}^2({\Bbb C})$ in a special position, namely the
$6$ points in ${\Bbb P}^2({\Bbb R})=S^2/\{ \pm 1 \}$ corresponding
to the $12$ vertices of an icosahedron inscribed in $S^2$.
Hirzebruch blowed up $10$ more points, namely those corresponding
to the $20$ vertices of the dual dodecahedron. The resulting
surface $Y$ can also be obtained from the Clebsch diagonal surface
by blowing up $10$ Eckhardt points, that is points, where $3$ of
the $27$ lines on the surface meet. Let ${\Cal O}_K \subset {\Bbb
Q}(\sqrt{5})$ be the ring of integers and $\Gamma \subset SL(2,
{\Cal O}_K)$ the congruence subgroup mod $2$. Hirzebruch proved
that the icosahedral surface $Y$ is the minimal resolution of
$\overline{{\Bbb H}^2/\Gamma}$.

  The significance of our cubic surface is that it is the
hypersurface of the hexahedral form of the cubic.

  There are $27$ lines on the cubic surface $S$:
$$l_1: \quad g_1=g_2+g_3=g_4=g_5+g_6=0,$$
$$l_2: \quad g_1=g_2+g_3=g_5=g_4+g_6=0,$$
$$l_3: \quad g_1=g_2+g_3=g_6=g_4+g_5=0,$$
$$l_4: \quad g_2=g_1+g_3=g_4=g_5+g_6=0,$$
$$l_5: \quad g_2=g_1+g_3=g_5=g_4+g_6=0,$$
$$l_6: \quad g_2=g_1+g_3=g_6=g_4+g_5=0,$$
$$l_7: \quad g_3=g_1+g_2=g_4=g_5+g_6=0,$$
$$l_8: \quad g_3=g_1+g_2=g_5=g_4+g_6=0,$$
$$l_9: \quad g_3=g_1+g_2=g_6=g_4+g_5=0,$$
and
$$l_{j, 1}: \quad g_1=\omega^j g_4, \quad g_2=\omega^j g_5, \quad g_3=\omega^j g_6,
            \quad g_1+g_2+g_3=0,$$
$$l_{j, 2}: \quad g_1=\omega^j g_4, \quad g_2=\omega^j g_6, \quad g_3=\omega^j g_5,
            \quad g_1+g_2+g_3=0,$$
$$l_{j, 3}: \quad g_1=\omega^j g_5, \quad g_2=\omega^j g_4, \quad g_3=\omega^j g_6,
            \quad g_1+g_2+g_3=0,$$
$$l_{j, 4}: \quad g_1=\omega^j g_5, \quad g_2=\omega^j g_6, \quad g_3=\omega^j g_4,
            \quad g_1+g_2+g_3=0,$$
$$l_{j, 5}: \quad g_1=\omega^j g_6, \quad g_2=\omega^j g_4, \quad g_3=\omega^j g_5,
            \quad g_1+g_2+g_3=0,$$
$$l_{j, 6}: \quad g_1=\omega^j g_6, \quad g_2=\omega^j g_5, \quad g_3=\omega^j g_4,
            \quad g_1+g_2+g_3=0,$$
where $j \equiv 0, 1, 2$(mod $3$).

  Let us consider the configuration of the $27$ lines (forgetting
the cubic surface). That means we consider simply the set of $27$
lines together with the incidence relations they satisfy. There is
a group of automorphisms of the configuration meaning the
permutations of the set of $27$ lines, preserving the incidence
relations.

  There are three difficulties: 1. The realization of $27$ lines
on the cubic surface by algebraic equations. 2. Under the action
of the subgroup of $W(E_6)$, $27$ lines on the cubic surface may
not lie on the same cubic surface. 3. Under the action of the
subgroup of $W(E_6)$, the incidence relations of $27$ lines may
not preserve.

  The double six is given by
$$N=\left(\matrix
    a_1 & a_2 & a_3 & a_4 & a_5 & a_6\\
    b_1 & b_2 & b_3 & b_4 & b_5 & b_6
     \endmatrix\right)
  =\left(\matrix
    l_{1, 1} & l_{2, 2} & l_{2, 3} & l_{1, 4} & l_{1, 5} & l_{2, 6}\\
    l_{2, 1} & l_{1, 2} & l_{1, 3} & l_{2, 4} & l_{2, 5} & l_{1, 6}
    \endmatrix\right).\tag 1.26$$

  We find that
$$c_{12}=l_1, \quad c_{13}=l_9, \quad c_{14}=l_{0, 5}, \quad c_{15}=l_{0, 4}, \quad c_{16}=l_5.\tag 1.27$$
$$c_{23}=l_{0, 6}, \quad c_{24}=l_6, \quad c_{25}=l_8, \quad c_{26}=l_{0, 3}.\tag 1.28$$
$$c_{34}=l_2, \quad c_{35}=l_4, \quad c_{36}=l_{0, 2}.\tag 1.29$$
$$c_{45}=l_{0, 1}, \quad c_{46}=l_7, \quad c_{56}=l_3.\tag 1.30$$

  We find that the action of the group $H$ generated by $A$, $B$,
$C$ and $E$ on the $27$ lines of our cubic surface $S$ gives the
permutations of the lines which preserve the intersection behavior
of the lines.

{\smc Theorem 1.9} (Main Theorem 4). {\it The group $H$ generated
by $A$, $B$, $C$ and $E$ is a subgroup of the group of
automorphisms of the configuration of $27$ lines on the cubuc
surface: $\text{Aut}({\Cal L})=W(E_6)$.}

  This gives the Galois representation:
$$\rho: G_{{\Bbb Q}}=\text{Gal}(\overline{{\Bbb Q}}/{\Bbb Q}) \to \text{Aut}({\Cal L})=W(E_6).\tag 1.31$$
Here,
$$\text{Im}(\rho)=H=\langle A, B, C, E \rangle, \quad
  \text{ker}(\rho)=\text{Gal}(\overline{{\Bbb Q}}/{\Cal K}).\tag 1.32$$
We have
$$\text{Im}(\rho) \cong \text{Gal}(\overline{{\Bbb Q}}/{\Bbb Q})/\text{ker}(\rho)
 =\text{Gal}(\overline{{\Bbb Q}}/{\Bbb Q})/\text{Gal}(\overline{{\Bbb Q}}/{\Cal K})
 \cong \text{Gal}({\Cal K}/{\Bbb Q}).\tag 1.33$$
Hence,
$$\text{Gal}({\Cal K}/{\Bbb Q}) \cong H.\tag 1.34$$

{\smc Conjecture 1.10} (Main Conjecture). {\it Let $\rho:
\text{Gal}(\overline{{\Bbb Q}}/{\Bbb Q}) \to W(E_6)$ be a Galois
representation. When the image of $\rho$: $\text{Im}(\rho) \cong H
\leq W(E_6)$, if $\rho$ is odd, i.e., $\det \rho(\sigma_c)=-1$,
where $\sigma_c \in \text{Gal}(\overline{{\Bbb Q}}/{\Bbb Q})$ is
the complex conjugate, then there exists a Picard modular form $f$
of weight one such that
$$L(\rho, s)=L(f, s)$$
up to finitely many Euler factors, where $L(\rho, s)$ is the Artin
$L$-function associated to the Galois representation $\rho$ and
$L(f, s)$ is the automorphic $L$-function associated to the Picard
modular form $f$. If $\rho$ is even, i.e., $\det
\rho(\sigma_c)=1$, then there exists an automorphic form $f$
associated to $U(2, 1):$ $\Delta f=s(s-2)f$ with $s=1$ and $($see
\cite{Y1}$)$
$$\aligned
  \Delta &=(z_1+\overline{z_1}-z_2 \overline{z_2}) \left[(z_1+\overline{z_1})
    \frac{\partial^2}{\partial z_1 \partial \overline{z_1}}+
    \frac{\partial^2}{\partial z_2 \partial \overline{z_2}}+
    z_2 \frac{\partial^2}{\partial \overline{z_1} \partial z_2}+
    \overline{z_2} \frac{\partial^2}{\partial z_1 \partial
    \overline{z_2}}\right],\\
        &(z_1, z_2) \in {\frak S}_2=\{ (z_1, z_2) \in {\Bbb C}^2:
    z_1+\overline{z_1}-z_2 \overline{z_2}>0\},
\endaligned$$
such that
$$L(\rho, s)=L(f, s)$$
up to finitely many Euler factors, where $L(\rho, s)$ is the Artin
$L$-function associated to the Galois representation $\rho$ and
$L(f, s)$ is the automorphic $L$-function associated to the
automorphic form $f$.}

  Note that the group $W(E_6)$ is combinatorially defined, while
the subgroup $H$ is defined on algebraic variety. Thus we give a
connection between combinatoric and algebraic geometry.
Furthermore, the Galois representation $\rho:
\text{Gal}(\overline{{\Bbb Q}}/{\Bbb Q}) \to H $ comes from number
theory, the Picard modular forms come from analysis and
representation theory of Lie groups. Here, Galois symmetry and
geometric symmetry meet together. Therefore, our theorem and
conjecture gives a connection which involves combinatoric,
algebraic geometry, number theory, analysis and representation
theory. More precisely, we find that Hessian polyhedra and their
invariants play a role of bridge which connects Galois
representations arising from $27$ lines on a cubic surface
(arithmetic) with Picard modular forms (analysis). In fact, there
are two aspects about Hessian polyhedra. One is Hessian group
which is related to analysis, the other is the system of
invariants which is related to algebra. Thus, Hessian polyhedra
give applications to number theory and algebraic geometry,
especially arithmetic algebraic geometry and noncommutative class
field theory.

{\smc Theorem 1.11} (Main Theorem 5). {\it The variety associated
with Hessian polyhedra
$$E: \quad y^2=x^3-27 C_{12} x+54 C_{18}\tag 1.31$$
can be expressed as the fiber product of two isogenous,
semi-stable, rational elliptic modular surfaces
$$E_{1, t}: \quad (t-1)(z_1^3+z_2^3+z_3^3)-3(t+2) z_1 z_2 z_3=0$$
and
$$E_{2, t}: \quad y^2=x^3-3t(t^3+8)x-2(t^6-20 t^3-8).$$
Hence, $E$ is a rigid Calabi-Yau threefold. There are only
finitely many rational points in $E_{1, t}$ and $E_{2, t}$.
However, there exist infinitely many nontrivial rational points in
$E$. More precisely, $P=(3 C_6, 108 C_9) \in E$ is of infinite
order.}

  The present paper consists of five sections. In section two, we
study Hessian polyhedra and cubic forms associated to the finite
simple group $G_{25, 920}$. In section three, we study Hessian
polyhedra and Picard modular forms. In section four, we
investigate Hessian polyhedra and Galois representations arising
from $27$ lines on a cubic surface. In the last section, we study
Hessian polyhedra and the arithmetic of rigid Calabi-Yau
threefolds.

{\bf Acknowledgement}. The author thanks Professor John McKay very
much for his stimulating discussion.

\vskip 0.5 cm

\centerline{\bf 2. Hessian polyhedra and cubic forms associated
                   to $G_{25, 920}$}

\vskip 0.5 cm

  Let us recall some basic facts about the unitary reflection
groups of order $25, 920$ (see \cite{Hu}). We have actions of
$G_{25, 920}$ on ${\Bbb P}^3$ and on ${\Bbb P}^4$, both of which
in fact are generated by unitary reflections.

   Burkhardt (see \cite{Bu1} and \cite{Bu2}) determined the invariants
of $G_{25, 920}$ acting on ${\Bbb P}^4$. Let the coordinates be
given as
$$Y_0=Y_{00}, Y_1=Y_{10}, Y_2=Y_{01}, Y_3=Y_{11}, Y_4=Y_{12}$$
where the $Y_{\alpha \beta}$ are the theta functions of
$${\Bbb P}^4=\{ Y_{\alpha \beta}=\frac{1}{2}(X_{\alpha
             \beta}+X_{-\alpha-\beta}) \}.$$
Here
$$X_{\alpha \beta}=\Theta \left[\matrix
  0 & 0\\ \frac{\alpha}{3} & \frac{\beta}{3}
  \endmatrix\right](\tau, z), \quad
  \tau \in {\Bbb S}_2, z \in {\Bbb C}^2,
  \alpha, \beta \in {\Bbb Z}/3 {\Bbb Z}.$$
Using the isomorphism $G_{25920} \cong PSp(4, {\Bbb Z}/3 {\Bbb
Z})$, it suffices to use generators of the symplectic group
$PSp(4, {\Bbb Z}/3 {\Bbb Z})$ to get generators of the action of
$G_{25920}$ on ${\Bbb P}^4$. However, Burkhardt used instead the
corresponding hyperelliptic curves, and a certain Weierstrass form
for them to describe the level $3$ structure. The generators of
the group $G_{25, 920}$ are transformations $B$, $C$, $D$ and
$S_2$, which act on ${\Bbb P}^4=\{ (Y_0, Y_1, Y_2, Y_3, Y_4)\}$ as
follows:
$$\aligned
 &B(Y_0, Y_1, Y_2, Y_3, Y_4)\\
=&\frac{1}{\sqrt{-3}}(Y_0+2 Y_1, Y_0-Y_1, Y_2+Y_3+Y_4,
  Y_2+\omega Y_3+\omega^2 Y_4, Y_2+\omega^2 Y_3+\omega Y_4),
\endaligned$$
$$C(Y_0, Y_1, Y_2, Y_3, Y_4)=(Y_0, Y_1, Y_4, Y_2, Y_3),$$
$$D(Y_0, Y_1, Y_2, Y_3, Y_4)=(-Y_0, -Y_2, -Y_1, -Y_3, -Y_4),$$
$$S_2(Y_0, Y_1, Y_2, Y_3, Y_4)=(Y_0, \omega^2 Y_1, Y_2, \omega^2 Y_3, \omega^2 Y_4).$$
In the form of matrices, we have
$$B=\frac{1}{\sqrt{-3}} \left(\matrix
    1 &  2 &   &   &   \\
    1 & -1 &   &   &   \\
      &    & 1 & 1 & 1\\
      &    & 1 & \omega & \omega^2\\
      &    & 1 & \omega^2 & \omega
   \endmatrix\right), \quad
 C=\left(\matrix
    1 &   &   &   &  \\
      & 1 &   &   &  \\
      &   & 0 & 0 & 1\\
      &   & 1 & 0 & 0\\
      &   & 0 & 1 & 0
   \endmatrix\right),$$
$$D=\left(\matrix
   -1 &    &    &    &   \\
      &  0 & -1 &    &   \\
      & -1 &  0 &    &   \\
      &    &    & -1 &   \\
      &    &    &    & -1
   \endmatrix\right), \quad
  S_2=\left(\matrix
    1 &          &   &          &         \\
      & \omega^2 &   &          &         \\
      &          & 1 &          &         \\
      &          &   & \omega^2 &         \\
      &          &   &          & \omega^2
   \endmatrix\right).$$
Note that
$$B^4=C^3=D^2=S_2^3=I$$
and
$$B^2=-\left(\matrix
       1 &   &   &   &  \\
         & 1 &   &   &  \\
         &   & 1 & 0 & 0\\
         &   & 0 & 0 & 1\\
         &   & 0 & 1 & 0
      \endmatrix\right).$$
In fact, $B, C, D, S_2 \in SL(5, {\Bbb C})$.

  Consider the hyperplane ${\Cal S}$ given by ${\Cal S}=\{ Y_0=0 \}$
and the plane ${\Cal J}=\{ Y_0=Y_1=0 \}$. The stabilizer of each
is a subgroup of order $648$, but these two subgroups of order
$648$ are not conjugate to each other; indeed,
$$N({\Cal S})=\langle C, D, S_2 \rangle, \quad
  N({\Cal J})=\langle B, C, S_2 \rangle.$$
Hence, $G_{25920}$ is generated by the subgroup of order $648$
acting on ${\Cal J}$, generated by $B$, $C$ and $S_2$, and by the
centralizer of ${\Cal J}$. Burkhardt shows this to be a homogenous
tetrahedral group, and its invariants are:
$$\left\{\aligned
  \Phi &=Y_0^4+8 Y_0 Y_1^3,\\
  \Psi &=Y_0^3 Y_1-Y_1^4,\\
     t &=Y_0^6-20 Y_0^3 Y_1^3-8 Y_1^6,
\endaligned\right.\tag 2.1$$
which satisfy the relation:
$$\Phi^3-64 \Psi^3=t^2.\tag 2.2$$
Furthermore, set
$$\left\{\aligned
  \varphi &=Y_2 Y_3 Y_4,\\
  \psi &=Y_2^3+Y_3^3+Y_4^3,\\
  u &=Y_0 \psi+6 Y_1 \varphi.
\endaligned\right.\tag 2.3$$
All invariants of $G_{25, 920}$ can be written as linear
combinations of the following:
$$\Phi, u, t, \Psi_1, C_6, C_9, \Phi_3, t_3, C_{12}, C_{18},$$
where
$$\left\{\aligned
  \Psi_1 &=\psi(-Y_0^3+4 Y_1^3)+18 \varphi Y_0^2 Y_1,\\
  \Psi_2 &=-\psi^2 Y_1^2-3 \varphi \psi Y_0^2+18 \varphi^2 Y_0 Y_1,\\
  \Phi_3 &=-\psi^3 Y_0+18 \varphi \psi^2 Y_1+108 \varphi^3 Y_0,\\
     t_3 &=\psi^3(Y_0^3+8 Y_1^3)-54 \varphi \psi^2 Y_0^2 Y_1+324
           \varphi^2 \psi Y_0 Y_1^2+216 \varphi^3 (Y_0^3-Y_1^3),
\endaligned\right.\tag 2.4$$
and
$$\left\{\aligned
  C_6 &=Y_2^6+Y_3^6+Y_4^6-10 (Y_2^3 Y_3^3+Y_3^3 Y_4^3+Y_4^3 Y_2^3),\\
  C_9 &=(Y_2^3-Y_3^3)(Y_3^3-Y_4^3)(Y_4^3-Y_2^3),\\
  C_{12} &=(Y_2^3+Y_3^3+Y_4^3)[(Y_2^3+Y_3^3+Y_4^3)^3+216 Y_2^3 Y_3^3 Y_4^3],\\
  {\frak C}_{12} &=Y_2 Y_3 Y_4 [27 Y_2^3 Y_3^3 Y_4^3-(Y_2^3+Y_3^3+Y_4^3)^3],\\
  C_{18} &=(Y_2^3+Y_3^3+Y_4^3)^6-540 Y_2^3 Y_3^3 Y_4^3 (Y_2^3+Y_3^3+Y_4^3)^3-5832 Y_2^6 Y_3^6 Y_4^6.
\endaligned\right.\tag 2.5$$
They satisfy the following relations:
$$\left\{\aligned
  \Psi_1^2 &=16 \Psi \Psi_2+u^2 \Phi,\\
  \Psi_1^3 &=2 \Phi^2 \Phi_3-t t_3-3 u^2 \Phi \Psi_1-u^3 t,
\endaligned\right.\tag 2.6$$
and
$$\left\{\aligned
  432 C_9^2 &=C_6^3-3 C_6 C_{12}+2 C_{18},\\
  1728 {\frak C}_{12}^3 &=C_{18}^{2}-C_{12}^{3}.
\endaligned\right.\tag 2.7$$
Burkhardt calculates the following expressions (see \cite{Bu1},
\cite{Bu2}, \cite{Hu}). There are invariants of degrees $4$, $6$,
$10$, $12$, $18$ and $45$. The invariant of degree $45$ is just
the product of the $45$ reflection hyperplanes defining the
arrangement in ${\Bbb P}^4$. The others are given by:
$$J_4=\Phi+8u,$$
$$J_6=t+20 \Psi_1-8 C_6,$$
$$J_{10}=\frac{1}{24}(\Phi \Psi_1+ut+2 \Phi C_6+2 u \Psi_1-2 \Phi_3-2 u C_6),$$
$$J_{12}=\frac{1}{24}(3t \Psi_1+3u \Phi^2+19 \Psi_1^2-9 u^2 \Phi-10 C_6 t-11 t_3+9 u^3-2 C_6 \Psi_1-4 C_{12}+4 C_6^2),$$
$$\aligned
  J_{18}
=&\frac{1}{864}(72 t \Psi \Psi_2+9 u \Phi^2 \Psi_1+9 u^2 \Phi t+288 C_6 \Psi^3+\\
 &+4 \Phi^2 \Phi_3-18 t t_3-42 u^2 \Phi \Psi_1-20 u^3 t-18 C_6 t \Psi_1+\\
 &-18 C_6 \Phi^2 u+84 C_{12} t-72 u \Phi \Phi_3+162 u^3 \Psi_1+\\
 &-240 C_6 \Psi \Psi_2+12 C_6 u^2 \Phi-6 C_6^2 t+24 \Psi_1 C_{12}-36 u^2 \Phi_3+\\
 &-6 C_6 t_3-18 C_6 u^3-12 C_6^2 \Psi_1-4 C_{18}+6 C_6 C_{12}-2 C_6^3).
\endaligned$$
There is a relation between these, which takes the form
$J_{45}^2=$ rational expression in the other invariants. There is
also an invariant of degree $40$: the product of the $40$ Steiner
primes. Denote this by $J_{40}$, the following relation holds:
$$3^{33} F_{40}=[J_4^2 (2^9 J_{12}-J_4^3)-3 \cdot 2^{18} J_{10}^2]^2
  -2^{19} [J_4 J_6-3 \cdot 2^8 J_{10}][J_4^3 J_{18}-2^{11} J_{10}^3].$$

  Now, we put
$$f_0=Y_0+2 Y_1, \quad f_1=Y_0-Y_1.\tag 2.8$$
Then
$$Y_0=\frac{1}{3}(f_0+2 f_1), \quad Y_1=\frac{1}{3}(f_0-f_1).$$
Recall that (see \cite{Y3})
$$H=\psi+6 \varphi, \quad K=\psi-3 \varphi.\tag 2.9$$
Thus,
$$\psi=\frac{1}{3}(H+2K), \quad \varphi=\frac{1}{9}(H-K).$$

  We have the following theorem which says that the above
invariants of $G_{25, 920}$ can be expressed as the polynomials of
the invariants $f_0$, $f_1$, $H$ and $K$.

{\smc Theorem 2.1} (Main Theorem 1). {\it The invariants $f_0$,
$f_1$, $H$ and $K$ satisfy the following algebraic equations,
which are the form-theoretic resolvents $($algebraic resolvents$)$
of $f_0$, $f_1$, $H$ and $K$$:$
$$u=\frac{1}{3}(f_0 H+2 f_1 K).\tag 2.10$$
$$\left\{\aligned
  \Phi &=\frac{1}{9} f_0 (f_0^3+8 f_1^3),\\
  \Psi &=\frac{1}{9} f_1 (f_0^3-f_1^3),\\
     t &=-\frac{1}{27} (f_0^6-20 f_0^3 f_1^3-8 f_1^6).
  \endaligned\right.\tag 2.11$$
$$\left\{\aligned
  \Psi_1 &=\frac{1}{9} [(f_0^3-4 f_1^3) H-6 f_0^2 f_1 K],\\
  \Psi_2 &=\frac{1}{9} (-f_0^2 H K+2 f_0 f_1 K^2-f_1^2 H^2),\\
  \Phi_3 &=\frac{1}{9} [f_0 (H^3-4 K^3)-6 f_1 H^2 K],\\
     t_3 &=\frac{1}{27} [f_0^3 (H^3+8 K^3)-18 f_0^2 f_1 H^2 K+
           36 f_0 f_1^2 H K^2+8 f_1^3 (H^3-K^3)]\\
         &=\frac{1}{27} [(f_0^3+8 f_1^3) H^3-18 f_0^2 f_1 H^2 K+
           36 f_0 f_1^2 H K^2+8 (f_0^3-f_1^3) K^3].
\endaligned\right.\tag 2.12$$}

{\it Proof}. We have
$$\aligned
  u &=\frac{1}{9}[(f_0+2f_1)(H+2K)+2(f_0-f_1)(H-K)]\\
    &=\frac{1}{3}(f_0 H+2 f_1 K).
\endaligned$$
$$\aligned
  \Phi &=\frac{1}{81}(f_0+2 f_1)[(f_0+2 f_1)^3+8(f_0-f_1)^3]\\
       &=\frac{1}{9}(f_0+2 f_1) f_0 (f_0^2-2 f_0 f_1+4 f_1^2)\\
       &=\frac{1}{9} f_0 (f_0^3+8 f_1^3).
\endaligned$$
$$\aligned
  \Psi &=\frac{1}{81}(f_0-f_1)[(f_0+2 f_1)^3-(f_0-f_1)^3]\\
       &=\frac{1}{9} (f_0-f_1) f_1 (f_0^2+f_0 f_1+f_1^2)\\
       &=\frac{1}{9} f_1 (f_0^3-f_1^3).
\endaligned$$
$$\aligned
  t &=\frac{1}{729}[(f_0+2 f_1)^6-20 (f_0+2 f_1)^3 (f_0-f_1)^3-8 (f_0-f_1)^6]\\
    &=-\frac{1}{27}(f_0^6-20 f_0^3 f_1^3-8 f_1^6).
\endaligned$$
$$\aligned
  \Psi_1
=&\frac{1}{9} \psi (f_0^3-6 f_0^2 f_1-4 f_1^3)+
  \frac{2}{3} \varphi (f_0^3+3 f_0^2 f_1-4 f_1^3)\\
=&\frac{1}{9} f_0^3 (\psi+6 \varphi)-
  \frac{2}{3} f_0^2 f_1 (\psi-3 \varphi)-
  \frac{4}{9} f_1^3 (\psi+6 \varphi)\\
=&\frac{1}{9} (f_0^3 H-6 f_0^2 f_1 K-4 f_1^3 H).
\endaligned$$
$$\aligned
  \Psi_2
=&\frac{1}{9} [-\psi^2 (f_0^2-2 f_0 f_1+f_1^2)-3 \varphi \psi
  (f_0^2+4 f_0 f_1+4 f_1^2)+18 \varphi^2 (f_0^2+f_0 f_1-2 f_1^2)]\\
=&-\frac{1}{9} f_0^2 (\psi+6 \varphi)(\psi-3 \varphi)
  +\frac{2}{9} f_0 f_1 (\psi-3 \varphi)^2
  -\frac{1}{9} f_1^2 (\psi+6 \varphi)^2\\
=&\frac{1}{9} (-f_0^2 HK+2 f_0 f_1 K^2-f_1^2 H^2).
\endaligned$$
$$\Phi_3=\frac{1}{3} f_0 (-\psi^3+18 \varphi \psi^2+108 \varphi^3)+
         \frac{2}{3} f_1 (-\psi^3-9 \varphi \psi^2+108 \varphi^3).$$
Here,
$$-\psi^3+18 \varphi \psi^2+108 \varphi^3
 =\frac{1}{3} (H^3-4 K^3),$$
$$-\psi^3-9 \varphi \psi^2+108 \varphi^3=-H^2 K.$$
Hence,
$$\Phi_3=\frac{1}{9}[f_0 (H^3-4K^3)-6 f_1 H^2 K].$$
$$\aligned
  t_3
=&\frac{1}{3} \psi^3 (f_0^3-2 f_0^2 f_1+4 f_0 f_1^2)-
  2 \varphi \psi^2 (f_0^3+3 f_0^2 f_1-4 f_1^3)+\\
 &+12 \varphi^2 \psi (f_0^3-3 f_0 f_1^2+2 f_1^3)+
  72 \varphi^3 (f_0^2 f_1+f_0 f_1^2+f_1^3)\\
=&\frac{1}{3} f_0^3 (\psi^3-6 \varphi \psi^2+36 \varphi^2 \psi)-
  \frac{2}{3} f_0^2 f_1 (\psi^3+9 \varphi \psi^2-108 \varphi^3)+\\
 &+\frac{4}{3} f_0 f_1^2 (\psi^3-27 \varphi^2 \psi+54 \varphi^3)+
  8 f_1^3 (\varphi \psi^2+3 \varphi^2 \psi+9 \varphi^3).
\endaligned$$
Here,
$$\psi^3-6 \varphi \psi^2+36 \varphi^2 \psi=\frac{1}{9} (H^3+8K^3),$$
$$\psi^3+9 \varphi \psi^2-108 \varphi^3=H^2 K,$$
$$\psi^3-27 \varphi^2 \psi+54 \varphi^3=H K^2,$$
$$\varphi \psi^2+3 \varphi^2 \psi+9 \varphi^3=\frac{1}{27}(H^3-K^3).$$
Hence,
$$t_3=\frac{1}{27} [f_0^3 (H^3+8 K^3)-18 f_0^2 f_1 H^2 K+
      36 f_0 f_1^2 H K^2+8 f_1^3 (H^3-K^3)].$$
$\qquad \qquad \qquad \qquad \qquad \qquad \qquad \qquad \qquad
 \qquad \qquad \qquad \qquad \qquad \qquad \qquad \qquad \qquad
 \quad \boxed{}$

  Note that (see \cite{Y3})
$$\left\{\aligned
  C_{12} &=\frac{1}{9} H(H^3+8 K^3),\\
  {\frak C}_{12} &=\frac{1}{27} K(K^3-H^3),\\
  C_{18} &=-\frac{1}{27}(H^6-20 H^3 K^3-8 K^6),
\endaligned\right.\tag 2.13$$
we find the following nonlinear duality:
$$f_0 \longleftrightarrow H, \quad f_1 \longleftrightarrow K.\tag 2.14$$
Here, $f_0$ and $f_1$ are linear, $H$ and $K$ are cubic. The
functions $u$ and $t_3$ are invariant under the above dual
transformation. Under the dual transformation (2.14),
$$C_{12} \longleftrightarrow \Phi, \quad
  {\frak C}_{12} \longleftrightarrow -\frac{1}{3} \Psi, \quad
  C_{18} \longleftrightarrow t, \quad
  \Psi_1 \longleftrightarrow \Phi_3.\tag 2.15$$

  Recall that (see \cite{Y3})
$$Z_1=432 \frac{C_9^2}{C_6^3}, \quad Z_2=3 \frac{C_{12}}{C_6^2}.$$
Now, we put
$$Z_3=\frac{\Phi^3}{64 \Psi^3}, \quad
  Z_4=\frac{\Psi_1^2}{16 \Psi \Psi_2}, \quad
  Z_5=\frac{2 \Phi \Phi_3}{3 u^2 \Psi_1}, \quad
  Z_6=\frac{\Psi_1^2}{3 u^2 \Phi}, \quad
  Z_7=\frac{ut}{3 \Phi \Psi_1}.\tag 2.16$$
Then
$$\frac{t^2}{64 \Psi^3}=Z_3-1, \quad
  \frac{u^2 \Phi}{16 \Psi \Psi_2}=Z_4-1, \quad
  \frac{t t_3}{3 u^2 \Phi \Psi_1}=Z_5-Z_6-Z_7-1.\tag 2.17$$

  Recall that (see \cite{Y3})
$$r_1=\frac{G^2}{C_6}, \quad
  r_2=\frac{H^2}{C_6}, \quad
  r_3=\frac{K^2}{C_6},$$
where $G=(Y_2-Y_3)(Y_3-Y_4)(Y_4-Y_2)$. Put
$$r_4=\frac{f_1}{f_0}.\tag 2.18$$
By Theorem 2.1, We find that
$$Z_3=\frac{1}{64} \frac{(1+8 r_4^3)^3}{r_4^3 (1-r_4^3)^3}.\tag 2.19$$
$$Z_4=-\frac{1}{16} \frac{[(1-4 r_4^3) \sqrt{r_2}-6 r_4 \sqrt{r_3}]^2}
      {r_4 (1-r_4^3) (\sqrt{r_2 r_3}-2 r_4 r_3+r_4^2 r_2)}.\tag 2.20$$
$$Z_5=\frac{2}{3} \frac{(1+8 r_4^3) (r_2 \sqrt{r_2}-4 r_3 \sqrt{r_3}-6
      r_4 r_2 \sqrt{r_3})}{(\sqrt{r_2}+2 r_4 \sqrt{r_3})^2 [(1-4
      r_4^3) \sqrt{r_2}-6 r_4 \sqrt{r_3}]}.\tag 2.21$$
$$Z_6=\frac{1}{3} \frac{[(1-4 r_4^3) \sqrt{r_2}-6 r_4 \sqrt{r_3}]^2}
      {(\sqrt{r_2}+2 r_4 \sqrt{r_3})^2 (1+8 r_4^3)}.\tag 2.22$$
$$Z_7=-\frac{1}{3} \frac{(\sqrt{r_2}+2 r_4 \sqrt{r_3})(1-20 r_4^3-8 r_4^6)}
      {(1+8 r_4^3) [(1-4 r_4^3) \sqrt{r_2}-6 r_4 \sqrt{r_3}]}.\tag 2.23$$

  We have
$$\left\{\aligned
  \sqrt{-3} f_0(B(Y_0, Y_1, Y_2, Y_3, Y_4)) &=f_0+2 f_1,\\
  \sqrt{-3} f_1(B(Y_0, Y_1, Y_2, Y_3, Y_4)) &=f_0-f_1.
\endaligned\right.\tag 2.24$$
$$\left\{\aligned
  f_0(C(Y_0, Y_1, Y_2, Y_3, Y_4)) &=f_0,\\
  f_1(C(Y_0, Y_1, Y_2, Y_3, Y_4)) &=f_1.
\endaligned\right.\tag 2.25$$

  Note that
$$CD=\left(\matrix
     -1 &    &    &    &   \\
        &  0 & -1 &  0 &  0\\
        &  0 &  0 &  0 & -1\\
        & -1 &  0 &  0 &  0\\
        &  0 &  0 & -1 &  0
     \endmatrix\right), \quad
  DC=\left(\matrix
     -1 &    &    &    &   \\
        &  0 &  0 &  0 & -1 \\
        & -1 &  0 &  0 &  0\\
        &  0 & -1 &  0 &  0\\
        &  0 &  0 & -1 &  0
     \endmatrix\right).$$
$$C^2 D=\left(\matrix
     -1 &    &    &    &   \\
        &  0 & -1 &  0 &  0\\
        &  0 &  0 & -1 &  0\\
        &  0 &  0 &  0 & -1\\
        & -1 &  0 &  0 &  0
     \endmatrix\right), \quad
  D C^2=\left(\matrix
     -1 &    &    &    &   \\
        &  0 &  0 & -1 &  0\\
        & -1 &  0 &  0 &  0\\
        &  0 &  0 &  0 & -1\\
        &  0 & -1 &  0 &  0
     \endmatrix\right).$$
$$(CD)^2=(DC^2)^2=\left(\matrix
      1 &   &   &   &  \\
        & 0 & 0 & 0 & 1\\
        & 0 & 0 & 1 & 0\\
        & 0 & 1 & 0 & 0\\
        & 1 & 0 & 0 & 0
     \endmatrix\right).$$
$$(DC)^2=(C^2 D)^2=\left(\matrix
      1 &   &   &   &  \\
        & 0 & 0 & 1 & 0\\
        & 0 & 0 & 0 & 1\\
        & 1 & 0 & 0 & 0\\
        & 0 & 1 & 0 & 0
     \endmatrix\right).$$
The above matrices are the permutation matrices.

  On the other hand, Maschke determined the invariants of $G_{25, 920}$ acting on
${\Bbb P}^3=\{ (z_0, z_1, z_2, z_3) \}$. According to Maschke
\cite{Mas}, $G_{25920}$ is generated by the following elements:
$M$, $N$, $P$, $Q$ and $R$, where
$$M=FAFC^2 EFA^2 C, \quad N=AEA^2,$$
$$P=FAFECFA^2 CJ, \quad Q=BEJ^3, \quad R=F.$$
Here,
$$A(z_0, z_1, z_2, z_3)=(z_0, z_2, z_3, z_1),$$
$$B(z_0, z_1, z_2, z_3)=(-z_0, z_1, z_3, z_2),$$
$$C(z_0, z_1, z_2, z_3)=(z_0, z_1, \omega z_2, \omega^2 z_3),$$
$$D(z_0, z_1, z_2, z_3)=(\omega z_0, z_1, \omega z_2, \omega z_3),$$
$$E(z_0, z_1, z_2, z_3)=(z_0, \frac{1}{\sqrt{-3}}(z_1+z_2+z_3),
  \frac{1}{\sqrt{-3}}(z_1+\omega z_2+\omega^2 z_3),
  \frac{1}{\sqrt{-3}}(z_1+\omega^2 z_2+\omega z_3)),$$
$$F(z_0, z_1, z_2, z_3)=(-z_2, z_1, -z_0, -z_3),$$
$$J(z_0, z_1, z_2, z_3)=(i z_0, i z_1, i z_2, i z_3).$$
In the form of matrices, we have
$$A=\left(\matrix
    1 &   &   &  \\
      &   & 1 &  \\
      &   &   & 1\\
      & 1 &   &
    \endmatrix\right), \quad
  B=\left(\matrix
    -1 &   &   &  \\
       & 1 &   &  \\
       &   &   & 1\\
       &   & 1 &
    \endmatrix\right),$$
$$C=\left(\matrix
    1 &   &        &         \\
      & 1 &        &         \\
      &   & \omega &         \\
      &   &        & \omega^2
    \endmatrix\right), \quad
  D=\left(\matrix
    \omega &   &        &       \\
           & 1 &        &       \\
           &   & \omega &       \\
           &   &        & \omega
    \endmatrix\right),$$
$$E=\frac{1}{\sqrt{-3}}
    \left(\matrix
    \sqrt{-3} &   &          &         \\
              & 1 & 1        & 1       \\
              & 1 & \omega   & \omega^2\\
              & 1 & \omega^2 & \omega
    \endmatrix\right), \quad
  F=\left(\matrix
       &   & -1 &   \\
       & 1 &    &   \\
    -1 &   &    &   \\
       &   &    & -1
    \endmatrix\right),$$
and $J=iI$. Note that $A$, $B$, $C$, $D$, $E$, $F$, $J \in SU(4)$.

  We find that
$$B^2=F^2=(AB)^2=(BA)^2=(BFB)^2=I.$$
$$(ABF)^2=(BFA)^2=(ABF \cdot BFA)^2=-I.$$
$$BFB=FBF, \quad ABF=-FAB.$$
Here,
$$AB=\left(\matrix
     -1 &   &   &  \\
        & 0 & 0 & 1\\
        & 0 & 1 & 0\\
        & 1 & 0 & 0
    \endmatrix\right), \quad
  BA=\left(\matrix
     -1 &   &   &  \\
        & 0 & 1 & 0\\
        & 1 & 0 & 0\\
        & 0 & 0 & 1
     \endmatrix\right),$$
$$ABF=\left(\matrix
      0 & 0 & 1 &  0\\
      0 & 0 & 0 & -1\\
     -1 & 0 & 0 &  0\\
      0 & 1 & 0 &  0
      \endmatrix\right), \quad
  BFA=\left(\matrix
      0 &  0 & 0 & 1\\
      0 &  0 & 1 & 0\\
      0 & -1 & 0 & 0\\
     -1 &  0 & 0 & 0
      \endmatrix\right),$$
$$BFB=\left(\matrix
      0 &  0 &  0 & 1\\
      0 &  1 &  0 & 0\\
      0 &  0 & -1 & 0\\
      1 &  0 &  0 & 0
      \endmatrix\right), \quad
  ABF \cdot BFA=\left(\matrix
                0 & -1 & 0 &  0\\
                1 &  0 & 0 &  0\\
                0 &  0 & 0 & -1\\
                0 &  0 & 1 &  0
                \endmatrix\right).$$
They play the similar role of the permutation matrices.

\vskip 0.5 cm

\centerline{\bf 3. Hessian polyhedra and Picard modular forms}

\vskip 0.5 cm

  In his celebrated lecture \cite{Kl}, Klein set
$${\bold A}_0=z_1 z_2, \quad {\bold A}_1=z_1^2, \quad
  {\bold A}_2=-z_2^2$$
and
$$A={\bold A}_0^2+{\bold A}_1 {\bold A}_2.$$
The invariants associated to the icosahedron: $f$, $H$ and $T$
(see \cite{Kl} or \cite{Y3}) can be expressed as the functions of
${\bold A}_0$, ${\bold A}_1$ and ${\bold A}_2$:
$$B=8 {\bold A}_0^4 {\bold A}_1 {\bold A}_2-2 {\bold A}_0^2 {\bold A}_1^2
    {\bold A}_2^2+{\bold A}_1^3 {\bold A}_2^3-{\bold A}_0 ({\bold A}_1^5+
    {\bold A}_2^5),$$
$$\aligned
  C&=320 {\bold A}_0^6 {\bold A}_1^2 {\bold A}_2^2-160 {\bold A}_0^4
     {\bold A}_1^3 {\bold A}_2^3+20 {\bold A}_0^2 {\bold A}_1^4
     {\bold A}_2^4+6 {\bold A}_1^5 {\bold A}_2^5+\\
   &-4 {\bold A}_0 ({\bold A}_1^5+{\bold A}_2^5)(32 {\bold A}_0^4-20
    {\bold A}_0^2 {\bold A}_1 {\bold A}_2+5 {\bold A}_1^2 {\bold A}_2^2)
    +{\bold A}_1^{10}+{\bold A}_2^{10},
\endaligned$$
$$\aligned
  D&=({\bold A}_1^5-{\bold A}_2^5)[-1024 {\bold A}_0^{10}+3840 {\bold A}_0^8
     {\bold A}_1 {\bold A}_2-3840 {\bold A}_0^6 {\bold A}_1^2 {\bold A}_2^2+\\
   &+1200 {\bold A}_0^4 {\bold A}_1^3 {\bold A}_2^3-100 {\bold A}_0^2
     {\bold A}_1^4 {\bold A}_2^4+{\bold A}_1^{10}+{\bold A}_2^{10}+
     2 {\bold A}_1^5 {\bold A}_2^5+\\
   &+{\bold A}_0 ({\bold A}_1^5+{\bold A}_2^5)(352 {\bold A}_0^4-160
     {\bold A}_0^2 {\bold A}_1 {\bold A}_2+10 {\bold A}_1^2 {\bold A}_2^2)].
\endaligned$$
They satisfy the identity
$$D^2=-1728 B^5+C^3+720 ACB^3-80 A^2 C^2 B+64 A^3 (5 B^2-AC)^2.$$
Moreover, the Jacobian of $A$, $B$ and $C$ satisfies the relation
(see \cite{F}):
$$\frac{\partial(A, B, C)}{\partial({\bold A}_0, {\bold A}_1, {\bold A}_2)}
 =\vmatrix \format \c \quad & \c \quad & \c\\
  \frac{\partial A}{\partial {\bold A}_0} &
  \frac{\partial A}{\partial {\bold A}_1} &
  \frac{\partial A}{\partial {\bold A}_2}\\
  \frac{\partial B}{\partial {\bold A}_0} &
  \frac{\partial B}{\partial {\bold A}_1} &
  \frac{\partial B}{\partial {\bold A}_2}\\
  \frac{\partial C}{\partial {\bold A}_0} &
  \frac{\partial C}{\partial {\bold A}_1} &
  \frac{\partial C}{\partial {\bold A}_2}\\
  \endvmatrix=-10 D.$$

  The general Jacobian equation of the sixth degree is the
following:
$$(z-A)^6-4A(z-A)^5+10B(z-A)^3-C(z-A)+(5B^2-AC)=0.$$
Let $z_0, z_1, z_2, z_3, z_4, z_{\infty}$ be the six roots, put
$$y=(z_{\infty}-z_0)(z_1-z_4)(z_2-z_3).$$
Brioschi remarked that the square root of this expression is
already rational in the ${\Bbb A}$'s, and gives rise to an
equation of the fifth degree. Let $x$ be this square root,
Brioschi found for the five values of which $x$ is susceptible the
following formula:
$$\aligned
  x_{\nu}=&-\epsilon^{\nu} {\bold A}_1 (4 {\bold A}_0^2-{\bold A}_1 {\bold A}_2)
          +\epsilon^{2 \nu} (2 {\bold A}_0 {\bold A}_1^2-{\bold A}_2^3)+\\
          &+\epsilon^{3 \nu} (-2 {\bold A}_0 {\bold A}_2^2+{\bold A}_1^3)
          +\epsilon^{4 \nu} {\bold A}_2 (4 {\bold A}_0^2-{\bold A}_1 {\bold A}_2),
\endaligned$$
with $\epsilon=\exp (2 \pi i/5)$, while for the corresponding
equation of the fifth degree he found this:
$$x^5+10 B x^3+5(9B^2-AC) x-D=0.$$
In fact, the Jacobian equation of the sixth degree with $A=0$ is
none other than that simplest resolvent of the sixth degree which
Klein established in the case of the icosahedron: put ${\bold
A}_0=z_1 z_2$, ${\bold A}_1=z_1^2$, ${\bold A}_2=-z_2^2$ and
correspondingly:
$$B=-f, \quad C=-H, \quad D=T.$$ The resolvent of
the fifth degree is transformed into the following:
$$x^5+10 B x^3+45 B^2 x-D=0.$$

  The actions of the icosahedral group on ${\bold A}_0$, ${\bold A}_1$
and ${\bold A}_2$ are given by
$$S: \quad {\bold A}_0^{\prime}={\bold A}_0, \quad
     {\bold A}_1^{\prime}=\epsilon^4 {\bold A}_1, \quad
     {\bold A}_2^{\prime}=\epsilon {\bold A}_2;$$
$$T: \left\{\aligned
  \sqrt{5} {\bold A}_0^{\prime} &={\bold A}_0+{\bold A}_1+{\bold A}_2,\\
  \sqrt{5} {\bold A}_1^{\prime} &=2 {\bold A}_0+(\epsilon^2+\epsilon^3)
                        {\bold A}_1+(\epsilon+\epsilon^4) {\bold A}_2,\\
  \sqrt{5} {\bold A}_2^{\prime} &=2 {\bold A}_0+(\epsilon+\epsilon^4)
                      {\bold A}_1+(\epsilon^2+\epsilon^3) {\bold A}_2;
\endaligned\right.$$
$$U: \quad {\bold A}_0^{\prime}=-{\bold A}_0, \quad
     {\bold A}_1^{\prime}=-{\bold A}_2, \quad
     {\bold A}_2^{\prime}=-{\bold A}_1.$$
This leads to the investigation of the binary quadratic forms:
$${\bold A}_1 z_1^2+2 {\bold A}_0 z_1 z_2-{\bold A}_2 z_2^2$$
and Hilbert modular surfaces (see \cite{Hi2}). Klein (see
\cite{Kl}) set
$$\delta_{\nu}=(\epsilon^{4 \nu} {\bold A}_1-\epsilon^{\nu} {\bold A}_2)
  [(1+\sqrt{5}) {\bold A}_0+\epsilon^{4 \nu} {\bold
  A}_1+\epsilon^{\nu} {\bold A}_2] [(1-\sqrt{5}) {\bold A}_0+\epsilon^{4 \nu} {\bold
  A}_1+\epsilon^{\nu} {\bold A}_2],$$
where $\epsilon=\exp(2 \pi i/5)$, $\nu=0, 1, 2, 3, 4$ and found
that
$$\left\{\aligned
  \delta_0+\delta_1+\delta_2+\delta_3+\delta_4 &=0,\\
  \delta_0^3+\delta_1^3+\delta_2^3+\delta_3^3+\delta_4^3 &=0.
\endaligned\right.$$
Moreover,
$$S(\delta_0)=\delta_1, S(\delta_1)=\delta_2, S(\delta_2)=\delta_3,
  S(\delta_3)=\delta_4, S(\delta_4)=\delta_0.$$
$$U(\delta_0)=\delta_0, U(\delta_1)=\delta_4, U(\delta_2)=\delta_3,
  U(\delta_3)=\delta_2, U(\delta_4)=\delta_1.$$
$T(\delta_0)=\delta_0$. However, there does not exist the
analogous transformation formulas for $T(\delta_i)$ $(i=1, 2, 3,
4)$.

  It is well-known that the group $A_5$ is isomorphic to the finite
group $I$ of those elements of $SO(3)$ which carry a given
icosahedron centered at the origin of the standard Euclidean space
${\Bbb R}^3$ to itself. The group $I$ operates linearly on ${\Bbb
R}^3$ (standard coordinates $x_0$, $x_1$, $x_2$) and thus also on
${\Bbb R} {\Bbb P}^2$ and ${\Bbb C} {\Bbb P}^2$. We are concerned
with the action on ${\Bbb C} {\Bbb P}^2$ (see \cite{Hi2} and
\cite{Hi3}). A curve in ${\Bbb C} {\Bbb P}^2$ which is mapped to
itself by all elements of $I$ is given by a homogeneous polynomial
in $x_0$, $x_1$, $x_2$ which is $I$-invariant up to constant
factors and hence $I$-invariant, because $I$ is a simple group.
The graded ring of all $I$-invariant polynomials in $x_0$, $x_1$,
$x_2$ is generated by homogeneous polynomials $A$, $B$, $C$, $D$
of degrees $2$, $6$, $10$, $15$ with $A=x_0^2+x_1^2+x_2^2$. The
action of $I$ on ${\Bbb C} {\Bbb P}^2$ has exactly one minimal
orbit where ``minimal'' means that the number of points in the
orbit is minimal. This orbit has six points, they are called
poles. These are the points of ${\Bbb R} {\Bbb P}^2 \subset {\Bbb
C} {\Bbb P}^2$ which are represented by the six axes through the
vertices of the icosahedron. Klein uses coordinates
$$A_0=x_0, \quad A_1=x_1+i x_2, \quad A_2=x_1-i x_2$$
and puts the icosahedron in such a position that the six poles are
given by
$$(A_0, A_1, A_2)=\left(\frac{\sqrt{5}}{2}, 0, 0\right), \quad
  \left(\frac{1}{2}, \varepsilon^{\nu}, \varepsilon^{-\nu}\right)$$
with $\varepsilon=\exp(2 \pi i/ 5)$ and $0 \leq \nu \leq 4$.

  The invariant curve $A=0$ does not pass through the poles. There is
exactly one invariant curve $B=0$ of degree $6$ which passes
through the poles, exactly one invariant curve $C=0$ of degree
$10$ which has higher multiplicity than the curve $B=0$ in the
poles and exactly one invariant curve $D=0$ of degree $15$. In
fact, $B=0$ has an ordinary double point (multiplicity $2$) in
each pole, $C=0$ has a double cusp (multiplicity $4$) in each pole
and $D=0$ is the union of the $15$ lines connecting the six poles.
Klein gives formulas for the homogeneous polynomials $A$, $B$,
$C$, $D$ (determined up to constant factors). They generate the
ring of all $I$-invariant polynomials. According to Klein the ring
of $I$-invariant polynomials is given as follows:
$${\Bbb C}[A_0, A_1, A_2]^{I}={\Bbb C}[A, B, C, D]/(R(A, B, C, D)=0).$$
The relation $R(A, B, C, D)=0$ is of degree $30$.

  The equations for $B$ and $C$ show that the two tangents of $B=0$
in the pole $(\sqrt{5}/2, 0, 0)$ are given by $A_1=0$, $A_2=0$.
They coincide with the tangents of $C=0$ in that pole. Therefore
the curves $B=0$ and $C=0$ have in each pole the intersection
multiplicity $10$. Thus they intersect only in the poles.

  When we restrict the action of $I$ to the conic $A=0$, we get
the well-known action of $I$ on ${\Bbb C} {\Bbb P}^1$. The curves
$B=0$, $C=0$, $D=0$ intersect $A=0$ transversally in $12$, $20$,
$30$ points respectively.

  Put
$$\Gamma(\sqrt{5})=\left\{ \left(\matrix \alpha & \beta\\
  \gamma & \delta \endmatrix\right) \in SL(2, {\Cal O}_K):
  \alpha \equiv \delta \equiv 1 (\text{mod $\sqrt{5}$}),
  \beta \equiv \gamma \equiv 0 (\text{mod $\sqrt{5}$}) \right\},$$
where $K={\Bbb Q}(\sqrt{5})$.

{\smc Theorem 3.1} (Hirzebruch) (see \cite{Hi2} and \cite{Eb}).
{\it The ring of symmetric Hilbert modular forms for
$\Gamma(\sqrt{5})$ is equal to ${\Bbb C}[A_0, A_1, A_2]$ where
$A_0$, $A_1$, $A_2$ have weight $1$.}

{\smc Corollary 3.2} (see \cite{Hi2} and \cite{Eb}). {\it The ring
of symmetric Hilbert modular forms for $SL(2, {\Cal O}_K)$ of even
weight is equal to ${\Bbb C}[A, B, C]$.}

  Let
$$\Gamma(2)=\left\{ \left(\matrix \alpha & \beta\\
  \gamma & \delta \endmatrix\right) \in SL(2, {\Cal O}_K):
  \alpha \equiv \delta \equiv 1 (\text{mod $2$}),
  \beta \equiv \gamma \equiv 0 (\text{mod $2$}) \right\},$$ where $K={\Bbb Q}(\sqrt{5})$.

{\smc Theorem 3.3} (see \cite{Hi1}). {\it The graded ring of
Hilbert modular forms of even weight relative to the group
$\Gamma(2)$ is the ring of polynomials in the five Eisenstein
series $E_0$, $\cdots$, $E_4$, with the relations
$$\sigma_2(E_0, \cdots, E_4)=0, \quad
  \sigma_4(E_0, \cdots, E_4)=0$$
and their consequences, but no other relations. Here, $\sigma_j$
is the $j$-th elementary symmetric function of $E_0, \cdots,
E_4$.}

{\smc Corollary 3.4} (see \cite{Hi1}). {\it The ring $M(SL(2,
{\Cal O}_K))$ is generated by the four modular forms:
$$\sigma_1=\sum_{0 \leq i \leq 4} E_i, \quad
  \sigma_3=\sum_{0 \leq i<j<k \leq 4} E_i E_j E_k, \quad
  \sigma_5=\prod_{0 \leq i \leq 4} E_i, \quad
  \Delta=\prod_{0 \leq i<j \leq 4} (E_i-E_j)$$
of weight $2$, $6$, $10$, and $20$, respectively, subject to the
single relation:
$$\aligned
  \Delta^2=&\sigma_5 (3125 \sigma_5^3+2000 \sigma_5^2 \sigma_3 \sigma_1^2+256 \sigma_5^2 \sigma_1^5+\\
           &-900 \sigma_5 \sigma_3^3 \sigma_1-128 \sigma_5 \sigma_3^2 \sigma_1^4+16 \sigma_3^4 \sigma_1^3
            +108 \sigma_3^5).
\endaligned$$}

  Moreover, Hirzebruch found the connection between lattices, theta
functions, coding theory and Hilbert modular surface associated
with ${\Bbb Q}(\sqrt{5})$ (see \cite{Hi}, pp.796-798). Let $p$ be
an odd prime and $\zeta=e^{2 \pi i/p}$. The field $K={\Bbb
Q}(\zeta)$ contains the totally real field $k={\Bbb
Q}(\zeta+\zeta^{-1})$. For the rings of integers ${\Cal O}_{K}$
and ${\Cal O}_k$ we have
$${\Cal O}_K={\Bbb Z} \zeta+\cdots+{\Bbb Z} \zeta^{p-1},$$
$${\Cal O}_k={\Bbb Z}(\zeta+\zeta^{-1})+\cdots+{\Bbb Z}(\zeta^{\frac{p-1}{2}}+\zeta^{\frac{p+1}{2}}).$$
The Galois group $G$ of $k$ over ${\Bbb Q}$ is cyclic of order
$\frac{p-1}{2}$. Let  $\text{Tr}: k \to {\Bbb Q}$ be the trace,
i.e.
$$\text{Tr}(y)=\sum_{\sigma \in G} \sigma(y), \quad \text{for} \quad y \in k.$$
We consider the following quadratic form $Q$ over the ${\Bbb
Z}$-lattice ${\Cal O}_K$
$$Q: {\Cal O}_K \to {\Bbb Q}, \quad Q(x)=2 \text{Tr}\left(\frac{x \overline{x}}{p}\right).$$
It has discriminant $\frac{1}{p}$.

  The ideal $(1-\zeta)={\frak p}$ in ${\Cal O}_K$ has norm $p$.
Let $\varrho: {\Cal O}_K \to {\Bbb F}_p=GF(p)$ be the
corresponding homomorphism. We also consider $\varrho: {\Cal
O}_K^n \to {\Bbb F}_p^n$. For a self-dual code $C \subset {\Bbb
F}_p^n$ ($n$ even) we define the lattice
$\Lambda=\varrho^{-1}(C)$. If we restrict $Q$ to $\Lambda$, then
we obtain an integral even unimodular form on the lattice
$\Lambda$. It is well-known that $n(p-1) \equiv 0$ $($mod $8)$.

  The theta functions in the Hilbert modular sense for the totally real
field $k$ depend on $\frac{p-1}{2}$ variables $z_{\sigma} \in
{\Bbb H}$ (upper half plane) indexed by the elements of the Galois
group $G$ of $k$ over ${\Bbb Q}$. For $y \in k$
$$\text{Tr}(zy) \triangleq \sum_{\sigma \in G} z_{\sigma} y_{\sigma}
  \quad (y_{\sigma}=\sigma(y)).$$
The fundamental theta functions are
$$\theta_j(z)=\sum_{x \in {\frak p}+j} e^{2 \pi i \text{Tr}
  \left(z \frac{x \overline{x}}{p}\right)} \quad \text{for}
  \quad j=0, 1, \cdots, \frac{p-1}{2}.$$

  The $\theta_j$ are Hilbert modular forms of weight $1$ for the
group
$$SL_2({\Cal O}_k, \widetilde{{\frak p}})=\{ A | A \in SL_2({\Cal O}_k)
  \quad \text{and} \quad \text{$A \equiv 1$ mod $\widetilde{{\frak p}}$} \}$$
where $\widetilde{{\frak p}}$ is the prime ideal in $k$ lying over
$p$
$$\widetilde{{\frak p}}=((\zeta-\overline{\zeta})^2), \quad
  \widetilde{{\frak p}}^{\frac{p-1}{2}}=(p).$$
In fact, the group $SL_2({\Bbb F}_p)$ acts on the
$\frac{p+1}{2}$-dimensional vector space over ${\Bbb C}$ generated
by the $\theta_j$. For $p \equiv 1$ mod $4$ this is an action of
$PSL_2({\Bbb F}_p)$.

  The lattice $\Lambda$ gives a theta function
$$\theta_C(z)=\sum_{x \in \Lambda} e^{2 \pi i \text{Tr} \left(z \frac{x \overline{x}}{p}\right)}.$$
Consider the Lee weight enumerator $L_C(X_0, \cdots,
X_{\frac{p-1}{2}})$.

{\smc Theorem 3.5} (see \cite{Hi}, p.797). {\it
$$\theta_C=L_C(\theta_0, \cdots, \theta_{\frac{p-1}{2}}).$$
$\theta_C$ is a Hilbert modular form of weight $n$ for the group
$SL_2({\Cal O}_k)$. The polynomial $L_C(X_0, \cdots,
X_{\frac{p-1}{2}})$ is an invariant polynomial for the above
mentioned representation of dimension $\frac{p+1}{2}$ of the group
$SL_2({\Bbb F}_p)$. All Hilbert modular forms $\theta_C$,
$\theta_j$ occurring are symmetric, i.e. invariant under the
action of the Galois group $G$ on the variables $z_{\sigma}$.}

  In his paper on ${\Bbb Q}(\sqrt{5})$, which is $k$ for $p=5$,
Hirzebruch proved that the ring of symmetric Hilbert modular forms
for $SL_2({\Cal O}_k, \widetilde{{\frak p}})$, $\widetilde{{\frak
p}}=(\sqrt{5})$, equals ${\Bbb C}[A_0, A_1, A_2]$ where $A_0$,
$A_1$, $A_2$ have weight $1$. He proved that the ring of symmetric
Hilbert modular forms for $SL_2({\Cal O}_k)$ of even weight equals
${\Bbb C}[A, B, C]$ where $A$, $B$, $C$ are the Klein invariants
of degrees $2$, $6$, $10$. Moreover, he found that
$$A_0=\theta_0, \quad A_1=2 \theta_1, \quad A_2=2 \theta_2.$$
Together with Theorem 3.5, this explains the above mentioned
connection.

  For $p=3$ (where $k={\Bbb Q}$). In this case $C$ is a ternary code.
According to \cite{Eb}, the Lee weight enumerator of the code $C$
satisfies the following identity:
$$W_C(\theta_0, \theta_1)=\theta_C.$$
We consider the functions $\theta_0$ and $\theta_1$. The lattice
$\Lambda$ is isomorphic to $A_2$. The lattice $A_2$ is the ${\Bbb
Z}$-module ${\Bbb Z}^2$ with the quadratic form
$$(x, y)^2=2x^2-2xy+2y^2.$$
Therefore,
$$\theta_0(z)=\sum_{(x, y) \in {\Bbb Z}^2} q^{x^2-xy+y^2}, \quad
  \theta_1(z)=q^{\frac{1}{3}} \sum_{(x, y) \in {\Bbb Z}^2}
              q^{x^2-xy+y^2+x-y}\tag 3.1$$
with $q=e^{2 \pi iz}$. The functions $\theta_0$ and $\theta_1$ are
modular forms of weight $1$ for the group
$$\Gamma(3)=\left\{ \left(\matrix a & b\\ c & d \endmatrix\right)
            \in SL_2({\Bbb Z}): a \equiv d \equiv 1 (\text{mod
            $3$}), b \equiv c \equiv 0 (\text{mod $3$}) \right\}.$$
We have
$$\theta_0 \left(-\frac{1}{z}\right)=\frac{z}{\sqrt{-3}}(\theta_0(z)+2
  \theta_1(z)), \quad
  \theta_1 \left(-\frac{1}{z}\right)=\frac{z}{\sqrt{-3}}(\theta_0(z)-
  \theta_1(z)).$$
$$\theta_0(z+1)=\theta_0(z), \quad \theta_1(z+1)=\omega \theta_1(z),
  \quad \omega=e^{\frac{2 \pi i}{3}}.$$
The Eisenstein series $E_4$ is the theta function of the
unimodular $E_8$-lattice and
$$E_4=\theta_0^4+8 \theta_0 \theta_1^3, \quad
  E_6=\theta_0^6-20 \theta_0^3 \theta_1^3-8 \theta_1^6.\tag 3.2$$

{\smc Theorem 3.6} (see \cite{Eb}). {\it The algebra of all
modular forms for the group $\Gamma(3)$ is isomorphic to the
polynomial algebra ${\Bbb C}[\theta_0, \theta_1]$.}

  The group $SL_2({\Bbb F}_3)$ acts on the polynomial algebra
${\Bbb C}[\theta_0, \theta_1]$. The ring of invariant polynomials
under this action is denoted by ${\Bbb C}[\theta_0,
\theta_1]^{SL_2({\Bbb F}_3)}$.

{\smc Theorem 3.7} (see \cite{Eb}). {\it
$${\Bbb C}[\theta_0, \theta_1]^{SL_2({\Bbb F}_3)}={\Bbb C}[E_4, E_6].$$}

  The Hamming weight enumerator $L_C=H_C=H_C(\theta_0, \theta_1)$
is a polynomial in the modular forms $E_4$ and $E_6^2$. For the
ternary Golay code of length $12$,
$$H_C=\frac{5}{8} E_4^3+\frac{3}{8} E_6^2.$$

  According to \cite{CS}, the Hamming weight enumerator classifies
codewords according to the number of nonzero coordinates. More
detailed information is supplied by the complete weight enumerator
(c.w.e.), which gives the number of codewords of each composition.
For example the c.w.e. of a ternary code $C$ is
$$c.w.e._{C}(x, y, z)=\sum_{u \in C} x^{n_0(u)} y^{n_1(u)} z^{n_{-1}(u)},$$
where $n_i(u)$ is the number of times $i \in {\Bbb F}_3$ occurs in
$u$. There is also a version of the MacWilliams identity for
c.w.e.'s.

{\smc Theorem 3.8} (Gleason) (see \cite{CS}). {\it If $C$ is a
Type III code then $W_C(x, y)$ is invariant under a group $G_3$ of
order $48$, and belongs to ${\Bbb C}[\psi_4, \xi_{12}]$, where
$$\psi_4=x^4+8 x y^3$$
is the weight enumerator of ${\Cal C}_{4}$ and
$$\xi_{12}=y^3 (x^3-y^3)^3.$$
Equivalently, $W_C(x, y)$ can be written as a polynomial in the
weight enumerators of ${\Cal C}_4$ and ${\Cal C}_{12}$.}

  To describe the complete weight enumerator, we first define some
further polynomials. Let
$$a=x^3+y^3+z^3, \quad p=3xyz,$$
$$b=x^3 y^3+y^3 z^3+z^3 x^3,$$
$$\beta_6=a^2-12 b=x^6+y^6+z^6-10(x^3 y^3+y^3 z^3+z^3 x^3),$$
$$\pi_9=(x^3-y^3)(y^3-z^3)(z^3-x^3),$$
$$\alpha_{12}=a(a^3+8p^3).$$

{\smc Theorem 3.9} (see \cite{CS}). {\it If $C$ is a Type III code
which contains the all-ones vector then the complete weight
enumerator is invariant under a group $G_4$ of order $2592$, and
belongs to
$$R \oplus \beta_6 \pi_9^2 R,$$
where
$$R={\Bbb C}[\beta_6^2, \alpha_{12}, \pi_9^4].$$
In other words the complete weight enumerator can be written
uniquely as a polynomial in $\beta_6^2$, $\alpha_{12}$ and
$\pi_9^4$, plus $\beta_6 \pi_9^2$ times another such polynomial.}

  Here, $G_3$ is generated by $\frac{1}{\sqrt{3}} \left(\matrix 1 & 2\\
1 & -1 \endmatrix\right)$, $\left(\matrix 1 & 0\\ 0 & \omega
\endmatrix\right)$ and is the reflection group $3[6]2$. $G_4$ is
generated by all permutations, $\text{diag}(1, 1, \omega)$ and
$$\frac{1}{\sqrt{3}} \left(\matrix
  1 & 1 & 1\\ 1 & \omega & \overline{\omega}\\
  1 & \overline{\omega} & \omega
  \endmatrix\right),$$
has a center $Z(G_4)$ of order $12$, and $G_4/Z(G_4)$ is the
Hessian group of order $216$.

  It is well-known that a complete system of invariants for the
Hessian groups (the corresponding geometric objects are Hessian
polyhedra) has degrees $6$, $9$, $12$, $12$ and $18$ and can be
given explicitly by the following forms:
$$\left\{\aligned
  C_6(z_1, z_2, z_3) &=z_1^6+z_2^6+z_3^6-10
                     (z_1^3 z_2^3+z_2^3 z_3^3+z_3^3 z_1^3),\\
  C_9(z_1, z_2, z_3) &=(z_1^3-z_2^3)(z_2^3-z_3^3)(z_3^3-z_1^3),\\
  C_{12}(z_1, z_2, z_3) &=(z_1^3+z_2^3+z_3^3)[(z_1^3+z_2^3+z_3^3)^3
                        +216 z_1^3 z_2^3 z_3^3],\\
  {\frak C}_{12}(z_1, z_2, z_3) &=z_1 z_2 z_3 [27 z_1^3 z_2^3 z_3^3
                                -(z_1^3+z_2^3+z_3^3)^3],\\
  C_{18}(z_1, z_2, z_3) &=(z_1^3+z_2^3+z_3^3)^6-540 z_1^3 z_2^3 z_3^3
                        (z_1^3+z_2^3+z_3^3)^3-5832 z_1^6 z_2^6 z_3^6.
\endaligned\right.\tag 3.3$$
They satisfy the following relations:
$$\left\{\aligned
  432 C_9^2 &=C_6^3-3 C_6 C_{12}+2 C_{18},\\
  1728 {\frak C}_{12}^3 &=C_{18}^{2}-C_{12}^{3}.
\endaligned\right.\tag 3.4$$
For the Hessian group $3[3]3[3]3$, of order $6 \times 9 \times
12$, the three forms are $C_6$, $C_9$, $C_{12}$. For the closely
related Hessian group $2[4]3[3]3$, of order $6 \times 12 \times
18$, they are $C_6$, $C_{12}$, $C_9^2$ (see \cite{Co}).

  Note that
$$\beta_6=C_6, \quad \pi_9=C_9, \quad \alpha_{12}=C_{12}, \quad
  R={\Bbb C}[C_6^2, C_{12}, C_9^4]\tag 3.5$$
under the identification
$$(x, y, z)=(z_1, z_2, z_3).$$
On the other hand,
$$\psi_4=9 C_{12}, \quad \xi_{12}=-27^3 {\frak C}_{12}^3\tag 3.6$$
under the identification
$$(x, y)=(H, K).$$
Thus, we find  a connection between Hessian invariants (Hessian
groups) and coding theory.

  Two theories apparently far apart, that of weight polynomials of
codes over the field ${\Bbb F}_2$ and that of theta functions of
lattices over ${\Bbb Z}$, are too similar to be strangers to each
other. The most striking similarity is no doubt the following: The
weight polynomials of certain codes $\cdots$ generate the graded
algebra of polynomials in two variables fixed by the operation of
a subgroup $H$ of $GL(2, {\Bbb C})$; the theta functions of even
unimodular lattices are modular forms which generate an algebra
visibly isomorphic to the preceding (see \cite{BrE}). Since the
times of Klein the analogies between invariant theory and modular
forms have emerged.

  Denote $E_2(z)$, $E_4(z)$, $E_6(z)$ the Eisenstein series
of order $2$, $4$, $6$ and the cusp form $\Delta(z)$ of weight
$12$ as
$$E_2(z)=1-24 \sum_{m=1}^{\infty} \sigma_1(m) q^m,$$
$$E_4(z)=1+240 \sum_{m=1}^{\infty} \sigma_3(m) q^m,$$
$$E_6(z)=1-504 \sum_{m=1}^{\infty} \sigma_5(m) q^m,$$
$$\Delta(z)=\frac{1}{1728}(E_4^3-E_6^2)=q \prod_{m=1}^{\infty} (1-q^m)^{24},$$
where $q=\exp(2 \pi i z)$, and the sum of $r$th power of divisor
function is $\sigma_r(m):=\sum_{d|m} d^r$. As is well known $E_4$,
$E_6$ are modular forms of weight $4$ and $6$ but $E_2$ is not.
The following result is due to Hecke.

{\smc Theorem 3.10} (see \cite{CMS}). {\it The algebra of modular
forms of weight multiple of $4$ is ${\Bbb C}[E_4, \Delta]$. The
algebra of modular forms of even weight is ${\Bbb C}[E_4, E_6]$.}

  Let $H_2=\langle M_2, N_2 \rangle$, $M_2=\frac{1}{\sqrt{2}}
\left(\matrix 1 & 1\\ 1 & -1 \endmatrix\right)$,
$N_2=\left(\matrix 1 & 0\\ 0 & i \endmatrix\right)$. $H_2$ is of
order $192$. There is a subgroup $G_2 \leq H_2$ of index $2$
defined as the kernel of the following linear character defined on
the generators of $H_2$ as
$$\chi(M_2)=-1, \quad \chi(N_2)=1.$$
Define the following invariants for the group $H_2$ of degree $8$
and $24$, respectively:
$$\psi_8=x^8+14 x^4 y^4+y^8,$$
$$\nu_{24}=x^4 y^4 (x^4-y^4)^4.$$
Then
$${\Bbb C}[x, y]^{H_2}={\Bbb C}[\psi_8, \nu_{24}],$$
$${\Bbb C}[x, y]^{G_2}={\Bbb C}[\psi_8, k_{12}],$$
where
$$k_{12}=x^{12}-33(x^8 y^4+x^4 y^8)+y^{12}.$$

{\smc Theorem 3.11} (see \cite{CMS}). {\it The map
$$\phi_1: {\Bbb C}[\psi_8, \nu_{24}] \to {\Bbb C}[E_4, \Delta]$$
defined by
$$\phi_1(h(\psi_8, \nu_{24}))=h(E_4, \Delta)$$
is an algebra isomorphism. The map
$$\phi_2: {\Bbb C}[\psi_8, k_{12}] \to {\Bbb C}[E_4, E_6]$$
defined by
$$\phi_2(h(\psi_8, k_{12}))=h(E_4, E_6)$$
is an algebra isomorphism.}

  Let $G_3=\langle M_3, N_3 \rangle$, $M_3=\frac{1}{\sqrt{3}}
\left(\matrix 1 & 2\\ 1 & -1 \endmatrix\right)$,
$N_3=\left(\matrix 1 & 0\\ 0 & \omega \endmatrix\right)$,
$\omega=\exp(2 \pi i/3)$. Then $G_3\cong SL_2({\Bbb F}_3)$.
Primary invariants for $G_3$ are
$$\psi_4=x^4+8 xy^3,$$
$$k_6=x^6-20 x^3 y^3-8 y^6.$$
$${\Bbb C}[x, y]^{G_3}={\Bbb C}[\psi_4, k_6].$$

{\smc Theorem 3.12} (see \cite{CMS}). {\it The map
$$\phi_3: {\Bbb C}[\psi_4, k_6] \to {\Bbb C}[E_4, E_6]$$
given by
$$\phi_3(h(\psi_4, k_6))=h(E_4, E_6)$$
is an algebra isomorphism.}

  In our case, we have the following:

{\smc Theorem 3.13}. {\it
$${\Bbb C}[H, K]^{G_3}={\Bbb C}[C_{12}, C_{18}].$$
The map
$$\varphi: {\Bbb C}[C_{12}, C_{18}] \to {\Bbb C}[E_4, E_6]$$
given by
$$\varphi(h(C_{12}, C_{18}))=h(E_4, E_6)$$
is an algebra isomorphism.}

  Consequently,
$$\varphi(h(C_{12}, -{\frak C}_{12}^3))=h(E_4, \Delta).$$

  Now, we will study the analogies between Hessian invariants and
Picard modular forms, which can be considered as the higher
dimensional counterpart of the analogies between regular
polyhedral invariants and elliptic modular forms studied by Klein.
Let us recall some basic facts about Picard modular forms. In his
monograph \cite{Ho1}, Holzapfel proved the following results:

{\smc Theorem 3.14} (see \cite{Ho1}). {\it \roster \item The
monodromy group of the hypergeometric differential equation system
considered by Picard is the principal congruence subgroup $U(2, 1;
{\Bbb Z}[\omega])(1-\omega)$ of the full Eisenstein lattice of the
ball $U(2, 1; {\Bbb Z}[\omega])$ with respect to the principal
ideal ${\Bbb Z}[\omega] \cdot (1-\omega)$, $\omega=\exp(2 \pi
i/3)$.

\item The rings of automorphic forms of $U(2, 1; {\Bbb Z}[\omega])
(1-\omega)$ and $U(2, 1; {\Bbb Z}[\omega])$ (of Nebentypus $\chi$)
are polynomial rings
$$\bigoplus_{m=0}^{\infty} [U(2, 1; {\Bbb Z}[\omega])(1-\omega),
  m]_{\chi}={\Bbb C}[\xi_1, \xi_2, \xi_3],$$
$$\bigoplus_{m=0}^{\infty} [U(2, 1; {\Bbb Z}[\omega]), m]_{\chi}
 ={\Bbb C}[G_2, G_3, G_4],$$
with $\xi_1$, $\xi_2$, $\xi_3$ of weight $1$, $G_k$ of weight $k$,
$k=2, 3, 4$.

  A key result is the following one:
$$\bigoplus_{m=0}^{\infty} [SU(2, 1; {\Bbb Z}[\omega])(1-\omega), m]
 ={\Bbb C}[\xi_1, \xi_2, \xi_3, \zeta],$$
with $\xi_k$ as above for $k=1, 2, 3$, $\zeta$ of weight $2$,
$$\zeta^3=(\xi_3^2-\xi_2^2)(\xi_3^2-\xi_1^2)(\xi_2^2-\xi_1^2).$$

\item ${\Bbb B}^2/U(2, 1; {\Bbb Z}[\omega])(1-\omega) \cong {\Bbb
P}^2 -\{ \text{four points in general position} \}$.

\item If $F_1$, $F_2$, $F_3$ are three fundamental solutions of
Picard's system of differential equations, then there is a basis
$\xi_1^{\prime}$, $\xi_2^{\prime}$, $\xi_3^{\prime}$ of the space
of automorphic forms of weight $1$ with respect to a group
isomorphic to $U(2, 1; {\Bbb Z}[\omega])(1-\omega)$, such that
$(\xi_1^{\prime}: \xi_2^{\prime}: \xi_3^{\prime})$ is the inverse
of the multi-valued map $(F_1: F_2: F_3)$:
$$\CD
 \widetilde{{\Bbb B}^2} @>=>> \widetilde{{\Bbb B}^2}\\
 @AA(F_1:F_2:F_3)A      @VV(\xi_1^{\prime}: \xi_2^{\prime}: \xi_3^{\prime})V \\
 {\Bbb P}^2-\{ \text{four points} \} @>=>> {\Bbb P}^2-\{
 \text{four points} \}
\endCD$$
where $\widetilde{{\Bbb B}^2}$ is projectively equivalent to
${\Bbb B}^2$ in ${\Bbb P}^2$.

\item Let $\widehat{{\Bbb B}^2/U(2, 1; {\Bbb Z}[\omega])}=
\overline{{\Bbb B}^2}/U(2, 1; {\Bbb Z}[\omega])$ be the
Baily-Borel compactification, here $\overline{{\Bbb B}^2}={\Bbb
B}^2 \cup \partial_{U(2, 1; {\Bbb Z}[\omega])} {\Bbb B}^2$,
$\partial_{U(2, 1; {\Bbb Z}[\omega])} {\Bbb B}^2 \subset \partial
{\Bbb B}^2$ the set of $U(2, 1; {\Bbb Z}[\omega])$-cusps. Then we
have the following commutative diagram:
$$\CD
  P @>>> C(P): Y^3=X^4+G_2(P) X^2+G_3(P) X+G_4(P) \\
  \overline{{\Bbb B}^2} @>>> \{\text{Picard curves}\}/
  \text{projective equivalence}\\
  @VV\text{quotient map}V    @VVV \\
  \overline{{\Bbb B}^2}/U(2, 1; {\Bbb Z}[\omega])
  @>=>> \widehat{{\Bbb B}^2/U(2, 1; {\Bbb Z}[\omega])}
\endCD$$
$(Y^3=X^4$ is excluded from the set of Picard curves. On ${\Bbb
B}^2-U(2, 1; {\Bbb Z}[\omega]) \cdot D \subset \overline{{\Bbb
B}^2}$, $D=\{(z, 0): |z|<1\}$ a subdisc of ${\Bbb B}^2$,
projective equivalence can be replaced by isomorphism.$)$ The
image of ${\Bbb B}^2-U(2, 1; {\Bbb Z}[\omega]) \cdot D$ is the set
of isomorphism classes of all smooth Picard curves.

\item The images of the $U(2, 1; {\Bbb Z}[\omega])$-fixed point
set of ${\Bbb B}^2$ is the set of classes of Picard curves with
larger automorphism groups $($larger than ${\Bbb Z}/3 {\Bbb Z})$.

\item The Lagrange resolvent map, which makes correspond to each
polynomial of degree four a polynomial of degree three, admits by
a change of fibres procedure over the moduli space of Picard
curves $($change Picard curves to elliptic curves$)$ a geometric
interpretation as a morphism ${\Bbb B}^2/U(2, 1; {\Bbb Z}[\omega])
\to {\Bbb H}/SL(2, {\Bbb Z})$. There is a lift along the
automorphic form $($quotient$)$ morphisms ${\Bbb B}^2 \to {\Bbb
B}^2/U(2, 1; {\Bbb Z}[\omega])$, ${\Bbb H} \to {\Bbb H}/SL(2,
{\Bbb Z})$.
\endroster}

  Moreover, Picard modular forms can be considered as theta
constants (see \cite{Ho2}). Theta functions $\vartheta
\left[\matrix a\\ b \endmatrix\right]$ with characteristics $a, b
\in {\Bbb Q}^g$ are holomorphic functions on ${\Bbb C}^g \times
{\Bbb H}_g$, ${\Bbb H}_g$ the generalized Siegel upper half plane
uniformizing the moduli space of principally polarized abelian
varieties of dimension $g$. Explicitly the theta functions
$$\vartheta \left[\matrix a\\ b \endmatrix\right]:
  {\Bbb C}^g \times {\Bbb H}_g \to {\Bbb C}$$
are defined by
$$\vartheta \left[\matrix a\\ b \endmatrix\right](z, \Omega)
 =\sum_{n \in {\Bbb Z}^{g}} \exp \left[\pi i {}^t(n+a)\Omega(n+a)
  +2 \pi i {}^t(n+a)(z+b)\right].$$
The restrictions $\vartheta|0 \times {\Bbb H}_g$,
$$\theta \left[\matrix a\\ b \endmatrix\right](\Omega)
 =\vartheta \left[\matrix a\\ b \endmatrix\right](0, \Omega)$$
are called theta constants (with characteristics $a$, $b$).

{\smc Theorem 3.15} (Feustel, Shiga) (see \cite{Ho2}). {\it Let
$\theta_i(\Omega)=\vartheta_i(0, \Omega)$, $i=0, 1, 2$, be the
theta constants on ${\Bbb H}_3$ restricting the theta functions
$$\vartheta_k(z, \Omega)=\vartheta \left[\matrix
  0 & \frac{1}{6} & 0\\
  \frac{k}{3} & \frac{1}{6} & \frac{k}{3}
  \endmatrix\right](z, \Omega), \quad
  k=0, 1, 2, \quad z \in {\Bbb C}^3.$$
Set
$$Th_1=\theta_0^3+\theta_1^3+\theta_2^3, \quad
  Th_2=-3 \theta_0^3+\theta_1^3+\theta_2^3,$$
$$Th_3=\theta_0^3-3 \theta_1^3+\theta_2^3, \quad
  Th_4=\theta_0^3+\theta_1^3-3 \theta_2^3,$$
and
$$th_i(\tau)=Th_i(*\tau), \quad i=1, 2, 3, 4, \quad \tau \in {\Bbb B}^2,$$
where the embedding $*: {\Bbb B}^2 \hookrightarrow {\Bbb H}_3$.
Then the functions $th_i(\tau)$ are the normalized Picard modular
forms satisfying
$$t_1+t_2+t_3+t_4=0, \quad \text{$t_1$, $t_2$, $t_3$ are linearly independent},$$
and the functional equations
$$\gamma^{*}(t_i)=(\det \gamma)^2 \cdot \text{sgn}(\overline{\gamma})
  \cdot j_{\gamma} \cdot t_{\overline{\gamma}(i)} \quad \text{for
  $i=1, 2, 3, 4$}.$$}

  The analogies between elliptic modular forms and Picard modular
forms are given as follows:

\roster

\item
$$\bigoplus_{m=0}^{\infty} [SL(2, {\Bbb Z}), m]={\Bbb C}[g_2, g_3],$$
$$\bigoplus_{m=0}^{\infty} [SL(2, {\Bbb Z})(2), m]={\Bbb C}[\varepsilon_1, \varepsilon_2],$$
where $g_i$ is of weight $i$, $i=2$, $3$, $\varepsilon_1$,
$\varepsilon_2$ are of weight $1$.

\item
$$\bigoplus_{m=0}^{\infty} [U(2, 1; {\Cal O}_K), m]_{\chi}={\Bbb C}[G_2, G_3, G_4],$$
$$\bigoplus_{m=0}^{\infty} [U(2, 1; {\Cal O}_K(1-\omega), m]_{\chi}={\Bbb C}[\xi_1, \xi_2, \xi_3],$$
where $G_k$ is of weight $k$, $k=2$, $3$, $4$,  $\xi_1$, $\xi_2$,
$\xi_3$ are of weight $1$.

\endroster

  This leads us to give the following conjecture, which is the
higher dimensional counterpart of Theorem 3.11 and Theorem 3.12.

{\smc Conjecture 3.16}. {\it
$${\Bbb C}[z_1, z_2, z_3]^{\text{Hessian groups}}
 ={\Bbb C}[C_6, C_9, C_{12}].\tag 3.7$$
The map
$$\phi: {\Bbb C}[C_6, C_9, C_{12}] \to {\Bbb C}[G_2, G_3, G_4]\tag 3.8$$
given by
$$\phi(h(C_6, C_9, C_{12}))=h(G_2, G_3, G_4)\tag 3.9$$
is an algebra isomorphism.}

  We have the following correspondences
$$(C_6, C_9, C_{12}) \longleftrightarrow (G_2, G_3, G_4),\tag 3.10$$
$$(G, H, K) \longleftrightarrow (\xi_1, \xi_2, \xi_3).\tag 3.11$$
$$(C_{12}, C_{18}) \longleftrightarrow (E_4, E_6)=(g_2, g_3),\tag 3.12$$
$$(H, K) \longleftrightarrow (\varepsilon_1, \varepsilon_2).\tag 3.13$$

  Note that
$$C_{12}=C_{12}(z_1, z_2, z_3) \longleftrightarrow G_4, \quad \text{three-variable function}.$$
$$C_{12}=C_{12}(H, K) \longleftrightarrow g_2, \quad \text{two-variable function}.$$

\roster

\item $GL(2)$: $\varepsilon_1$, $\varepsilon_2$: theta
      functions of weight one, elliptic modular functions

\item $GL(3)$: $\xi_1$, $\xi_2$, $\xi_3$: theta functions of
      weight one, Picard modular functions

\endroster

$$GL(2): \quad (x, y) \longleftrightarrow (\theta_0(z), \theta_1(z)), z \in {\Bbb H}.$$
$$E_4(z)=\theta_0(z)^4+8 \theta_0(z) \theta_1(z)^3, \quad
  E_6(z)=\theta_0(z)^6-20 \theta_0(z)^3 \theta_1(z)^3-8 \theta_1(z)^6.$$
$\theta_0(z)$, $\theta_1(z)$ are of weight $1$.
$$U(2, 1): \quad (z_1, z_2, z_3) \longleftrightarrow (\vartheta_0(w_1, w_2),
  \vartheta_1(w_1, w_2), \vartheta_2(w_1, w_2)), (w_1, w_2) \in {\frak S}_2,\tag 3.14$$
where ${\frak S}_2=\{ (z_1, z_2) \in {\Bbb C}^2:
z_1+\overline{z_1}>z_2 \overline{z_2} \}.$ $\vartheta_0(w_1,
w_2)$, $\vartheta_1(w_1, w_2)$ and $\vartheta_2(w_1, w_2)$ are of
weight $\frac{1}{3}$. By
$$9 C_{12}=H(H^3+8K^3), \quad -27 C_{18}=H^6-20 H^3 K^3-8 K^6,$$
we have
$$(H, K) \longleftrightarrow (\theta_0(z), \theta_1(z)),$$
$$(9 C_{12}, -27 C_{18}) \longleftrightarrow (E_4(z), E_6(z)).$$
Note that
$$H=z_1^3+z_2^3+z_3^3+6 z_1 z_2 z_3, \quad
  K=z_1^3+z_2^3+z_3^3-3 z_1 z_2 z_3.$$
Thus,
$$\theta_0(z)=\vartheta_0(w_1, w_2)^3+\vartheta_1(w_1,
  w_2)^3+\vartheta_2(w_1, w_2)^3+6 \vartheta_0(w_1, w_2)
  \vartheta_1(w_1, w_2) \vartheta_2(w_1, w_2),$$
$$\theta_1(z)=\vartheta_0(w_1, w_2)^3+\vartheta_1(w_1,
  w_2)^3+\vartheta_2(w_1, w_2)^3-3 \vartheta_0(w_1, w_2)
  \vartheta_1(w_1, w_2) \vartheta_2(w_1, w_2).$$
Hence,
$$\vartheta_0(w_1, w_2)^3+\vartheta_1(w_1, w_2)^3+\vartheta_2(w_1, w_2)^3
 =\frac{1}{3}[\theta_0(z)+2 \theta_1(z)],\tag 3.15$$
$$\vartheta_0(w_1, w_2) \vartheta_1(w_1, w_2) \vartheta_2(w_1, w_2)
 =\frac{1}{9}[\theta_0(z)-\theta_1(z)].\tag 3.16$$

  The group $PSL(2, {\Bbb Z})/\Gamma(3)$ is the tetrahedral group.
The theta functions $\theta_0(z)$, $\theta_1(z)$ define a mapping
$$\frac{\theta_0}{\theta_1}: \overline{{\Bbb H}/\Gamma(3)} \to {\Bbb P}^1({\Bbb C}).$$
This mapping is a bijection. The group $PSL(2, {\Bbb
Z})/\Gamma(3)$ acts on $\overline{{\Bbb H}/\Gamma(3)}$. Under the
above mapping this action corresponds to the action of the
tetrahedral group on a tetrahedron lying inside the Riemann sphere
${\Bbb P}^1({\Bbb C})$ (see \cite{Eb}).

  In our case, the theta functions $\vartheta_0(w_1, w_2)$, $\vartheta_1(w_1, w_2)$
and $\vartheta_2(w_1, w_2)$ define a mapping
$$\left(\frac{\vartheta_1}{\vartheta_0}, \frac{\vartheta_2}{\vartheta_0}\right):
  \overline{{\frak S}_2/U(2, 1; {\Cal O}_K)(1-\omega)} \to {\Bbb
  P}^2({\Bbb C}),$$
where $K={\Bbb Q}(\sqrt{-3})$ and ${\Cal O}_K={\Bbb Z}[\omega]$.

  Let ${\frak S}_2$ be the Siegel domain
$${\frak S}_2=\{ (z_1, z_2) \in {\Bbb C}^2: z_1+
  \overline{z_1}>|z_2|^2 \}.$$
The imaginary quadratic field $K={\Bbb Q}(\sqrt{-3})$ is called
the field of Eisenstein numbers. Its ring of integers
$${\Cal O}={\Cal O}_K=\{ a/2+b \sqrt{-3}/2: a, b \in {\Bbb Z},
  a \equiv b (\text{mod} 2) \}={\Bbb Z}+{\Bbb Z} [\omega]$$
with $\omega=\frac{-1+\sqrt{-3}}{2}$ is called the ring of
Eisenstein integers. We define the Cayley transform $C: {\Bbb
B}^2:=\{(w_1, w_2) \in {\Bbb C}^2: |w_1|^2+|w_2|^2<1 \} \to {\frak
S}_2$, where
$$C=\left(\matrix
     -\overline{\omega} &   & -\omega\\
                        & 1 &        \\
                     -1 &   & 1
    \endmatrix\right) \in GL(3, {\Cal O_K}).$$
It is known that
$$U(2, 1)=\{ g \in GL(3, {\Bbb C}): g I_{2, 1} g^{*}=I_{2, 1} \},$$
where $I_{2, 1}=\text{diag}\{ 1, 1, -1 \}$.

  Let $G=C U(2, 1) C^{-1}$, $K=C (U(2) \times U(1)) C^{-1}$ and
$G({\Cal O})=C U(2, 1; {\Cal O}_K) C^{-1}$, then $G({\Cal O})$ is
an arithmetic subgroup of $G$ and $S:=G({\Cal O}) \backslash G/K$
is called a Picard modular surface. Denote $J=\left(\matrix
            0 &  0 & 1\\
            0 & -1 & 0\\
            1 &  0 & 0
           \endmatrix\right)$, then
$C I_{2, 1} C^{*}=J$. It follows that
$$G=\{ g \in GL(3, {\Bbb C}): g J g^{*}=J \}.$$
This is the other realization of $U(2, 1)$. For simplicity, from
now on, we use the same symbol $U(2, 1)$ to denote as $G$. In
fact, $G=\{ g \in GL(3, {\Bbb C}): g^{*} J g=J \}$.

  The principal congruence subgroup $U(2, 1; {\Cal O}_K)(1-\omega)$
of $U(2, 1; {\Cal O}_K)$ is the monodromy group of the Appell
hypergeometric partial differential equations with parameters $(a,
b, b^{\prime}, c)=(\frac{1}{3}, \frac{1}{3}, \frac{1}{3},
\frac{1}{3})$. Moreover, $U(2, 1; {\Cal O}_K)(1-\omega)$ is the
modular group of the algebraic family of projective plane curves
$C$ of equation type
$$C: \quad y^3=p_4(x), \quad \deg p_4(x)=4$$
(Picard curves) endowed with a fixed order of the ramification
points of the $3$-sheeted Galois covering $C \to {\Bbb P}^1$
induced by the projection $(x, y) \mapsto x$. More precisely, the
noncompact surface ${\frak S}_2/U(2, 1; {\Cal O}_K)(1-\omega)$
classifies the isomorphy classes of these curves endowed with four
ordered points. ${\frak S}_2/U(2, 1; {\Cal O}_K)(1-\omega)$ is
isomorphic to ${\Bbb P}^2-\{\text{$4$ points}\}$. The four points
come from four cusps lying on $\partial {\frak S}_2$. $SU(2, 1;
{\Cal O}_K)(1-\omega)$ has up to $SU(2, 1; {\Cal
O}_K)(1-\omega)$-equivalence exactly four cusps and three elliptic
points $\widetilde{S_0}$, $\widetilde{S_1}$, $\widetilde{S_2} \in
{\frak S}_2$ with stationary groups of order $3$. For the smooth
compactification $\overline{{\frak S}_2/SU(2, 1; {\Cal
O}_K)(1-\omega)}$ of ${\frak S}_2/SU(2, 1; {\Cal O}_K)(1-\omega)$
one needs four elliptic curves $E_0$, $E_1$, $E_2$, $E_3$. So
${\frak S}_2/SU(2, 1; {\Cal O}_K)(1-\omega)$ or $(\overline{{\frak
S}_2/SU(2, 1; {\Cal O}_K)(1-\omega)}, E_0+E_1+E_2+E_3)$ is an
elliptic bounded surface (see \cite{Ho4}).

  According to Picard (see \cite{Pi} or \cite{Ho1}), the five
generators of the principal congruence subgroup $U(2, 1; {\Cal
O}_{K})(1-\omega)$ are given as follows:
$$g_1=\left(\matrix
  1 & \overline{\omega}-\omega & 1-\overline{\omega}\\
  0 & \overline{\omega}        & 1-\omega\\
  0 & 0                        & 1
  \endmatrix\right), \quad
  g_2=\left(\matrix
  1 & \overline{\omega}-1 & 1-\overline{\omega}\\
  0 & \overline{\omega}   & 1-\overline{\omega}\\
  0 & 0                   & 1
  \endmatrix\right), \quad
  g_3=\left(\matrix
  1 & 0      & 0\\
  0 & \omega & 0\\
  0 & 0      & 1
  \endmatrix\right),$$
$$g_4=\left(\matrix
  \omega              & 0  & 1-\overline{\omega}\\
  0                   & 1  & 0\\
  1-\overline{\omega} & 0  & -2 \omega
  \endmatrix\right), \quad
  g_5=\left(\matrix
  1 & 0 & 0\\
  \overline{\omega}-\omega & \overline{\omega} & 0\\
  1-\overline{\omega} & 1-\omega & 1
  \endmatrix\right).$$
In fact (see \cite{Y2}),
$$\left\{\aligned
  g_1 &=U_1^{-4} T_1^{-1} T_2^{-2}, \\
  g_2 &=U_1^{-4} T_1^{-2} T_2^{-1}, \\
  g_3 &=U_1^4,\\
  g_4 &=U_2^3 S^3 [T_1, T_2] S^3 (S^4 U_2)^{-1} [T_1, T_2],\\
  g_5 &=S^3 U_1^{-4} T_1^{-1} T_2 S^3.
\endaligned\right.$$
Here,
$$T_1=\left(\matrix
      1 & 1 & -\omega\\
        & 1 &       1\\
        &   &       1
      \endmatrix\right), \quad
  T_2=\left(\matrix
      1 & \omega & -\omega\\
        &      1 & \overline{\omega}\\
        &        & 1
      \endmatrix\right), \quad
  S=-\overline{\omega} J
   =\left(\matrix
                       &                   & -\overline{\omega}\\
                       & \overline{\omega} &     \\
    -\overline{\omega} &                   &
    \endmatrix\right),$$
$$U_1=\left(\matrix
      1 &         &  \\
        & -\omega &  \\
        &         & 1
      \endmatrix\right), \quad
  U_2=\left(\matrix
      -1 &         &   \\
         & -\omega &   \\
         &         & -1
      \endmatrix\right).$$
$$[T_1, T_2]:=T_1 T_2 T_1^{-1} T_2^{-1}
 =\left(\matrix
  1 & 0 & \overline{\omega}-\omega\\
    & 1 & 0\\
    &   & 1
  \endmatrix\right).$$
Note that $g_3^3=g_4^3=I$ and $g_3 g_4=g_4 g_3$. Hence,
$$\langle g_3, g_4 \rangle \cong {\Bbb Z}/3 {\Bbb Z} \times {\Bbb Z}/3 {\Bbb Z}.$$

  In \cite{Y3}, we set
$$\varphi=z_1 z_2 z_3, \quad
  \psi=z_1^3+z_2^3+z_3^3, \quad
  \chi=z_1^3 z_2^3+z_2^3 z_3^3+z_3^3 z_1^3.$$
Put
$$X(z_1, z_2, z_3)=z_1^3, \quad Y(z_1, z_2, z_3)=z_2^3,
  \quad Z(z_1, z_2, z_3)=z_3^3.\tag 3.17$$
$$Q_1(z_1, z_2, z_3)=z_1 z_2^2+z_2 z_3^2+z_3 z_1^2, \quad
  Q_2(z_1, z_2, z_3)=z_1^2 z_2+z_2^2 z_3+z_3^2 z_1.\tag 3.18$$
The functions $\varphi$, $\psi$, $\chi$, $X$, $Y$, $Z$ and $Q_1$,
$Q_2$ satisfy the relations:
$$\left\{\aligned
  \varphi^3 &=XYZ,\\
       \psi &=X+Y+Z,\\
       \chi &=XY+YZ+ZX,\\
  \chi+3 \varphi^2+\varphi \psi &=Q_1 Q_2,\\
  \psi \chi+6 \varphi \chi+6 \varphi^2 \psi+9 \varphi^3 &=Q_1^3+Q_2^3,\\
  (X-Y)(Y-Z)(Z-X) &=Q_1^3-Q_2^3.
\endaligned\right.\tag 3.19$$

  Let
$$\left\{\aligned
  W_2 &=(X+Y+Z)^2-12(XY+YZ+ZX),\\
  W_3 &=(X-Y)(Y-Z)(Z-X),\\
  {\frak W}_3 &=XYZ-\varphi^3,\\
  W_4 &=(X+Y+Z)[(X+Y+Z)^3+216 XYZ],\\
  {\frak W}_4 &=\varphi [27 XYZ-(X+Y+Z)^3],\\
  W_6 &=(X+Y+Z)^6-540 XYZ (X+Y+Z)^3-5832 X^2 Y^2 Z^2,\\
  {\frak V}_3 &=Q_1^3+Q_2^3-(X+Y+Z+6 \varphi)(XY+YZ+ZX)-6 \varphi^2 (X+Y+Z)-9 \varphi^3,\\
  {\frak V}_2 &=Q_1 Q_2-(XY+YZ+ZX)-\varphi(X+Y+Z)-3 \varphi^2,\\
  {\frak U}_3 &=Q_1^3-Q_2^3-(X-Y)(Y-Z)(Z-X).
\endaligned\right.\tag 3.20$$

{\smc Theorem 3.17}. {\it The invariants $W_2$, $W_3$, ${\frak
W}_3$, $W_4$, ${\frak W}_4$ and $W_6$ satisfy the following
algebraic relations:
$$\left\{\aligned
  W_2^3-3 W_2 W_4+2 W_6 &=432 W_3^2,\\
  8 U^3 (W_6^2-W_4^3-1728 {\frak W}_4^3) &=27 {\frak W}_3 (W_4-9 U^4)^3,
\endaligned\right.\tag 3.21$$
where
$$27 U^8-18 W_4 U^4-8 W_6 U^2-W_4^2=0.\tag 3.22$$}

{\it Proof}. By the definition of $W_2$, $W_3$, $W_4$ and $W_6$,
we have
$$W_2^3-3 W_2 W_4+2 W_6=432 W_3^2.$$
We find that
$$W_6^2-W_4^3-1728 {\frak W}_4^3=1728 {\frak W}_3 [27 XYZ-(X+Y+Z)^3]^3.\tag 3.23$$
Put
$$U=X+Y+Z, \quad V=XYZ.\tag 3.24$$
Then
$$\left\{\aligned
  U(U^3+216V) &=W_4,\\
  U^6-540 U^3 V-5832 V^2 &=W_6.
\endaligned\right.\tag 3.25$$
By the first equation in (3.25), we have
$$V=\frac{W_4-U^4}{216 U}.\tag 3.26$$
Substituting the above expression (3.26) in the second equation in
(3.25), we obtain that
$$27 U^8-18 W_4 U^4-8 W_6 U^2-W_4^2=0.$$
Set $t=U^2$. Let
$$t^4-\frac{2}{3} W_4 t^2-\frac{8}{27} W_6 t-\frac{1}{27} W_4^2
 =(t^2+kt+l)(t^2-kt+m).\tag 3.27$$
Then
$$\left\{\aligned
  l+m-k^2 &=-\frac{2}{3} W_4,\\
  k(m-l) &=-\frac{8}{27} W_6,\\
  lm &=-\frac{1}{27} W_4^2.
\endaligned\right.\tag 3.28$$
By the first two equations in (3.28), we have
$$\left\{\aligned
  l &=\frac{1}{2k}(k^3-\frac{2}{3} W_4 k+\frac{8}{27} W_6),\\
  m &=\frac{1}{2k}(k^3-\frac{2}{3} W_4 k-\frac{8}{27} W_6).
\endaligned\right.\tag 3.29$$
Substituting these two expressions in the third equation in
(3.28), we have
$$k^6-\frac{4}{3} W_4 k^4+\frac{16}{27} W_4^2 k^2-\frac{64}{729} W_6^2=0,\tag 3.30$$
i.e.,
$$(k^2-\frac{4}{9} W_4)^3=\frac{64}{729}(W_6^2-W_4^3).$$
Hence,
$$k=\pm \frac{2}{3} \sqrt{W_4+\root 3 \of{W_6^2-W_4^3}}.\tag 3.31$$
Consequently,
$$l=-\frac{1}{9} W_4+\frac{2}{9} \root 3 \of{W_6^2-W_4^3}
    \pm \frac{2}{9} \frac{W_6}{\sqrt{W_4+\root 3 \of{W_6^2-W_4^3}}}.\tag 3.32$$
$$m=-\frac{1}{9} W_4+\frac{2}{9} \root 3 \of{W_6^2-W_4^3}
    \mp \frac{2}{9} \frac{W_6}{\sqrt{W_4+\root 3 \of{W_6^2-W_4^3}}}.\tag 3.33$$
Thus,
$$t=\frac{1}{2}(-k \pm \sqrt{k^2-4l}), \quad \text{or} \quad
  t=\frac{1}{2}(k \pm \sqrt{k^2-4m}).\tag 3.34$$
We have
$$W_6^2-W_4^3-1728 {\frak W}_4^3
 =1728 {\frak W}_3 (27V-U^3)^3
 =\frac{27}{8} {\frak W}_3 \frac{(W_4-9 U^4)^3}{U^3}.$$
Therefore,
$$8 U^3 (W_6^2-W_4^3-1728 {\frak W}_4^3)=27 {\frak W}_3 (W_4-9 U^4)^3.$$
$\qquad \qquad \qquad \qquad \qquad \qquad \qquad \qquad \qquad
 \qquad \qquad \qquad \qquad \qquad \qquad \qquad \qquad \qquad
 \quad \boxed{}$

{\smc Theorem 3.18}. {\it The invariants $W_2$, $W_3$, ${\frak
W}_3$, $W_4$, ${\frak W}_4$, ${\frak U}_3$, ${\frak V}_2$  and
${\frak V}_3$ satisfy the following differential relations:
$$\frac{\partial(W_2, W_3, W_4, {\frak W}_3)}{\partial(X, Y, Z, \varphi)}
 =288 {\frak W}_4^2.\tag 3.35$$
$$\frac{\partial({\frak V}_3, {\frak V}_2)}{\partial(Q_1, Q_2)}=3({\frak U}_3+W_3).\tag 3.36$$}

{\it Proof}. We have
$$\left\{\aligned
  \frac{\partial W_2}{\partial X} &=2X-10Y-10Z,\\
  \frac{\partial W_2}{\partial Y} &=2Y-10Z-10X,\\
  \frac{\partial W_2}{\partial Z} &=2Z-10X-10Y,\\
  \frac{\partial W_2}{\partial \varphi} &=0.
  \endaligned\right.$$
$$\left\{\aligned
  \frac{\partial W_3}{\partial X} &=(Y-Z)(Y+Z-2X),\\
  \frac{\partial W_3}{\partial Y} &=(Z-X)(Z+X-2Y),\\
  \frac{\partial W_3}{\partial Z} &=(X-Y)(X+Y-2Z),\\
  \frac{\partial W_3}{\partial \varphi} &=0.
  \endaligned\right.$$
$$\left\{\aligned
  \frac{\partial W_4}{\partial X} &=4(X+Y+Z)^3+216 XYZ+216 (X+Y+Z)YZ,\\
  \frac{\partial W_4}{\partial Y} &=4(X+Y+Z)^3+216 XYZ+216 (X+Y+Z)ZX,\\
  \frac{\partial W_4}{\partial Z} &=4(X+Y+Z)^3+216 XYZ+216 (X+Y+Z)XY,\\
  \frac{\partial W_4}{\partial \varphi} &=0.
  \endaligned\right.$$
$$\left\{\aligned
  \frac{\partial {\frak W}_3}{\partial X} &=YZ,\\
  \frac{\partial {\frak W}_3}{\partial Y} &=ZX,\\
  \frac{\partial {\frak W}_3}{\partial Z} &=XY,\\
  \frac{\partial {\frak W}_3}{\partial \varphi} &=-3 \varphi^2.
  \endaligned\right.$$
Hence,
$$\frac{\partial(W_2, W_3, W_4, {\frak W}_3)}{\partial(X, Y, Z, \varphi)}
 =-24 \varphi^2 \Delta,$$
where
$$\aligned
  \Delta
=&[(X+Y+Z)^3+54 YZ(2X+Y+Z)] \times\\
 &\times [-(X+Y+Z)^3-9X(X+Y+Z)^2+27 X^2(Y+Z)+27 X(Y^2+Z^2)]+\\
+&[(X+Y+Z)^3+54 ZX(2Y+Z+X)] \times\\
 &\times [-(X+Y+Z)^3-9Y(X+Y+Z)^2+27 Y^2(Z+X)+27 Y(Z^2+X^2)]+\\
+&[(X+Y+Z)^3+54 XY(2Z+X+Y)] \times\\
 &\times [-(X+Y+Z)^3-9Z(X+Y+Z)^2+27 Z^2(X+Y)+27 Z(X^2+Y^2)].
\endaligned$$
By the proof of Proposition 4.3 in \cite{Y3}, we have
$$\Delta=-12 [27 XYZ-(X+Y+Z)^3]^2.$$
Thus,
$$\frac{\partial(W_2, W_3, W_4, {\frak W}_3)}{\partial(X, Y, Z, \varphi)}
 =288 \varphi^2 [27 XYZ-(X+Y+Z)^3]^2=288 {\frak W}_4^2.$$
On the other hand, we find that
$$\frac{\partial({\frak V}_3, {\frak V}_2)}{\partial(Q_1, Q_2)}
 =3(Q_1^3-Q_2^3)=3({\frak U}_3+W_3).$$
$\qquad \qquad \qquad \qquad \qquad \qquad \qquad \qquad \qquad
 \qquad \qquad \qquad \qquad \qquad \qquad \qquad \qquad \qquad
 \quad \boxed{}$

  Combine Theorem 3.17 with Theorem 3.18, we get the proof of Theorem 1.4 (Main
Theorem 2).

  According to \cite{Ho1}, we denote the elementary symmetric functions of
$X_1$, $X_2$, $X_3$, $X_4$ by
$$\aligned
  \Sigma_1 &=-(X_1+X_2+X_3+X_4),\\
  \Sigma_2 &=X_1 X_2+X_1 X_3+X_1 X_4+X_2 X_3+X_2 X_4+X_3 X_4,\\
  \Sigma_3 &=-(X_1 X_2 X_3+X_1 X_2 X_4+X_1 X_3 X_4+X_2 X_3 X_4),\\
  \Sigma_4 &=X_1 X_2 X_3 X_4,
\endaligned$$
and
$$\Delta=\prod_{1 \leq i<j \leq 4} (X_i-X_j).$$
Then
$$27 \Delta^2=4 (\Sigma_2^2-3 \Sigma_1 \Sigma_3+12 \Sigma_4)^3
  -(27 \Sigma_1^2 \Sigma_4+27 \Sigma_3^2+2 \Sigma_2^3-72 \Sigma_2
  \Sigma_4+9 \Sigma_1 \Sigma_2 \Sigma_3)^2.$$
Their images in ${\Bbb C}[x_1, x_2, x_3, x_4]={\Bbb C}[X_1, X_2,
X_3, X_4]/(X_1+X_2+X_3+X_4)$ are denoted by $G_1=0$, $G_2$, $G_3$
and $G_4$, respectively. Hence,
$$\Delta^2=16 G_2^4 G_4-128 G_2^2 G_4^2-4 G_2^3 G_3^2+144 G_2 G_3^2 G_4-27 G_3^4+256 G_4^3.$$

{\smc Proposition 3.19} (see \cite{Ho1}). {\it There are two
isomorphisms of graded rings:
$$\bigoplus_{m=0}^{\infty} [U(2, 1; {\Cal O}_K), m]_{\chi}
  \cong {\Bbb C} [G_2, G_3, G_4]$$
and
$$\bigoplus_{m=0}^{\infty} [SU(2, 1; {\Cal O}_K), m]_{\chi}
  \cong {\Bbb C} [G_2, G_3, G_4, z^2].$$
The fundamental relation between the generators is
$z^6=\Delta^2$.}

  We find the following correspondence:
$$(G_2, G_3, G_4, z^2) \longleftrightarrow (W_2, W_3, W_4, {\frak W}_4).$$
Here, $W_2$, $W_3$, $W_4$ and ${\frak W}_4$ satisfy the single
relation:
$$6912 {\frak W}_4^3=W_2^6+9 W_2^2 W_4^2+432^2 W_3^4-4 W_4^3
  -6 W_2^4 W_4-864 W_2^3 W_3^2+2592 W_2 W_3^2 W_4.$$

  The equations $W_4=0$ and ${\frak W}_4=0$ define two surfaces
of order $4$. The intersection of the cubic surface ${\frak
W}_3=0$ with the quartic surface $W_4=0$ gives four elliptic
curves
$$\left\{\aligned
  \varphi^3+6 \varphi XY+(X+Y)XY &=0, \\
  \varphi^3+6 \omega \varphi XY+(X+Y)XY &=0, \\
  \varphi^3+6 \overline{\omega} \varphi XY+(X+Y)XY &=0, \\
  \varphi^3+(X+Y)XY &=0.
\endaligned\right.$$
The intersection of the cubic surface ${\frak W}_3=0$ with the
quartic surface ${\frak W}_4=0$ gives the other four elliptic
curves
$$\left\{\aligned
  \varphi^3-3 \varphi XY+(X+Y)XY &=0, \\
  \varphi^3-3 \omega \varphi XY+(X+Y)XY &=0, \\
  \varphi^3-3 \overline{\omega} \varphi XY+(X+Y)XY &=0, \\
  XYZ &=0.
\endaligned\right.$$
The equation $W_2=0$ defines a quadric surface, the equation
$W_3=0$ is the union of three planes. The intersection of $W_3=0$
with ${\frak W}_3=0$ gives three elliptic curves
$$Y^2 Z-\varphi^3=0, \quad Z^2 X-\varphi^3=0, \quad X^2 Y-\varphi^3=0.$$
The intersection of $W_2=0$ with ${\frak W}_3=0$ gives an
algebraic curve of degree six
$$\varphi^6-10 XY(X+Y) \varphi^3+(X^4 Y^2-10 X^3 Y^3+X^2 Y^4)=0.$$

  It is known that the Hessian group of order $216$ is generated
by the five generators (see \cite{Y3}):
$$A=\left(\matrix
  0 & 1 & 0\\
  0 & 0 & 1\\
  1 & 0 & 0
  \endmatrix\right), \quad
  B=\left(\matrix
  1 & 0 & 0\\
  0 & 0 & 1\\
  0 & 1 & 0
  \endmatrix\right),$$
$$C=\left(\matrix
  1 & 0 & 0\\
  0 & \omega & 0\\
  0 & 0 & \omega^2
  \endmatrix\right), \quad
  D=\left(\matrix
  1 & 0 & 0\\
  0 & \omega & 0\\
  0 & 0 & \omega
  \endmatrix\right),$$
$$E=\frac{1}{\sqrt{-3}} \left(\matrix
  1 & 1 & 1\\
  1 & \omega & \omega^2\\
  1 & \omega^2 & \omega
  \endmatrix\right).$$
We have
$$\left\{\aligned
  \varphi(A(z_1, z_2, z_3)) &=\varphi,\\
        X(A(z_1, z_2, z_3)) &=Y,\\
        Y(A(z_1, z_2, z_3)) &=Z,\\
        Z(A(z_1, z_2, z_3)) &=X,\\
      Q_1(A(z_1, z_2, z_3)) &=Q_1,\\
      Q_2(A(z_1, z_2, z_3)) &=Q_2.
\endaligned\right.$$
$$\left\{\aligned
  \varphi(B(z_1, z_2, z_3)) &=\varphi,\\
        X(B(z_1, z_2, z_3)) &=X,\\
        Y(B(z_1, z_2, z_3)) &=Z,\\
        Z(B(z_1, z_2, z_3)) &=Y,\\
      Q_1(B(z_1, z_2, z_3)) &=Q_2,\\
      Q_2(B(z_1, z_2, z_3)) &=Q_1.
\endaligned\right.$$
$$\left\{\aligned
  \varphi(C(z_1, z_2, z_3)) &=\varphi,\\
        X(C(z_1, z_2, z_3)) &=X,\\
        Y(C(z_1, z_2, z_3)) &=Y,\\
        Z(C(z_1, z_2, z_3)) &=Z,\\
      Q_1(C(z_1, z_2, z_3)) &=\overline{\omega} Q_1,\\
      Q_2(C(z_1, z_2, z_3)) &=\omega Q_2.
\endaligned\right.$$
$$\left\{\aligned
  \varphi(D(z_1, z_2, z_3)) &=\overline{\omega} \varphi,\\
        X(D(z_1, z_2, z_3)) &=X,\\
        Y(D(z_1, z_2, z_3)) &=Y,\\
        Z(D(z_1, z_2, z_3)) &=Z,\\
      Q_1(D(z_1, z_2, z_3)) &=\overline{\omega} z_1 z_2^2+z_2 z_3^2+\omega z_3 z_1^2,\\
      Q_2(D(z_1, z_2, z_3)) &=\omega z_1^2 z_2+z_2^2 z_3+\overline{\omega} z_3^2 z_1.
\endaligned\right.$$
$$\left\{\aligned
  (\sqrt{-3})^3 \varphi(E(z_1, z_2, z_3)) &=X+Y+Z-3 \varphi,\\
  (\sqrt{-3})^3 X(E(z_1, z_2, z_3)) &=X+Y+Z+6 \varphi+3 Q_1+3 Q_2,\\
  (\sqrt{-3})^3 Y(E(z_1, z_2, z_3)) &=X+Y+Z+6 \varphi+3 \overline{\omega} Q_1+3 \omega Q_2,\\
  (\sqrt{-3})^3 Z(E(z_1, z_2, z_3)) &=X+Y+Z+6 \varphi+3 \omega Q_1+3 \overline{\omega} Q_2,\\
  (\sqrt{-3})^3 Q_1(E(z_1, z_2, z_3)) &=3X+3 \overline{\omega} Y+3 \omega Z,\\
  (\sqrt{-3})^3 Q_2(E(z_1, z_2, z_3)) &=3X+3 \omega Y+3 \overline{\omega} Z.
\endaligned\right.$$

  In the form of matrices, we have
$$\aligned
 &A(X, Y, Z, \varphi, Q_1, Q_2)\\
=&\left(\matrix
   & 1 &   &   &   &   \\
   &   & 1 &   &   &   \\
 1 &   &   &   &   &   \\
   &   &   & 1 &   &   \\
   &   &   &   & 1 &   \\
   &   &   &   &   & 1
 \endmatrix\right)
 \left(\matrix
   X\\
   Y\\
   Z\\
   \varphi\\
   Q_1\\
   Q_2
 \endmatrix\right),
\endaligned\tag 3.37$$
$$\aligned
 &B(X, Y, Z, \varphi, Q_1, Q_2)\\
=&\left(\matrix
 1 &   &   &   &   &   \\
   &   & 1 &   &   &   \\
   & 1 &   &   &   &   \\
   &   &   & 1 &   &   \\
   &   &   &   &   & 1 \\
   &   &   &   & 1 &
 \endmatrix\right)
 \left(\matrix
   X\\
   Y\\
   Z\\
   \varphi\\
   Q_1\\
   Q_2
 \endmatrix\right),
\endaligned\tag 3.38$$
$$\aligned
 &C(X, Y, Z, \varphi, Q_1, Q_2)\\
=&\left(\matrix
 1 &   &   &   &   &   \\
   & 1 &   &   &   &   \\
   &   & 1 &   &   &   \\
   &   &   & 1 &   &   \\
   &   &   &   & \overline{\omega} &   \\
   &   &   &   &   & \omega
 \endmatrix\right)
 \left(\matrix
   X\\
   Y\\
   Z\\
   \varphi\\
   Q_1\\
   Q_2
 \endmatrix\right),
\endaligned\tag 3.39$$
$$\aligned
 &(\sqrt{-3})^3 E(X, Y, Z, \varphi, Q_1, Q_2)\\
=&\left(\matrix
 1 & 1 & 1 & 6 & 3 & 3 \\
 1 & 1 & 1 & 6 & 3 \overline{\omega} & 3 \omega\\
 1 & 1 & 1 & 6 & 3 \omega & 3 \overline{\omega}\\
 1 & 1 & 1 & -3 & 0 & 0\\
 3 & 3 \overline{\omega} & 3 \omega & 0 & 0 & 0\\
 3 & 3 \omega & 3 \overline{\omega} & 0 & 0 & 0
 \endmatrix\right)
 \left(\matrix
   X\\
   Y\\
   Z\\
   \varphi\\
   Q_1\\
   Q_2
 \endmatrix\right),
\endaligned\tag 3.40$$

  We find that
$$A^3=B^2=C^3=E^4=I.\tag 3.41$$

{\smc Theorem 3.19} (see \cite{Y3}, Main Theorem 4). {\it The
invariants $G$, $H$ and $K$ satisfy the following algebraic
equations, which are the form-theoretic resolvents $($algebraic
resolvents$)$ of $G$, $H$, $K$$:$
$$\left\{\aligned
  4 G^3+H^2 G-C_6 G-4 C_9 &=0,\\
  H (H^3+8 K^3)-9 C_{12} &=0,\\
  K (K^3-H^3)-27 {\frak C}_{12} &=0.
\endaligned\right.$$}
Consequently,
$$C_{18}=-\frac{1}{27}(H^6-20 H^3 K^3-8 K^6).$$

  Note that
$$W_2=C_6, W_3=C_9, W_4=C_{12}, {\frak W}_4={\frak C}_{12}, W_6=C_{18}.$$
Set $(x, y, z)=(G, H, K)$. Then
$$x=Q_1-Q_2, \quad y=X+Y+Z+6 \varphi, \quad z=X+Y+Z-3 \varphi.\tag 3.42$$
We have
$$\frac{1}{9}y(y^3+8z^3)=(X+Y+Z)[(X+Y+Z)^3+216 \varphi^3]=W_4-216 (X+Y+Z) {\frak W}_3.\tag 3.43$$
$$\frac{1}{27}z(z^3-y^3)=\varphi[27 \varphi^3-(X+Y+Z)^3]={\frak W}_4-27 \varphi {\frak W}_3.\tag 3.44$$
Thus,
$$X+Y+Z=\frac{9 W_4-y(y^3+8z^3)}{1944 {\frak W}_3}.$$
$$\varphi=\frac{27 {\frak W}_4-z(z^3-y^3)}{729 {\frak W}_3}.$$
$$XYZ={\frak W}_3+\varphi^3={\frak W}_3+\frac{[27 {\frak W}_4-z(z^3-y^3)]^3}
      {729^3 {\frak W}_3^3}.$$
Note that $U=X+Y+Z$ and $V=XYZ$ satisfy the equations:
$$U(U^3+216V)=W_4, \quad U^6-540 U^3 V-5832 V^2=W_6.$$
We obtain the following relations:
$$\aligned
 &[9W_4-y(y^3+8z^3)]^4+2^{12} [9W_4-y(y^3+8z^3)] \{3^{18} {\frak W}_3^4
  +[27 {\frak W}_4-z(z^3-y^3)]^3\}\\
=&2^{12} \cdot 3^{20} {\frak W}_3^4 W_4.
\endaligned$$
$$\aligned
 &[9W_4-y(y^3+8z^3)]^6-2^{11} \cdot 5 [9W_4-y(y^3+8z^3)]^3 \{3^{18} {\frak W}_3^4
  +[27 {\frak W}_4-z(z^3-y^3)]^3\}\\
 &-2^{21} \{3^{18} {\frak W}_3^4+[27 {\frak W}_4-z(z^3-y^3)]^3\}^2
  =2^{18} \cdot 3^{30} {\frak W}_3^6 W_6.
\endaligned$$

  Similarly, we have
$$4 x^3+x y^2-W_2 x-4 W_3=4 {\frak U}_3-12 (Q_1-Q_2) {\frak V}_2.$$
Hence,
$$Q_1-Q_2=\frac{4 {\frak U}_3-4 x^3-x y^2+W_2 x+4 W_3}{12 {\frak V}_2}.$$
By the identity:
$$(Q_1-Q_2)^3=Q_1^3-Q_2^3-3 Q_1 Q_2 (Q_1-Q_2),$$
where
$$Q_1^3-Q_2^3={\frak U}_3+W_3,$$
$$Q_1 Q_2={\frak V}_2+(XY+YZ+ZX)+\varphi (X+Y+Z)+3 \varphi^2.$$
Note that
$$\aligned
 &XY+YZ+ZX=\frac{1}{12}[(X+Y+Z)^2-W_2]\\
=&\frac{1}{2^8 \cdot 3^{11} {\frak W}_3^2}
  [(9 W_4-y (y^3+8 z^3))^2-2^6 \cdot 3^{10} W_2 {\frak W}_3^2].
\endaligned$$
We obtain the following relation:
$$\aligned
 &2^4 \cdot 3^8 {\frak W}_3^2 (4 {\frak U}_3-4 x^3-x y^2+W_2 x+4 W_3)^3\\
=&2^{10} \cdot 3^{11} {\frak W}_3^2 {\frak V}_2^3 ({\frak U}_3+W_3)
  -(4{\frak U}_3-4 x^3-x y^2+W_2 x+4 W_3) {\frak V}_2^2 \times\\
 &\times \{2^8 \cdot 3^{11} {\frak W}_3^2 {\frak V}_2+[(9 W_4-y
 (y^ 3+8 z^3))^2-2^6 \cdot 3^{10} W_2 {\frak W}_3^2]+\\
 &+2^5 [27 {\frak W}_4-z (z^3-y^3)][9 W_4-y (y^3+8 z^3)]+2^8
  [27 {\frak W}_4-z (z^3-y^3)]^2\}.
\endaligned$$
Therefore, we finish the proof of Theorem 1.5 (Main Theorem 3).
\flushpar $\qquad \qquad \qquad \qquad \qquad \qquad \qquad
\qquad\qquad \qquad \qquad \qquad \qquad \qquad \qquad \qquad
\qquad \qquad \quad \boxed{}$

{\smc Proposition 3.20}. {\it The space curves $W_2={\frak
W}_3={\frak V}_3={\frak V}_2=0$, $W_3={\frak W}_3={\frak
V}_3={\frak V}_2=0$, $W_4={\frak W}_3={\frak V}_3={\frak V}_2=0$,
$W_6={\frak W}_3={\frak V}_3={\frak V}_2=0$ and ${\frak
W}_4={\frak W}_3={\frak V}_3={\frak V}_2=0$ are invariant curves
on the invariant surface ${\frak W}_3={\frak V}_3={\frak V}_2=0$
under the action of the subgroup of Hessian group generated by
$A$, $B$, $C$ and $E$.}

{\it Proof}. Note that on the surface ${\frak W}_3={\frak
V}_3={\frak V}_2=0$, $XYZ=\varphi^3$,
$$Q_1 Q_2=XY+YZ+ZX+\varphi(X+Y+Z)+3 \varphi^2,$$
$$Q_1^3+Q_2^3=(X+Y+Z+6 \varphi)(XY+YZ+ZX)+6 \varphi^2 (X+Y+Z)+9 \varphi^3.$$
Hence,
$$\aligned
 &(\sqrt{-3})^9 {\frak W}_3(E(X, Y, Z, \varphi, Q_1, Q_2))\\
=&(X+Y+Z+6 \varphi)^3+27 Q_1^3+27 Q_2^3-27 (X+Y+Z+6 \varphi) Q_1 Q_2+\\
 &-(X+Y+Z-3 \varphi)^3\\
=&0.
\endaligned$$
Similarly,
$$\aligned
 &(\sqrt{-3})^9 {\frak V}_3(E(X, Y, Z, \varphi, Q_1, Q_2))\\
=&54 (X^3+Y^3+Z^3+6XYZ)-81(X^2 Y+Y^2 Z+Z^2 X+X Y^2+Y Z^2+Z X^2)+\\
 &-9(X+Y+Z)[3(X+Y+Z+6 \varphi)^2-27 Q_1 Q_2]+\\
 &-18 (X+Y+Z-3 \varphi)^2 (X+Y+Z+6 \varphi)-9(X+Y+Z-3 \varphi)^3\\
=&243 (X+Y+Z) [Q_1 Q_2-(X+Y+Z) \varphi-3 \varphi^2]-729 \varphi^3+\\
 &-243 (X^2 Y+Y^2 Z+Z^2 X+X Y^2+Y Z^2+Z X^2)\\
=&729 (XYZ-\varphi^3)\\
=&0.
\endaligned$$
$$\aligned
 &(\sqrt{-3})^6 {\frak V}_2(E(X, Y, Z, \varphi, Q_1, Q_2))\\
=&9(X^2+Y^2+Z^2-XY-YZ-ZX)-3(X+Y+Z+6 \varphi)^2+27 Q_1 Q_2+\\
 &-3(X+Y+Z+6 \varphi)(X+Y+Z-3 \varphi)-3(X+Y+Z-3 \varphi)^2\\
=&27[Q_1 Q_2-(XY+YZ+ZX)-\varphi (X+Y+Z)-3 \varphi^2]\\
=&0.
\endaligned$$
On the other hand, we have
$$\left\{\aligned
  {\frak W}_3(A(X, Y, Z, \varphi, Q_1, Q_2)) &={\frak W}_3,\\
  {\frak W}_3(B(X, Y, Z, \varphi, Q_1, Q_2)) &={\frak W}_3,\\
  {\frak W}_3(C(X, Y, Z, \varphi, Q_1, Q_2)) &={\frak W}_3.
\endaligned\right.$$
$$\left\{\aligned
  {\frak V}_3(A(X, Y, Z, \varphi, Q_1, Q_2)) &={\frak V}_3,\\
  {\frak V}_3(B(X, Y, Z, \varphi, Q_1, Q_2)) &={\frak V}_3,\\
  {\frak V}_3(C(X, Y, Z, \varphi, Q_1, Q_2)) &={\frak V}_3.
\endaligned\right.$$
$$\left\{\aligned
  {\frak V}_2(A(X, Y, Z, \varphi, Q_1, Q_2)) &={\frak V}_2,\\
  {\frak V}_2(B(X, Y, Z, \varphi, Q_1, Q_2)) &={\frak V}_2,\\
  {\frak V}_2(C(X, Y, Z, \varphi, Q_1, Q_2)) &={\frak V}_2.
\endaligned\right.$$
Hence, ${\frak W}_3={\frak V}_3={\frak V}_2=0$ is an invariant
surface. On this surface, we have
$$Q_1^3-Q_2^3=\pm (X-Y)(Y-Z)(Z-X).$$
Hence,
$$(\sqrt{-3})^9 W_3(E(X, Y, Z, \varphi, Q_1, Q_2))
 =81 \sqrt{-3} (Q_1^3-Q_2^3)=\pm 81 \sqrt{-3} W_3.$$
Furthermore,
$$\left\{\aligned
  W_3(A(X, Y, Z, \varphi, Q_1, Q_2)) &=(Y-Z)(Z-X)(X-Y)=W_3,\\
  W_3(B(X, Y, Z, \varphi, Q_1, Q_2)) &=(X-Z)(Z-Y)(Y-X)=-W_3,\\
  W_3(C(X, Y, Z, \varphi, Q_1, Q_2)) &=(X-Y)(Y-Z)(Z-X)=W_3.
\endaligned\right.$$
This implies that $W_3={\frak W}_3={\frak V}_3={\frak V}_2=0$ is
an invariant curve on the above surface.

  We have
$$\aligned
 &(\sqrt{-3})^6 W_2(E(X, Y, Z, \varphi, Q_1, Q_2))\\
=&-27 (X+Y+Z+6 \varphi)^2+324 Q_1 Q_2\\
=&-27[(X+Y+Z)^2-12(XY+YZ+ZX)]\\
=&-27 W_2.
\endaligned$$
Moreover,
$$\left\{\aligned
  W_2(A(X, Y, Z, \varphi, Q_1, Q_2))&=W_2,\\
  W_2(B(X, Y, Z, \varphi, Q_1, Q_2))&=W_2,\\
  W_2(C(X, Y, Z, \varphi, Q_1, Q_2))&=W_2.
\endaligned\right.$$
This implies that $W_2={\frak W}_3={\frak V}_3={\frak V}_2=0$ is
an invariant curve on the above surface.

  We have
$$\aligned
 &(\sqrt{-3})^{12} W_4(E(X, Y, Z, \varphi, Q_1, Q_2))\\
=&729(X+Y+Z+6 \varphi)[(X+Y+Z+6 \varphi)^3+24(Q_1^3+Q_2^3)+\\
 &-24(X+Y+Z+6 \varphi) Q_1 Q_2]\\
=&729(X+Y+Z+6 \varphi)(X+Y+Z)[(X+Y+Z)^2-6(X+Y+Z) \varphi+36 \varphi^2]\\
=&729(X+Y+Z)[(X+Y+Z)^3+216 \varphi^3]\\
=&729(X+Y+Z)[(X+Y+Z)^3+216 XYZ]\\
=&729 W_4.
\endaligned$$
Moreover,
$$\left\{\aligned
  W_4(A(X, Y, Z, \varphi, Q_1, Q_2))&=W_4,\\
  W_4(B(X, Y, Z, \varphi, Q_1, Q_2))&=W_4,\\
  W_4(C(X, Y, Z, \varphi, Q_1, Q_2))&=W_4.
\endaligned\right.$$
This implies that $W_4={\frak W}_3={\frak V}_3={\frak V}_2=0$ is
an invariant curve on the above surface.

  On the above invariant surface one has
$$W_2^3-3 W_2 W_4+2 W_6=432 W_3^2, \quad
  W_6^2-W_4^3=1728 {\frak W}_4^3.$$
Thus, $W_6={\frak W}_3={\frak V}_3={\frak V}_2=0$ and ${\frak
W}_4={\frak W}_3={\frak V}_3={\frak V}_2=0$ are invariant curves
on the above surface.

\flushpar $\qquad \qquad \qquad \qquad \qquad \qquad \qquad \qquad
\qquad \qquad \qquad \qquad \qquad \qquad \qquad \qquad \qquad
\qquad \quad \boxed{}$

\vskip 0.5 cm

\centerline{\bf 4. Hessian polyhedra and Galois representations
                   associated with cubic surfaces}

\vskip 0.5 cm

  Manin's book \cite{Ma} (see also \cite{MT}) is devoted to the
arithmetic of smooth cubic surfaces. We give an account of some of
the results of this book that are characteristic of cubic surfaces
only. A cubic surface $X$ has a very specific property: every
point $x \in X$ generates a birational automorphism $t_x$ of $X$
by reflection in $x$. This enables us to describe the group of
birational automorphisms of $X$. On the other hand, it allows us
to introduce on $X(k)$ the structure of the ternary relation of
collinearity of points. We recall that a symmetric quasigroup $C$
is a set with a binary composition $(x, y) \mapsto x \circ y$ such
that the ternary relation $x \circ y=z$ is symmetric under all
permutations of $x$, $y$ and $z$, and such that $x \circ y=y \circ
x$ and $x \circ (x \circ y)=y$ for any pair $x, y \in C$. Suppose
that an equivalence relation $S$ is given on $X(k)$, we say that
$S$ is admissible if collinearity of points induces the structure
of a symmetric quasigroup on $X(k)/S$. We say that a symmetric
quasigroup $C$ is abelian if for any element $u$ the composition
law $xy=u \circ (x \circ y)$ turns $C$ into an abelian group. If
any three elements of $C$ generate an abelian subquasigroup, we
say that $C$ is a CH-quasigroup. We fix an element $u$ in $C$ and
set $xy=u \circ (x \circ y)$, this composition law defines on $C$
the structure of a commutative Moufang loop (CML), that is, it is
commutative, has a unit element $u$ and inverse, and moreover
satisfies the three weak associativity properties: $x(xy)=x^2 y$,
$(xy)(xz)=x^2 (yz)$ and $x(y(xz))=(x^2 y) z$.

{\smc Theorem 4.1} (see \cite{Ma}, \cite{MT}). {\it Let $S$ be an
admissible equivalence relation on a smooth cubic surface $X$.
Then $X(k)/S$ is a CH-quasigroup, and in the corresponding CML the
relation $x^6=1$ is satisfied identically; this CML is finite if
and only if it is finitely generated, and in this case its order
is $2^a \cdot 3^b$.}

  In our case, let $H=\langle A, B, C, E \rangle$. Note that
the order of $H$ is $2^3 \cdot 3^2$.

  Let $k$ be a perfect field, and let $V$ run through the smooth
cubic surfaces over $k$. The classes of the lines on $V \otimes
\overline{k}$ generate the group $N(V)=\text{Pic}(V \otimes
\overline{k})$ and the action of the Galois group
$G=\text{Gal}(\overline{k}/k)$ on $N(V)$ preserves symmetry and it
implicitly contains an extremely large amount of information on
the arithmetic and geometry of $V$.

  Let $k$ be a global field. Let $\chi$ denote the character of the
representation of $G$ on $N(V)$, and let $L(s, \chi, k)$ be the
Artin $L$-function of the field $k$ corresponding to the character
$\chi$.

{\smc Theorem 4.2} (Weil) (see \cite{Ma}). {\it The Hasse-Weil
zeta function of the surface $V$ coincides $($up to a finite
number of Euler factors$)$ with the product
$$\zeta(s, k) \zeta(s-2, k) L(s-1, \chi, k).$$}

{\it Definition} 4.3 (see \cite{Ma}). Let $r \geq 1$ be an
integer. We consider a composed object $\{ N_r, \omega_r, ( ,
)\}$, where

\roster

\item $N_r={\Bbb Z}^{r+1}=\bigoplus_{i=0}^{r} {\Bbb Z} l_i$;
$(l_i)$ a chosen basis.

\item $\omega_r=(-3, 1, \cdots, 1) \in N_r$.

\item $( , )$ is a bilinear form $N_r \times N_r \to {\Bbb Z}$
given by the formula
$$(l_0, l_0)=1, \quad (l_i, l_i)=-1 \quad \text{if $i \geq 1$},
  \quad (l_i, l_j)=0 \quad \text{if $i \neq j$}.$$

\item $R_r :=\{ l \in N_r: \quad (l, \omega_r)=0, \quad (l, l)=-2
\}$.

\item $I_r :=\{ l \in N_r: \quad (l, \omega_r)=(l, l)=-1 \}$.

\endroster

{\smc Theorem 4.4} (see \cite{Ma}). {\it Let $V$ be a smooth cubic
surface over a perfect field $k$.

\roster

\item The triple $\{ N(V), \omega_V, \text{intersection number}
\}$ is isomorphic to $\{ N_6, \omega_6, ( , ) \}$ described in
Definition $4.3$.

\item The Galois group $G$ acts on $N(V)$ and preserves $\omega_V$
and the intersection number.

\item The set of classes of lines on $V \otimes \overline{k}$ goes
over into $I_r$ under the isomorphism $(1)$.

\endroster}

{\smc Theorem 4.5} (see \cite{Ma}). {\it Let $3 \leq r \leq 8$.

\roster

\item The scalar product $($with opposite sign$)$ on ${\Bbb R}
\otimes_{{\Bbb Z}} N_r \cong {\Bbb R}^{r+1}$ induces on the
orthogonal complement of $\omega_r$ the structure of a Euclidean
space. The set $R_r$ is a root system in it of type $A_1 \times
A_2$, $A_4$, $D_5$, $E_6$, $E_7$, $E_8$, respectively $($the sum
of the indices is equal to $r)$.

\item The following groups coincide:

\item"(a)" the group of automorphisms of the lattice $N_r$
preserving $\omega_r$ and the scalar product;

\item"(b)" the group of permutations of the vectors from $I_r$
preserving their pairwise scalar products;

\item"(c)" the Weyl group $W(R_r)$ of the system $R_r$ generated
by the reflections with respect to the roots.

\endroster}

  In fact, Fano varieties in dimension two are called Del Pezzo
surfaces. Over an algebraically closed field, these are: ${\Bbb
P}^2$, ${\Bbb P}^1 \times {\Bbb P}^1$ and $S_d$, where $S_d$ is
the blow up of ${\Bbb P}^2$ at $9-d$ points, and the degree $d=1,
\cdots, 8$. However, starting with $r=9$, the system $R_r$ becomes
infinite. Let $V$ be a Del Pezzo surface of degree $1 \leq d \leq
7$, and let $f: V \to {\Bbb P}^2$ be its representation in the
form of a monoidal transformation of the plane with as centre the
union of $r=9-d$ points $x_1, \cdots, x_r$. It is shown that (see
\cite{Ma}) the map $D \mapsto$ (class of ${\Cal O}_V(D)$) $\in
\text{Pic}(V)$ establishes a one-one onto correspondence between
exceptional curves on $V$ and exceptional classes in the Picard
group. These classes generate the Picard group. When $r=6$, the
number of exceptional curves on $V$ is $27$.

  The significance of cubic surfaces comes from the following: the
higher the degree of a rational surface, the smaller the rank of
its Picard module, and therefore the less freedom there is for the
action of the Galois group. Del Pezzo surfaces of degree three are
cubic surfaces in ${\Bbb P}^3$ and the Weyl group $W(E_6)$ has
enough large order. We have the following theorem (see \cite{MT}):
Any minimal rational surface $X$ over $k$ of degree $d \leq 4$ is
$k$-birationally nontrivial. The triple $\{ N(V), \omega_V,
\text{intersection number on $N(V)$}\}$ can be constructed for
every surface $V$. The action of the Galois group and the group of
automorphisms of $V$ on this triple yields essential invariants of
the surface.

  Note that when $r=6$, we have $|I_6|=27$, $|R_6|=72$, $|W(E_6)|=2^7 \cdot 3^4 \cdot 5$,
$|S_6|=720$, $|W(E_6)/S_6|=72$. On the other hand, in our case,
for $H=\langle A, B, C, E \rangle$, we have $|H|=72$.

  According to \cite{Hu}, the modern point of view is to consider
cubic surfaces as del Pezzo surfaces of degree $3$. The
combinatorics of the $27$ lines is then encoded in the Picard
group of the del Pezzo surface. In fact, the complement in
$\text{Pic}(S)$ to the hyperplane section, call it
$\text{Pic}^{0}(S)$, is isomorphic to the root lattice of $E_6$.
The equation of the surface is given by the embedding of $S$ by
means of the linear system of elliptic curves through the six
given points.

  Fix six points in ${\Bbb P}^2$, say $x=(p_1, \cdots, p_6)$, such
that the $p_i$ are in general position, i.e., no three lie on a
line, and not all six lie on a conic. Let $\widehat{{\Bbb P}_x^2}$
denote the blow up of ${\Bbb P}^2$ at all six points, $\varrho_x:
\widehat{{\Bbb P}_x^2} \to {\Bbb P}^2$. Consider the following
curves as classes in $\text{Pic}(\widehat{{\Bbb P}_x^2})$:

\roster

\item $a_1, \cdots, a_6$, the exceptional divisors over $(p_1,
      \cdots, p_6)$;

\item $b_1, \cdots, b_6$, $b_i$ the proper transform of the conic
      $q_i$ passing through all points $p_j$, $j \neq i$;

\item $c_{ik}$, the proper transform of the line $\overline{p_i
      p_k}$.

\endroster
If we consider the surface $\widehat{{\Bbb P}_x^2}$, we have
$H^2(\widehat{{\Bbb P}_x^2}, {\Bbb Z})=[l] {\Bbb Z} \oplus_i {\Bbb
Z} a_i$. Let $Q$ be the intersection form on $H^2(\widehat{{\Bbb
P}_x^2}, {\Bbb Z})$, then the classes $a_i$, $b_i$, $c_{ij}$
fulfill $Q(a_i, a_i)=Q(b_i, b_i)=Q(c_{ij}, c_{ij})=-1$. In a
well-known manner one takes a rank $6$ subset, which is isomorphic
to the root lattice of type $E_6$.

  Consider the orthocomplement of the canonical class on
$\widehat{{\Bbb P}_x^2}$, and denote this by
$\text{Pic}^0(\widehat{{\Bbb P}_x^2})$. Recall that the
anti-canonical class is $3l+\sum_{i=1}^{6} a_i$, and that the
anti-canonical embedding of $\widehat{{\Bbb P}_x^2}$ is as a cubic
surface. Consequently we may view $\text{Pic}^0(\widehat{{\Bbb
P}_x^2})$ as the orthocomplement of the hyperplane section class
of $\text{Pic}(S_x)$, where $S_x$ is the cubic surface which is
the anti-canonical embedding. The following elements $\lambda$
with $Q(\lambda, \lambda)=-2$ form a basis of
$\text{Pic}^{0}(S_x)$:
$$\aligned
  \alpha_0 &=l-a_1-a_2-a_3,\\
  \alpha_1 &=a_1-a_2,\\
  \alpha_2 &=a_2-a_3,\\
  \alpha_3 &=a_3-a_4,\\
  \alpha_4 &=a_4-a_5,\\
  \alpha_5 &=a_5-a_6.
\endaligned$$
These also form a base of a root system of type $E_6$, by taking
$\alpha_1, \cdots, \alpha_5$ as the sub-root system of type $A_5$.
Since the classes $a_i$, $b_i$, $c_{ij}$ are exceptional, they all
represent elements of $\text{Pic}^{0}(S_x)$.

  According to \cite{Ma}, we have the realization problems.
\roster

\item Which subgroups of $W(E_6)$ are realized as an image of the
      Galois group by its representation on $\text{Pic}(V \otimes
      \overline{k})$ for some $V$?

\item Which extensions of the field $k$ correspond with kernels of
      such representations?

\endroster

  Shioda (see \cite{Shi2} and \cite{Shi3}) considered Galois
representations and algebraic equations arising from Mordell-Weil
lattices. Let $E=E_{\lambda}$ be the elliptic curve, defined by
one of the following equations, over $K_0=k_0(t)$ where $k_0={\Bbb
Q}(\lambda)={\Bbb Q}(p_i, q_j)$.
$$\aligned
  (E_8) \quad &y^2=x^3+x \left(\sum_{i=0}^{3} p_i t^i\right)+\left(\sum_{i=0}^{3} q_i t^i+t^5\right)\\
  &\lambda=(p_0, p_1, p_2, p_3, q_0, q_1, q_2, q_3) \in {\Bbb A}^8
\endaligned$$
$$\aligned
  (E_7) \quad &y^2=x^3+x (p_0+p_1 t+t^3)+\left(\sum_{i=0}^{4} q_i t^i\right)\\
  &\lambda=(p_0, p_1, q_0, q_1, q_2, q_3, q_4) \in {\Bbb A}^7
\endaligned$$
$$\aligned
  (E_6) \quad &y^2=x^3+x \left(\sum_{i=0}^{2} p_i t^i\right)+\left(\sum_{i=0}^{2} q_i t^i+t^4\right)\\
  &\lambda=(p_0, p_1, p_2, q_0, q_1, q_2) \in {\Bbb A}^6
\endaligned$$
Let $f: S_{\lambda} \to {\Bbb P}^1$ be the associated elliptic
surface. The fibre $f^{-1}(\infty)$ is an additive singular fibre
of type II, III or IV according to the case $(E_r)$ for $r=8, 7$
or $6$. Assume that $\lambda$ satisfies the condition:
$$\text{every fibre of $f$ over $t \neq \infty$ is irreducible.}$$
Then the Mordell-Weil lattice $E(k(t))$ is isomorphic to
$E_r^{*}$, the dual lattice of $E_r$, where $k$ is the algebraic
closure of $k_0$, and we get the Galois representation
$$\varrho=\varrho_{\lambda}: \text{Gal}(k/k_0) \to \text{Aut}(E_r^{*})=\text{Aut}(E_r).$$
Recall that $\text{Aut}(E_r)=W(E_r)$ for $r=8$ or $7$, and $W(E_r)
\{ \pm 1 \}$ for $r=6$, where $W(E_r)$ denotes the Weyl group of
type $E_r$. In any case, we have $\text{Im}(\varrho_{\lambda})
\subset W(E_r)$.

  We assume that $\lambda$ is generic, i.e. $p_i$, $q_j$ are algebraically
independent over ${\Bbb Q}$. Then the above condition holds.

{\smc Theorem 4.6} (Shioda). {\it Let $\lambda$ be generic over
${\Bbb Q}$. Then

$(i)$ the image of the Galois representation $\varrho_{\lambda}$
is the full Weyl group $W(E_r)$:
$\text{Im}(\varrho_{\lambda})=W(E_r)$.

$(ii)$ Let ${\Cal K}_{\lambda}/k_0$ be the Galois extension
corresponding to $\text{Ker}(\varrho_{\lambda})$, and let $\{ P_1,
\cdots, P_r \}$ be a basis of $E(k(t)) \simeq E_r^{*}$ consisting
of minimal vectors. Further let $u_i=sp_{\infty}^{\prime}(P_i) \in
{\Cal K}_{\lambda} \subset k$. Then $u_1, \cdots, u_r$ are
algebraically independent over ${\Bbb Q}$, and we have
$${\Cal K}_{\lambda}=k_0(u_1, \cdots, u_r)={\Bbb Q}(u_1, \cdots, u_r),$$
$$\text{Gal}({\Bbb Q}(u_1, \cdots, u_r)/{\Bbb Q}(p_i, q_j))=W(E_r).$$

$(iii)$ $W(E_r)$ acts on the vector space ${\Bbb Q} u_1 \oplus
\cdots \oplus {\Bbb Q} u_r$ and hence on the polynomial ring
${\Bbb Q}[u_1, \cdots, u_r]$, and the ring of the invariants is:
$${\Bbb Q}[u_1, \cdots, u_r]^{W(E_r)}={\Bbb Q}[p_i, q_j].$$
In particular, $p_i$ and $q_j$ form the fundamental invariants of
the Weyl group $W(E_r)$ and we can explicitly write
$$\text{$p_i$ or $q_j=J_d(u_1, \cdots, u_r)$,}$$
where $J_d$ denotes a $W(E_r)$-invariant of degree $d$, $d \in \{
2, 5, 6, 8, 9, 12 \}$, $\{ 2, 6, 8, 10, 12,$ $14, 18 \}$ or $\{ 2,
8, 12, 14, 18, 20, 24, 30 \}$ for $r=6$, $7$ or $8$.}

{\smc Theorem 4.7} (Shioda). {\it Every Galois extension of ${\Bbb
Q}$ with Galois group $W(E_r)$ is obtained as ${\Cal K}_{\lambda}$
for some $\lambda \in \Lambda$. In other words, every such
extension arises from the Mordell-Weil lattice of the elliptic
curve $E_{\lambda}/{\Bbb Q}(t)$ for some $\lambda \in {\Bbb
Q}^r$.}

  In the equation $(E_6)$, let $y=y^{\prime} \pm t^2$ and $(x: y^{\prime}:
t: 1)=(X: Y: Z: W)$. Then we have a cubic surface
$V=V_{\lambda}^{\pm}$ in ${\Bbb P}^3$:
$$Y^2 W \pm 2 Y Z^2=X^3+X(p_0 W^2+p_1 ZW+p_2 Z^2)+q_0 W^3+q_1
  ZW^2+q_2 Z^2 W.$$
It is smooth if and only if the above condition about $\lambda$
holds. Under this assumption, the narrow Mordell-Weil lattice
$E_{\lambda}(k(t))^{0} \simeq E_6$ is isomorphic to the primitive
part of $NS(V)$. The $27$ minimal sections of the form $P=(at+b,
t^2+dt+e)$ in the Mordell-Weil lattice $E_{\lambda}(k(t))$ are
transformed into the $27$ lines on $V=V_{\lambda}^{+}$ defined by
the equation
$$X=aZ+bW, \quad Y=dZ+eW.$$
(Similarly, the $27$ sections $-P$ are mapped to the $27$ lines on
$V_{\lambda}^{-}$.)

  Therefore, Shioda answered the problems raised by Manin when the
image is surjective, i.e. $\text{Im}(\varrho)=W(E_6)$.

  The Galois representations arising from Mordell-Weil lattices
$$\varrho: G=\text{Gal}(k/k_0) \to \text{Aut}(E(k(C)), \langle , \rangle)=\text{a finite group},$$
where $k(C)$ is the function field of a smooth projective curve
$C$ over an algebraically closed field $k$, are quite different
from those arising from the torsion points of an elliptic curve or
an abelian variety (e.g. the Tate modules), because we are dealing
with points of infinite order here.

  In fact, Mordell-Weil lattices $E_6$, $E_7$ and $E_8$ come from Klein
singularities $E_6$, $E_7$ and $E_8$ which are intimately related
to tetrahedron, octahedron and icosahedron. In contrast to
Shioda's approach, we study Galois representations and the
algebraic equations arising from Hessian polyhedra. It should be
pointed out that Shioda used elliptic cubic surfaces, which are
not the cubic surfaces of general type. However, we will use the
general type cubic surfaces of type $F_2$.

  For $(z_1, z_2, z_3) \in {\Bbb P}^2$, $X$, $Y$, $Z$, $\varphi$,
$Q_1$ and $Q_2$ are cubic polynomials of $z_1$, $z_2$ and $z_3$.
We will study the invariants $x$, $y$ and $z$ as functions of $X$,
$Y$, $Z$, $\varphi$, $Q_1$ and $Q_2$.

  We have
$$A(x)=x, \quad B(x)=-x, \quad B^2(x)=x.$$
$$C(x)=\overline{\omega} Q_1-\omega Q_2, \quad
  C^2(x)=\omega Q_1-\overline{\omega} Q_2, \quad C^3(x)=x.$$
$$E(x)=Y-Z, \quad E^2(x)=x.$$
$$EC(x)=X-Y, \quad EC^2(x)=Z-X.$$
$$E^2 C(x)=\omega Q_1-\overline{\omega} Q_2, \quad
  E^2 C^2(x)=\overline{\omega} Q_1-\omega Q_2.$$
$$A(y)=y, \quad B(y)=y, \quad C(y)=y.$$
$$E(y)=\frac{\sqrt{-3}}{3}(y+2z), \quad E^2(y)=-y, \quad
  E^3(y)=-\frac{\sqrt{-3}}{3}(y+2z), \quad E^4(y)=y.$$
$$A(z)=z, \quad B(z)=z, \quad C(z)=z.$$
$$E(z)=\frac{\sqrt{-3}}{3}(y-z), \quad E^2(z)=-z, \quad
  E^3(z)=-\frac{\sqrt{-3}}{3}(y-z), \quad E^4(z)=z.$$
Put
$$g_1=E(x), \quad g_2=EC(x), \quad g_3=EC^2(x), \quad
  g_4=x, \quad g_5=C(x), \quad g_6=C^2(x).\tag 4.1$$
We find that
$$g_1+g_2+g_3=0, \quad g_4+g_5+g_6=0.\tag 4.2$$
$$g_1^3+g_2^3+g_3^3=3W_3=3 C_9, \quad
  g_4^3+g_5^3+g_6^3=3(Q_1^3-Q_2^3)=3 C_9.\tag 4.3$$
Hence, we obtain the following variety:
$$S: \left\{\aligned
  g_1+g_2+g_3 &=0,\\
  g_4+g_5+g_6 &=0,\\
  g_1^3+g_2^3+g_3^3 &=g_4^3+g_5^3+g_6^3,
\endaligned\right.\tag 4.4$$
which is a cubic surface in ${\Bbb P}^5$. Moreover, we get the
following map:
$$\aligned
  {\Bbb P}^2 &\to {\Bbb P}^5\\
  (z_1, z_2, z_3) &\mapsto (g_1, g_2, g_3, g_4, g_5, g_6)
\endaligned$$
by which the projective plane ${\Bbb C} {\Bbb P}^2$ is mapped
biholomorphically to the cubic surface $S \subset {\Bbb C} {\Bbb
P}^5$.

  We find that
$$\left\{\aligned
  E(g_1) &=g_4,\\
  E(g_2) &=g_6,\\
  E(g_3) &=g_5,
\endaligned\right. \quad
  \left\{\aligned
  E(g_4) &=g_1,\\
  E(g_5) &=g_2,\\
  E(g_6) &=g_3.
\endaligned\right.\tag 4.5$$
$$\left\{\aligned
  A(g_1) &=g_3,\\
  A(g_2) &=g_1,\\
  A(g_3) &=g_2,
\endaligned\right. \quad
  \left\{\aligned
  A(g_4) &=g_4,\\
  A(g_5) &=g_5,\\
  A(g_6) &=g_6.
\endaligned\right.\tag 4.6$$
$$\left\{\aligned
  B(g_1) &=-g_1,\\
  B(g_2) &=-g_3,\\
  B(g_3) &=-g_2,
\endaligned\right. \quad
  \left\{\aligned
  B(g_4) &=-g_4,\\
  B(g_5) &=-g_6,\\
  B(g_6) &=-g_5.
\endaligned\right.\tag 4.7$$
$$\left\{\aligned
  C(g_1) &=g_1,\\
  C(g_2) &=g_2,\\
  C(g_3) &=g_3,
\endaligned\right. \quad
  \left\{\aligned
  C(g_4) &=g_5,\\
  C(g_5) &=g_6,\\
  C(g_6) &=g_4.
\endaligned\right.\tag 4.8$$
Hence, the cubic surface $S$ is invariant under the action of the
subgroup generated by $A$, $B$, $C$ and $E$.

  There are $27$ lines on the cubic surface $S$:
$$l_1: \quad g_1=g_2+g_3=g_4=g_5+g_6=0,$$
$$l_2: \quad g_1=g_2+g_3=g_5=g_4+g_6=0,$$
$$l_3: \quad g_1=g_2+g_3=g_6=g_4+g_5=0,$$
$$l_4: \quad g_2=g_1+g_3=g_4=g_5+g_6=0,$$
$$l_5: \quad g_2=g_1+g_3=g_5=g_4+g_6=0,$$
$$l_6: \quad g_2=g_1+g_3=g_6=g_4+g_5=0,$$
$$l_7: \quad g_3=g_1+g_2=g_4=g_5+g_6=0,$$
$$l_8: \quad g_3=g_1+g_2=g_5=g_4+g_6=0,$$
$$l_9: \quad g_3=g_1+g_2=g_6=g_4+g_5=0,$$
and
$$l_{j, 1}: \quad g_1=\omega^j g_4, \quad g_2=\omega^j g_5, \quad g_3=\omega^j g_6,
            \quad g_1+g_2+g_3=0,$$
$$l_{j, 2}: \quad g_1=\omega^j g_4, \quad g_2=\omega^j g_6, \quad g_3=\omega^j g_5,
            \quad g_1+g_2+g_3=0,$$
$$l_{j, 3}: \quad g_1=\omega^j g_5, \quad g_2=\omega^j g_4, \quad g_3=\omega^j g_6,
            \quad g_1+g_2+g_3=0,$$
$$l_{j, 4}: \quad g_1=\omega^j g_5, \quad g_2=\omega^j g_6, \quad g_3=\omega^j g_4,
            \quad g_1+g_2+g_3=0,$$
$$l_{j, 5}: \quad g_1=\omega^j g_6, \quad g_2=\omega^j g_4, \quad g_3=\omega^j g_5,
            \quad g_1+g_2+g_3=0,$$
$$l_{j, 6}: \quad g_1=\omega^j g_6, \quad g_2=\omega^j g_5, \quad g_3=\omega^j g_4,
            \quad g_1+g_2+g_3=0,$$
where $j \equiv 0, 1, 2$(mod $3$).

  We find that
$$E(l_1)=l_1, \quad A(l_1)=l_7, \quad B(l_1)=l_1, \quad C(l_1)=l_2.$$
$$E(l_2)=l_4, \quad A(l_2)=l_8, \quad B(l_2)=l_3, \quad C(l_2)=l_3.$$
$$E(l_3)=l_7, \quad A(l_3)=l_9, \quad B(l_3)=l_2, \quad C(l_3)=l_1.$$
$$E(l_4)=l_3, \quad A(l_4)=l_1, \quad B(l_4)=l_7, \quad C(l_4)=l_5.$$
$$E(l_5)=l_6, \quad A(l_5)=l_2, \quad B(l_5)=l_9, \quad C(l_5)=l_6.$$
$$E(l_6)=l_9, \quad A(l_6)=l_3, \quad B(l_6)=l_8, \quad C(l_6)=l_4.$$
$$E(l_7)=l_2, \quad A(l_7)=l_4, \quad B(l_7)=l_4, \quad C(l_7)=l_8.$$
$$E(l_8)=l_5, \quad A(l_8)=l_5, \quad B(l_8)=l_6, \quad C(l_8)=l_9.$$
$$E(l_9)=l_8, \quad A(l_9)=l_6, \quad B(l_9)=l_5, \quad C(l_9)=l_7.$$
$$E(l_{j, 1})=l_{2j, 2}, \quad A(l_{j, 1})=l_{j, 4}, \quad
  B(l_{j, 1})=l_{j, 1}, \quad C(l_{j, 1})=l_{j, 4}.$$
$$E(l_{j, 2})=l_{2j, 1}, \quad A(l_{j, 2})=l_{j, 6}, \quad
  B(l_{j, 2})=l_{j, 2}, \quad C(l_{j, 2})=l_{j, 3}.$$
$$E(l_{j, 3})=l_{2j, 5}, \quad A(l_{j, 3})=l_{j, 2}, \quad
  B(l_{j, 3})=l_{j, 6}, \quad C(l_{j, 3})=l_{j, 6}.$$
$$E(l_{j, 4})=l_{2j, 3}, \quad A(l_{j, 4})=l_{j, 5}, \quad
  B(l_{j, 4})=l_{j, 5}, \quad C(l_{j, 4})=l_{j, 5}.$$
$$E(l_{j, 5})=l_{2j, 6}, \quad A(l_{j, 5})=l_{j, 1}, \quad
  B(l_{j, 5})=l_{j, 4}, \quad C(l_{j, 5})=l_{j, 1}.$$
$$E(l_{j, 6})=l_{2j, 4}, \quad A(l_{j, 6})=l_{j, 3}, \quad
  B(l_{j, 6})=l_{j, 3}, \quad C(l_{j, 6})=l_{j, 2}.$$

  It is well-known that the Segre cubic threefold is given by
(see \cite{Hu})
$${\Cal S}_3: \quad \left\{\aligned
  x_0+x_1+x_2+x_3+x_4+x_5 &=0,\\
  x_0^3+x_1^3+x_2^3+x_3^3+x_4^3+x_5^3 &=0.
\endaligned\right.$$
The hyperplane sections $\{ x_i=0 \}$ of ${\Cal S}_3$ are Clebsch
diagonal cubic surfaces with equation
$$S_3: \quad \left\{\aligned
  x_0+x_1+x_2+x_3+x_4 &=0,\\
  x_0^3+x_1^3+x_2^3+x_3^3+x_4^3 &=0.
\endaligned\right.$$
The relation between $S_3$ and the icosahedral group was studied
by Hirzebruch. It turns out that $S_3$ is $A_5$-equivariantly
birational to the Hilbert modular surface for ${\Cal O}_K$ of
level $\sqrt{5}$ where $K={\Bbb Q}(\sqrt{5})$ (see \cite{Hu}).

  We find that the hyperplane sections $\{ x_i+x_j+x_k=0 \}$ of
${\Cal S}_3$ are our cubic surfaces with equation
$$S: \quad \left\{\aligned
  g_1+g_2+g_3 &=0,\\
  (-g_4)+(-g_5)+(-g_6) &=0,\\
  g_1^3+g_2^3+g_3^3+(-g_4)^3+(-g_5)^3+(-g_6)^3 &=0.
\endaligned\right.$$

  Our cubic surface can be expressed as
$$S: \quad x_1^2 x_2+x_1 x_2^2+x_3^2 x_4+x_3 x_4^2=0.\tag 4.9$$
The Hessian variety of $S$ (up to a constant) is given by
$$(x_1^2+x_1 x_2+x_2^2)(x_3^2+x_3 x_4+x_4^2)=0,\tag 4.10$$
which are four planes.

  Let us consider the group of permutations (see \cite{Hu}),
$\text{Aut}({\Cal L})$, of the $27$ lines, by which we mean the
permutations of the lines preserving the intersection behavior of
the lines. For this it is useful to consider the famous double
sixes and the notation for the $27$ lines introduced by
Schl\"{a}fli. A double six is an array
$$N=\left[\matrix
    a_1 & a_2 & a_3 & a_4 & a_5 & a_6\\
    b_1 & b_2 & b_3 & b_4 & b_5 & b_6
   \endmatrix\right]$$
of $12$ of the $27$ lines with the property that two of these $12$
meet if and only if they are in different rows and columns. (This
notation distinguishes this particular set of $12$, although any
such double six is equivalent to it under $\text{Aut}({\Cal L})$.)
The other lines are given by the $\left(\matrix 6\\ 2
\endmatrix\right)=15$ $c_{ij}=a_i b_j \cap a_j b_i$, where $a_i
b_j$ denotes the tritangent spanned by those two lines. There are
$36$ double sixes, namely the $N$ above, $15$ $N_{ij}$ and $20$
$N_{ijk}$:
$$N_{ij}=\left[\matrix
         a_i & b_i & c_{jk} & c_{jl} & c_{jm} & c_{jn}\\
         a_j & b_j & c_{ik} & c_{il} & c_{im} & c_{in}
        \endmatrix\right],$$
$$N_{ijk}=\left[\matrix
          a_i & a_j & a_k & c_{mn} & c_{ln} & c_{lm}\\
          c_{jk} & c_{ik} & c_{ij} & b_l & b_m & b_n
         \endmatrix\right].$$
Since a double six describes the intersection behavior of the
lines, we see that $S_6$ (the symmetric group on six letters) acts
by permutations on a double six and a ${\Bbb Z}_2$ acts by
exchanging rows. Since there are $36$ double sixes, we have
$|\text{Aut}({\Cal L})|=|S_6| \cdot 2 \cdot 36=51840$. In fact
$\text{Aut}({\Cal L})$ is the Weyl group $W(E_6)$, and the $36$
double sixes correspond to the positive roots. The $27$ lines
correspond to the $27$ fundamental weights of $E_6$, and many
other sets of objects (lines, tritangents, etc.) correspond to
natural sets of objects (roots, weights, etc.) of $E_6$.

  According to \cite{Seg}, a nonsingular cubic surface can only be
of one of the five types: $F_1$, $F_2$, $F_3$, $F_4$ and $F_5$.
Our cubic surface is of $F_2$ which has $15$ real lines and $12$
complex lines of the 2nd kind.

  The $15$ real lines of a cubic surface $F_2$ are of two different
sorts: $6$ of them, $l_{0, 1}$, $\cdots$, $l_{0, 6}$, constituting
$2$ complementary triplets, are elliptic; the $9$ others $l_1$,
$\cdots$, $l_9$ are hyperbolic of the 1st kind and constitute a
Steiner set. The remaining lines $l_{1, 1}$, $\cdots$, $l_{1, 6}$,
$l_{2, 1}$, $\cdots$, $l_{2, 6}$ of $F_2$ form a self-conjugate
double-six, which we call `of the 9th kind', in which each pair of
corresponding lines is a pair of conjugate complex lines of the
2nd kind; hence the $\sigma$-transformation inherent to it changes
each line of $F_2$ in its conjugate.

  $F_2$ has $15$ other self-conjugate double-sixes, which are the
$15$ double-sixes permutable with the one just considered. Each
sextuplet of such a double-six is consequently self-conjugate, and
consists of $4$ real lines and $2$ conjugate complex lines of the
2nd kind; but, while in $9$ of these $15$ double-sixes (which we
call `of the 4th kind') the $4$ real lines of each sextuplet are
$2$ elliptic and $2$ hyperbolic, in the $6$ others (which we call
`of the 5th kind', and which constitute $2$ permutable triads of
associate double-sixes) the $4$ real lines of each sextuplet are
$1$ elliptic and $3$ hyperbolic. If $\omega$ is any one of the
$15$ double-sixes considered, its $8$ real lines are the only real
lines of $F_2$ skew to a well-determined real line of $F_2$ which
is hyperbolic or elliptic according as $\omega$ is of the 4th or
5th kind.

  In fact, the Clebsch diagonal cubic surface is of $F_1$ which has
its $27$ lines all real: $12$ of them are elliptic, and $15$ are
hyperbolic of the 1st kind.

  Let $\Gamma_i$ be the group of the lines of a real cubic surface
$F_i$ and ${\frak G}$ be the group of the $27$ lines of a cubic
surface.

  Segre \cite{Seg} proved that the group $\Gamma_1$ inherent to a
cubic surface $F_1$ is icosahedral, and can be defined as the
group of the $60$ substitutions of ${\frak G}$ which transform
into themselves each of the $2$ sextuplets of elliptic lines and
also the system of $5$ principal tritangent planes of $F_1$,
inducing among these planes a substitution of even class.

  The group $\Gamma_2$ inherent to a cubic surface $F_2$ is of
order $36$; it is simply isomorphic with the direct product of $2$
symmetric groups of degree $3$, its transformations being uniquely
defined by the property of performing an arbitrary substitution
among the $3$ elliptic left-handed lines of $F_2$ and an arbitrary
substitution among the $3$ elliptic right-handed lines of $F_2$.
The group operates transitively among the $9$ pairs of
complementary triplets of the 1st kind of $F_2$, any $2$ such
pairs being transformed one into the other by $4$ transformations
of $\Gamma_2$; in particular, the subgroup of $\Gamma_2$
transforming into itself one of these pairs is trirectangular, and
coincides with the subgroup of $\Gamma_2$ which leaves unchanged
each of the $2$ elliptic lines of the pair.

  In our case, we find that the set of $15$ real lines, i.e.,
$\{ l_1, \cdots, l_9 \}$ and $\{ l_{0, 1}, \cdots, l_{0, 6} \}$
are invariant under the actions of $A$, $B$, $C$ and $E$,
respectively. For the complex lines, we have
$$E(l_{1, 1})=l_{2, 2}, E(l_{1, 2})=l_{2, 1}, E(l_{1, 3})=l_{2, 5},
  E(l_{1, 4})=l_{2, 3}, E(l_{1, 5})=l_{2, 6}, E(l_{1, 6})=l_{2, 4},$$
$$E(l_{2, 1})=l_{1, 2}, E(l_{2, 2})=l_{1, 1}, E(l_{2, 3})=l_{1, 5},
  E(l_{2, 4})=l_{1, 3}, E(l_{2, 5})=l_{1, 6}, E(l_{2, 6})=l_{1, 4}.$$

  We find that
$$\left\{\aligned
  l_{1, 1} \cap l_{1, 2}=\overline{l_{2, 1} \cap l_{2, 2}} &=(-2 \omega, \omega, \omega, -2, 1, 1),\\
  l_{1, 1} \cap l_{1, 3}=\overline{l_{2, 1} \cap l_{2, 3}} &=(\omega, \omega, -2 \omega, 1, 1, -2),\\
  l_{1, 1} \cap l_{1, 4}=\overline{l_{2, 1} \cap l_{2, 4}} &=\emptyset,\\
  l_{1, 1} \cap l_{1, 5}=\overline{l_{2, 1} \cap l_{2, 5}} &=\emptyset,\\
  l_{1, 1} \cap l_{1, 6}=\overline{l_{2, 1} \cap l_{2, 6}} &=(\omega, -2 \omega, \omega, 1, -2, 1).
\endaligned\right.$$
$$\left\{\aligned
  l_{1, 2} \cap l_{1, 3}=\overline{l_{2, 2} \cap l_{2, 3}} &=\emptyset,\\
  l_{1, 2} \cap l_{1, 4}=\overline{l_{2, 2} \cap l_{2, 4}} &=(\omega, -2 \omega, \omega, 1, 1, -2),\\
  l_{1, 2} \cap l_{1, 5}=\overline{l_{2, 2} \cap l_{2, 5}} &=(\omega, \omega, -2 \omega, 1, -2, 1),\\
  l_{1, 2} \cap l_{1, 6}=\overline{l_{2, 2} \cap l_{2, 6}} &=\emptyset.
\endaligned\right.$$
$$\left\{\aligned
  l_{1, 3} \cap l_{1, 4}=\overline{l_{2, 3} \cap l_{2, 4}} &=(-2 \omega, \omega, \omega, 1, -2, 1),\\
  l_{1, 3} \cap l_{1, 5}=\overline{l_{2, 3} \cap l_{2, 5}} &=(\omega, -2 \omega, \omega, -2, 1, 1),\\
  l_{1, 3} \cap l_{1, 6}=\overline{l_{2, 3} \cap l_{2, 6}} &=\emptyset.
\endaligned\right.$$
$$\left\{\aligned
  l_{1, 4} \cap l_{1, 5}=\overline{l_{2, 4} \cap l_{2, 5}} &=\emptyset,\\
  l_{1, 4} \cap l_{1, 6}=\overline{l_{2, 4} \cap l_{2, 6}} &=(\omega, \omega, -2 \omega, -2, 1, 1).
\endaligned\right.$$
$$l_{1, 5} \cap l_{1, 6}=\overline{l_{2, 5} \cap l_{2, 6}}=(-2 \omega, \omega, \omega, 1, 1, -2).$$
$$\left\{\aligned
  l_{1, 1} \cap l_{2, 1} &=\emptyset,\\
  l_{1, 1} \cap l_{2, 2} &=\emptyset,\\
  l_{1, 1} \cap l_{2, 3} &=\emptyset,\\
  l_{1, 1} \cap l_{2, 4} &=(1, \overline{\omega}, \omega, \overline{\omega}, \omega, 1),\\
  l_{1, 1} \cap l_{2, 5} &=(\overline{\omega}, 1, \omega, \omega, \overline{\omega}, 1),\\
  l_{1, 1} \cap l_{2, 6} &=\emptyset.
\endaligned\right.$$
$$\left\{\aligned
  l_{1, 2} \cap l_{2, 1} &=\emptyset,\\
  l_{1, 2} \cap l_{2, 2} &=\emptyset,\\
  l_{1, 2} \cap l_{2, 3} &=(1, \omega, \overline{\omega}, \overline{\omega}, \omega, 1),\\
  l_{1, 2} \cap l_{2, 4} &=\emptyset,\\
  l_{1, 2} \cap l_{2, 5} &=\emptyset,\\
  l_{1, 2} \cap l_{2, 6} &=(\overline{\omega}, \omega, 1, \omega, \overline{\omega}, 1).
\endaligned\right.$$
$$\left\{\aligned
  l_{1, 3} \cap l_{2, 1} &=\emptyset,\\
  l_{1, 3} \cap l_{2, 2} &=(1, \overline{\omega}, \omega, \omega, \overline{\omega}, 1),\\
  l_{1, 3} \cap l_{2, 3} &=\emptyset,\\
  l_{1, 3} \cap l_{2, 4} &=\emptyset,\\
  l_{1, 3} \cap l_{2, 5} &=\emptyset,\\
  l_{1, 3} \cap l_{2, 6} &=(\overline{\omega}, 1, \omega, \overline{\omega}, \omega, 1).
\endaligned\right.$$
$$\left\{\aligned
  l_{1, 4} \cap l_{2, 1} &=(1, \omega, \overline{\omega}, \omega, \overline{\omega}, 1),\\
  l_{1, 4} \cap l_{2, 2} &=\emptyset,\\
  l_{1, 4} \cap l_{2, 3} &=\emptyset,\\
  l_{1, 4} \cap l_{2, 4} &=\emptyset,\\
  l_{1, 4} \cap l_{2, 5} &=(\overline{\omega}, \omega, 1, \overline{\omega}, \omega, 1),\\
  l_{1, 4} \cap l_{2, 6} &=\emptyset.
\endaligned\right.$$
$$\left\{\aligned
  l_{1, 5} \cap l_{2, 1} &=(\omega, 1, \overline{\omega}, \overline{\omega}, \omega, 1),\\
  l_{1, 5} \cap l_{2, 2} &=\emptyset,\\
  l_{1, 5} \cap l_{2, 3} &=\emptyset,\\
  l_{1, 5} \cap l_{2, 4} &=(\omega, \overline{\omega}, 1, \omega, \overline{\omega}, 1),\\
  l_{1, 5} \cap l_{2, 5} &=\emptyset,\\
  l_{1, 5} \cap l_{2, 6} &=\emptyset.
\endaligned\right.$$
$$\left\{\aligned
  l_{1, 6} \cap l_{2, 1} &=\emptyset,\\
  l_{1, 6} \cap l_{2, 2} &=(\omega, \overline{\omega}, 1, \overline{\omega}, \omega, 1),\\
  l_{1, 6} \cap l_{2, 3} &=(\omega, 1, \overline{\omega}, \omega, \overline{\omega}, 1),\\
  l_{1, 6} \cap l_{2, 4} &=\emptyset,\\
  l_{1, 6} \cap l_{2, 5} &=\emptyset,\\
  l_{1, 6} \cap l_{2, 6} &=\emptyset.
\endaligned\right.$$

  The double six is given by
$$N=\left(\matrix
    a_1 & a_2 & a_3 & a_4 & a_5 & a_6\\
    b_1 & b_2 & b_3 & b_4 & b_5 & b_6
     \endmatrix\right)
  =\left(\matrix
    l_{1, 1} & l_{2, 2} & l_{2, 3} & l_{1, 4} & l_{1, 5} & l_{2, 6}\\
    l_{2, 1} & l_{1, 2} & l_{1, 3} & l_{2, 4} & l_{2, 5} & l_{1, 6}
    \endmatrix\right).\tag 4.11$$

  We find that
$$\aligned
 l_{1, 1} \cap l_{1, 2}=\overline{l_{2, 1} \cap l_{2, 2}} &=(-2 \omega, \omega, \omega, -2, 1, 1),\\
 l_{1} \cap l_{1, 1}=\overline{l_{1} \cap l_{2, 1}} &=(0, -\omega, \omega, 0, -1, 1),\\
 l_{1} \cap l_{1, 2}=\overline{l_{1} \cap l_{2, 2}} &=(0, \omega, -\omega, 0, -1, 1).
\endaligned$$
Note that $\overline{l_1}=l_1$. Hence, $l_1 \cap l_{1, 1} \in
l_1$, $l_1 \cap l_{2, 1} \in l_1$, $l_1 \cap l_{1, 2} \in l_1$,
$l_1 \cap l_{2, 2} \in l_1$. Thus, $c_{12}=l_1.$

  By the same method as above, we obtain that
$$c_{12}=l_1, \quad c_{13}=l_9, \quad c_{14}=l_{0, 5}, \quad c_{15}=l_{0, 4}, \quad c_{16}=l_5.\tag 4.12$$
$$c_{23}=l_{0, 6}, \quad c_{24}=l_6, \quad c_{25}=l_8, \quad c_{26}=l_{0, 3}.\tag 4.13$$
$$c_{34}=l_2, \quad c_{35}=l_4, \quad c_{36}=l_{0, 2}.\tag 4.14$$
$$c_{45}=l_{0, 1}, \quad c_{46}=l_7, \quad c_{56}=l_3.\tag 4.15$$

  Therefore, we have
$$E(a_1)=a_2, \quad A(a_1)=a_4, \quad B(a_1)=a_1, \quad C(a_1)=a_4.$$
$$E(a_2)=a_1, \quad A(a_2)=a_6, \quad B(a_2)=a_2, \quad C(a_2)=a_3.$$
$$E(a_3)=a_5, \quad A(a_3)=a_2, \quad B(a_3)=a_6, \quad C(a_3)=a_6.$$
$$E(a_4)=a_3, \quad A(a_4)=a_5, \quad B(a_4)=a_5, \quad C(a_4)=a_5.$$
$$E(a_5)=a_6, \quad A(a_5)=a_1, \quad B(a_5)=a_4, \quad C(a_5)=a_1.$$
$$E(a_6)=a_4, \quad A(a_6)=a_3, \quad B(a_6)=a_3, \quad C(a_6)=a_2.$$
$$E(b_1)=b_2, \quad A(b_1)=b_4, \quad B(b_1)=b_1, \quad C(b_1)=b_4.$$
$$E(b_2)=b_1, \quad A(b_2)=b_6, \quad B(b_2)=b_2, \quad C(b_2)=b_3.$$
$$E(b_3)=b_5, \quad A(b_3)=b_2, \quad B(b_3)=b_6, \quad C(b_3)=b_6.$$
$$E(b_4)=b_3, \quad A(b_4)=b_5, \quad B(b_4)=b_5, \quad C(b_4)=b_5.$$
$$E(b_5)=b_6, \quad A(b_5)=b_1, \quad B(b_5)=b_4, \quad C(b_5)=b_1.$$
$$E(b_6)=b_4, \quad A(b_6)=b_3, \quad B(b_6)=b_3, \quad C(b_6)=b_2.$$
$$E(c_{12})=c_{12}, \quad A(c_{12})=c_{46}, \quad B(c_{12})=c_{12}, \quad C(c_{12})=c_{34}.$$
$$E(c_{13})=c_{25}, \quad A(c_{13})=c_{24}, \quad B(c_{13})=c_{16}, \quad C(c_{13})=c_{46}.$$
$$E(c_{14})=c_{23}, \quad A(c_{14})=c_{45}, \quad B(c_{14})=c_{15}, \quad C(c_{14})=c_{45}.$$
$$E(c_{15})=c_{26}, \quad A(c_{15})=c_{14}, \quad B(c_{15})=c_{14}, \quad C(c_{15})=c_{14}.$$
$$E(c_{16})=c_{24}, \quad A(c_{16})=c_{34}, \quad B(c_{16})=c_{13}, \quad C(c_{16})=c_{24}.$$
$$E(c_{23})=c_{15}, \quad A(c_{23})=c_{26}, \quad B(c_{23})=c_{26}, \quad C(c_{23})=c_{36}.$$
$$E(c_{24})=c_{13}, \quad A(c_{24})=c_{56}, \quad B(c_{24})=c_{25}, \quad C(c_{24})=c_{35}.$$
$$E(c_{25})=c_{16}, \quad A(c_{25})=c_{16}, \quad B(c_{25})=c_{24}, \quad C(c_{25})=c_{13}.$$
$$E(c_{26})=c_{14}, \quad A(c_{26})=c_{36}, \quad B(c_{26})=c_{23}, \quad C(c_{26})=c_{23}.$$
$$E(c_{34})=c_{35}, \quad A(c_{34})=c_{25}, \quad B(c_{34})=c_{56}, \quad C(c_{34})=c_{56}.$$
$$E(c_{35})=c_{56}, \quad A(c_{35})=c_{12}, \quad B(c_{35})=c_{46}, \quad C(c_{35})=c_{16}.$$
$$E(c_{36})=c_{45}, \quad A(c_{36})=c_{23}, \quad B(c_{36})=c_{36}, \quad C(c_{36})=c_{26}.$$
$$E(c_{45})=c_{36}, \quad A(c_{45})=c_{15}, \quad B(c_{45})=c_{45}, \quad C(c_{45})=c_{15}.$$
$$E(c_{46})=c_{34}, \quad A(c_{46})=c_{35}, \quad B(c_{46})=c_{35}, \quad C(c_{46})=c_{25}.$$
$$E(c_{56})=c_{46}, \quad A(c_{56})=c_{13}, \quad B(c_{56})=c_{34}, \quad C(c_{56})=c_{12}.$$

  For
$$N=\left(\matrix
    a_1 & a_2 & a_3 & a_4 & a_5 & a_6\\
    b_1 & b_2 & b_3 & b_4 & b_5 & b_6
     \endmatrix\right),$$
we have
$$E(N)=\left(\matrix
       a_2 & a_1 & a_5 & a_3 & a_6 & a_4\\
       b_2 & b_1 & b_5 & b_3 & b_6 & b_4
\endmatrix\right).$$
$$A(N)=\left(\matrix
       a_4 & a_6 & a_2 & a_5 & a_1 & a_3\\
       b_4 & b_6 & b_2 & b_5 & b_1 & b_3
\endmatrix\right).$$
$$B(N)=\left(\matrix
       a_1 & a_2 & a_6 & a_5 & a_4 & a_3\\
       b_1 & b_2 & b_6 & b_5 & b_4 & b_3
\endmatrix\right).$$
$$C(N)=\left(\matrix
       a_4 & a_3 & a_6 & a_5 & a_1 & a_2\\
       b_4 & b_3 & b_6 & b_5 & b_1 & b_2
\endmatrix\right).$$

  For
$$N_{12}=\left(\matrix
    a_1 & b_1 & c_{23} & c_{24} & c_{25} & c_{26}\\
    a_2 & b_2 & c_{13} & c_{14} & c_{15} & c_{16}
     \endmatrix\right),$$
we have
$$E(N_{12})=\left(\matrix
    a_2 & b_2 & c_{15} & c_{13} & c_{16} & c_{14}\\
    a_1 & b_1 & c_{25} & c_{23} & c_{26} & c_{24}
\endmatrix\right).$$
$$A(N_{12})=\left(\matrix
    a_4 & b_4 & c_{26} & c_{56} & c_{16} & c_{36}\\
    a_6 & b_6 & c_{24} & c_{54} & c_{14} & c_{34}
\endmatrix\right).$$
$$B(N_{12})=\left(\matrix
    a_1 & b_1 & c_{26} & c_{25} & c_{24} & c_{23}\\
    a_2 & b_2 & c_{16} & c_{15} & c_{14} & c_{13}
\endmatrix\right).$$
$$C(N_{12})=\left(\matrix
    a_4 & b_4 & c_{36} & c_{35} & c_{31} & c_{32}\\
    a_3 & b_3 & c_{46} & c_{45} & c_{41} & c_{42}
\endmatrix\right).$$

  For
$$N_{123}=\left(\matrix
  a_1 & a_2 & a_3 & c_{56} & c_{46} & c_{45}\\
  c_{23} & c_{13} & c_{12} & b_4 & b_5 & b_6
\endmatrix\right),$$
we have
$$E(N_{123})=\left(\matrix
  a_2 & a_1 & a_5 & c_{46} & c_{34} & c_{36}\\
  c_{15} & c_{25} & c_{12} & b_3 & b_6 & b_4
\endmatrix\right).$$
$$A(N_{123})=\left(\matrix
  a_4 & a_6 & a_2 & c_{13} & c_{35} & c_{15}\\
  c_{26} & c_{24} & c_{46} & b_5 & b_1 & b_3
\endmatrix\right).$$
$$B(N_{123})=\left(\matrix
  a_1 & a_2 & a_6 & c_{34} & c_{35} & c_{45}\\
  c_{26} & c_{16} & c_{12} & b_5 & b_4 & b_3
\endmatrix\right).$$
$$C(N_{123})=\left(\matrix
  a_4 & a_3 & a_6 & c_{12} & c_{25} & c_{15}\\
  c_{36} & c_{46} & c_{34} & b_5 & b_1 & b_2
\endmatrix\right).$$

  The other $N_{ij}$ and $N_{ijk}$ can be calculated similarly.
We find that the action of the group $H$ generated by $A$, $B$,
$C$ and $E$ on the $27$ lines gives the permutations of the lines
which preserve the intersection behavior of the lines. Hence, we
conclude that the group $H$ is a subgroup of $\text{Aut}({\Cal
L})=W(E_6)$.

{\smc Theorem 4.8} (Main Theorem 4). {\it The group $H$ generated
by $A$, $B$, $C$ and $E$ is a subgroup of $\text{Aut}({\Cal
L})=W(E_6)$.}

  For the complex conjugate $\sigma_c \in \text{Gal}(\overline{{\Bbb Q}}/{\Bbb Q})$,
we have
$$\sigma_c(N)=\left(\matrix
  b_1 & b_2 & b_3 & b_4 & b_5 & b_6\\
  a_1 & a_2 & a_3 & a_4 & a_5 & a_6
\endmatrix\right).$$
$$\sigma_c(N_{12})=\left(\matrix
  b_1 & a_1 & c_{23} & c_{24} & c_{25} & c_{26}\\
  b_2 & a_2 & c_{13} & c_{14} & c_{15} & c_{16}
\endmatrix\right).$$
$$\sigma_c(N_{123})=\left(\matrix
  b_1 & b_2 & b_3 & c_{56} & c_{46} & c_{45}\\
  c_{23} & c_{13} & c_{12} & a_4 & a_5 & a_6
\endmatrix\right).$$

  This leads us to study the representation:
$$\rho: G_{{\Bbb Q}}=\text{Gal}(\overline{{\Bbb Q}}/{\Bbb Q}) \to \text{Aut}({\Cal L}),\tag 4.16$$
i.e.,
$$\rho: \text{Gal}(\overline{{\Bbb Q}}/{\Bbb Q}) \to W(E_6).\tag 4.17$$
Here,
$$\text{Im}(\rho)=H=\langle A, B, C, E \rangle, \quad
  \text{ker}(\rho)=\text{Gal}(\overline{{\Bbb Q}}/{\Cal K}).\tag 4.18$$
We have
$$\text{Im}(\rho) \cong \text{Gal}(\overline{{\Bbb Q}}/{\Bbb Q})/\text{ker}(\rho)
 =\text{Gal}(\overline{{\Bbb Q}}/{\Bbb Q})/\text{Gal}(\overline{{\Bbb Q}}/{\Cal K})
 \cong \text{Gal}({\Cal K}/{\Bbb Q}).\tag 4.19$$
Hence,
$$\text{Gal}({\Cal K}/{\Bbb Q}) \cong H.\tag 4.20$$
The (geometrical, analytical and algebraical) properties of the
Hessian group $G=\langle A, B, C, D, E \rangle$ have been studied
in our previous paper \cite{Y3}. Note that $A^{*} A=B^{*} B=C^{*}
C=I$, but $E^{*} E \neq I$. We have $\text{Tr}(A)=3$, $\det(A)=1$.
$\text{Tr}(B)=2$, $\det(B)=1$. $\text{Tr}(C)=3$, $\det(C)=1$.
$\text{Tr}(E)=0$.
$$\det(E)=\frac{1}{(\sqrt{-3})^{18}} \det
          \left(\matrix
          1 & 1 & 1 & 9 & 3 & 3\\
          1 & 1 & 1 & 9 & 3 \overline{\omega} & 3 \omega\\
          1 & 1 & 1 & 9 & 3 \omega & 3 \overline{\omega}\\
          1 & 1 & 1 & 0 & 0 & 0\\
          3 & 3 \overline{\omega} & 3 \omega & 0 & 0 & 0\\
          3 & 3 \omega & 3 \overline{\omega} & 0 & 0 & 0
          \endmatrix\right)=-1.$$
Note that
$$E^2=\left(\matrix
           -1 &    &    &    &    &   \\
              &  0 & -1 &    &    &   \\
              & -1 &  0 &    &    &   \\
              &    &    & -1 &    &   \\
              &    &    &    &  0 & -1\\
              &    &    &    & -1 &  0
      \endmatrix\right).$$
In fact, the elements of $H$ have integer determinant and trace.

  Now let us recall the definition of Artin $L$-functions (see
\cite{Ar1}, \cite{G}). Suppose $K$ is a number field, and $E$ is a
finite Galois extension of $K$ with Galois group
$G=\text{Gal}(E/K)$. By a representation of $G$ we understand a
homomorphism $\sigma$ of $G$ into $GL(V)$, the group of invertible
linear transformations of a complex vector space of dimension $n$.
Given a finite place ${\frak p}$ of $K$, we say a prime ${\frak
P}$ of ${\Cal O}_E$ lies over ${\frak p}$ if ${\frak P}$ appears
in the factorization of ${\frak p} {\Cal O}_E$ into prime ideals
of ${\Cal O}_E$. Given such a pair ${\frak P}/{\frak p}$, we have:

\roster

\item The decomposition group $D_{{\frak P}}=\{ g \in G: g({\frak
      P})={\frak P} \};$

\item The inertia subgroup $I_{{\frak P}}=\{ g \in D_{{\frak P}}:
      \text{$g(x)=x$(mod ${\frak P}$) for all $x \in {\Cal
      O}_E$} \}$;

\item The Frobenius automorphism $\text{Fr}_{{\frak P}}$
      generating the cyclic group $D_{{\frak P}}/I_{{\frak P}}$ which
      is isomorphic to $\text{Gal}({\Cal O}_E/{\frak P}:{\Cal O}_K/{\frak p})$.

\endroster

  Let $V_{{\frak P}}$ denote the subspace of $V$ formed by vectors
invariant by $\sigma(I_{{\frak P}})$. Then
$\sigma(\text{Fr}_{{\frak P}})$ is defined on $V_{{\frak P}}$, and
the Euler factor
$$L_{{\frak p}}(\sigma, s)=[\det(I-\sigma(\text{Fr}_{{\frak P}})
                           \cdot N {\frak p}^{-s}|_{V_{{\frak P}}})]^{-1}$$
depends only on ${\frak p}$, not ${\frak P}$. The Artin
$L$-function attached to $\sigma$ is given by
$$L(\sigma, s)=L_{E/K}(\sigma, s)=\prod_{{\frak p}} L_{{\frak p}}(\sigma, s)$$
which is convergent for $\text{Re}(s)>1$.

  In \cite{Ar2}, Artin studied representations with a finite image
$\text{Im}(\rho)$ (see also \cite{MP}). In this case the
representation $\rho$ is always rational for some number field and
is semisimple. Hence $\rho$ is uniquely determined (up to
equivalence) by its character $\chi_{\rho}$. If $\rho$ is an
arbitrary rational representation then the function $L(\rho, s)$
is uniquely defined by its character $\chi_{\rho}$. In view of
this fact one often uses the notation $L(\chi_{\rho}, s)=L(\rho,
s)$.

  For arbitrary $n$, there is the following remarkable Reciprocity
Conjecture:

{\smc Conjecture 4.9} (Langlands). {\it Suppose $E$ is a finite
Galois extension of $F$ with Galois group $G=\text{Gal}(E/F)$, and
$\sigma: G \to GL_n({\Bbb C})$ is an irreducible representation of
$G$. Then there exists an automorphic cuspidal representation
$\pi_{\sigma}$ on $GL_n$ over $F$ such that $L(s,
\pi_{\sigma})=L(s, \sigma)$.}

  The truth of this conjecture implies the truth of Artin conjecture
on the entirety of the Artin $L$-function $L(s, \sigma)$.

{\smc Conjecture 4.10} (Main Conjecture). {\it Let $\rho:
\text{Gal}(\overline{{\Bbb Q}}/{\Bbb Q}) \to W(E_6)$ be a Galois
representation. When the image of $\rho$: $\text{Im}(\rho) \cong H
\leq W(E_6)$, if $\rho$ is odd, i.e., $\det \rho(\sigma_c)=-1$,
where $\sigma_c \in \text{Gal}(\overline{{\Bbb Q}}/{\Bbb Q})$ is
the complex conjugate, then there exists a Picard modular form $f$
of weight one such that
$$L(\rho, s)=L(f, s)$$
up to finitely many Euler factors, where $L(\rho, s)$ is the Artin
$L$-function associated to the Galois representation $\rho$ and
$L(f, s)$ is the automorphic $L$-function associated to the Picard
modular form $f$. If $\rho$ is even, i.e., $\det
\rho(\sigma_c)=1$, then there exists an automorphic form $f$
associated to $U(2, 1):$ $\Delta f=s(s-2)f$ with $s=1$ and $($see
\cite{Y1}$)$
$$\aligned
  \Delta &=(z_1+\overline{z_1}-z_2 \overline{z_2}) \left[(z_1+\overline{z_1})
    \frac{\partial^2}{\partial z_1 \partial \overline{z_1}}+
    \frac{\partial^2}{\partial z_2 \partial \overline{z_2}}+
    z_2 \frac{\partial^2}{\partial \overline{z_1} \partial z_2}+
    \overline{z_2} \frac{\partial^2}{\partial z_1 \partial
    \overline{z_2}}\right],\\
        &(z_1, z_2) \in {\frak S}_2=\{ (z_1, z_2) \in {\Bbb C}^2:
    z_1+\overline{z_1}-z_2 \overline{z_2}>0\},
\endaligned$$
such that
$$L(\rho, s)=L(f, s)$$
up to finitely many Euler factors, where $L(\rho, s)$ is the Artin
$L$-function associated to the Galois representation $\rho$ and
$L(f, s)$ is the automorphic $L$-function associated to the
automorphic form $f$.}

  Here, we find that three topics: $L$-functions, Galois representations
and automorphic forms or, equivalently, representations, which
figure prominently in the modern higher arithmetic, appear
naturally in our conjecture. The $L$-functions are attached to
both the Galois representations and the automorphic
representations and are the link that joins them. Our conjecture,
like Langlands reciprocity law and Artin's abelian reciprocity
law, gives a correspondence independent of any particular
geometric construction.

\roster

\item tetrahedron, $A_4$ (Langlands) (see \cite{L2}), elliptic
modular forms of
      weight $1$.

\item octahedron, $S_4$ (Tunnell), elliptic modular forms of
      weight $1$.

\item icosahedron, $A_5$ (Artin conjecture), elliptic modular
      forms of weight $1$.

\item Hessian polyhedra, $H$ (our conjecture), Picard modular
      forms of weight $1$.

\endroster

  Note that the group $W(E_6)$ is combinatorially defined, while
the subgroup $H$ is defined on algebraic variety. Thus we give a
connection between combinatoric and algebraic geometry.
Furthermore, the Galois representation $\rho:
\text{Gal}(\overline{{\Bbb Q}}/{\Bbb Q}) \to H $ comes from number
theory, the Picard modular forms come from analysis and
representation theory of Lie groups. Here, Galois symmetry and
geometric symmetry meet together. Therefore, our theorem and
conjecture gives a connection which involves combinatoric,
algebraic geometry, number theory, analysis and representation
theory.

\vskip 0.5 cm

\centerline{\bf 5. Hessian polyhedra and the arithmetic of
                   rigid Calabi-Yau threefolds}

\vskip 0.5 cm

  A complete system of invariants for the Hessian groups (the
corresponding geometric objects are Hessian polyhedra) has degrees
$6$, $9$, $12$, $12$ and $18$ and can be given explicitly by the
following forms:
$$\left\{\aligned
  C_6(z_1, z_2, z_3) &=z_1^6+z_2^6+z_3^6-10
                     (z_1^3 z_2^3+z_2^3 z_3^3+z_3^3 z_1^3),\\
  C_9(z_1, z_2, z_3) &=(z_1^3-z_2^3)(z_2^3-z_3^3)(z_3^3-z_1^3),\\
  C_{12}(z_1, z_2, z_3) &=(z_1^3+z_2^3+z_3^3)[(z_1^3+z_2^3+z_3^3)^3
                        +216 z_1^3 z_2^3 z_3^3],\\
  {\frak C}_{12}(z_1, z_2, z_3) &=z_1 z_2 z_3 [27 z_1^3 z_2^3 z_3^3
                                -(z_1^3+z_2^3+z_3^3)^3],\\
  C_{18}(z_1, z_2, z_3) &=(z_1^3+z_2^3+z_3^3)^6-540 z_1^3 z_2^3 z_3^3
                        (z_1^3+z_2^3+z_3^3)^3-5832 z_1^6 z_2^6 z_3^6.
\endaligned\right.$$
They satisfy the following relations:
$$\left\{\aligned
  432 C_9^2 &=C_6^3-3 C_6 C_{12}+2 C_{18},\\
  1728 {\frak C}_{12}^3 &=C_{18}^{2}-C_{12}^{3}.
\endaligned\right.$$
The first relation can be considered as an elliptic curve with two
parameters $C_{12}$ and $C_{18}$:
$$432 C_9^2=C_6^3-3 C_{12} C_6+2 C_{18}.\tag 5.1$$
The second one can be considered as an elliptic curve with a
parameter ${\frak C}_{12}$:
$$C_{18}^2=C_{12}^3+1728 {\frak C}_{12}^3.\tag 5.2$$

  When $C_9=0$, the equation (5.1) reduces to
$$C_6^3-3 C_6 C_{12}+2 C_{18}=0.$$
In fact, Wirtinger studied the algebraic function $z$ of two
variables $x$, $y$ defined by the equation (see \cite{BE})
$$z^3+3xz+2y=0.$$
In modern terms: the projection of the surface $X$ with this
equation to the $(x, y)$-plane is the semiuniversal unfolding of
the $0$-dimensional $A_2$-type singularity $z^3=0$. The
discriminant curve $D \subset {\Bbb C}^2$ has the equation
$$x^3+y^2=0, \quad \text{i.e.}, \quad 1728 {\frak C}_{12}^3=0.$$
The fundamental group of ${\Bbb C}^2-D$ operates on the fibre over
the base point by the monodromy representation. Wirtinger
calculates $\pi_1({\Bbb C}^2-D)$ and finds a presentation with two
generators and one relation $sts=tst$. In modern terms: $\pi_1$ is
the braid group on $3$ strings. The monodromy representation is
the canonical homomorphism of this group to the symmetric group
$S_3$. In his computation of $\pi_1({\Bbb C}^2-D)$, Wirtinger used
an idea of Heegaard. Heegaard reduced the complex geometry of an
algebroid covering $(X, x) \to ({\Bbb C}^2, 0)$ with a singularity
$(D, 0)$ of the discriminant to a situation of $3$-dimensional
topology.  He considered a small $4$-ball $B \subset {\Bbb C}^2$
centered at $0$ with boundary $\partial B=S^3$, a $3$-sphere. The
intersection $L=D \cap S^3$ is a knot or link in $S^3$. In
Wirtinger's example it is the trefoil knot. This established a
link between the geometry of singularities of complex surfaces and
$3$-dimensional topology.

  Put
$$x=3 C_6, \quad y=108 C_9$$
in the elliptic curve (5.1), we have
$$E: \quad y^2=x^3-27 C_{12} x+54 C_{18}. \tag 5.3$$
The $j$-invariant is given by
$$j=-\frac{C_{12}^3}{{\frak C}_{12}^3}.\tag 5.4$$

  According to \cite{Y3}, put
$$\varphi=z_1 z_2 z_3, \quad
  \psi=z_1^3+z_2^3+z_3^3, \quad
  \chi=z_1^3 z_2^3+z_2^3 z_3^3+z_3^3 z_1^3.$$
$$G=(z_1-z_2)(z_2-z_3)(z_3-z_1).$$
$$H=\psi+6 \varphi, \quad K=\psi-3 \varphi.$$
$$\left\{\aligned
     C_6 &=\psi^2-12 \chi,\\
  C_{12} &=\psi (\psi^3+216 \varphi^3),\\
  {\frak C}_{12} &=\varphi (27 \varphi^3-\psi^3),\\
  C_{18} &=\psi^6-540 \varphi^3 \psi^3-5832 \varphi^6.
\endaligned\right.$$

{\smc Theorem 5.1} (see \cite{Y3}, Main Theorem 4). {\it The
invariants $G$, $H$ and $K$ satisfy the following algebraic
equations, which are the form-theoretic resolvents $($algebraic
resolvents$)$ of $G$, $H$, $K$$:$
$$\left\{\aligned
  4 G^3+H^2 G-C_6 G-4 C_9 &=0,\\
  H (H^3+8 K^3)-9 C_{12} &=0,\\
  K (K^3-H^3)-27 {\frak C}_{12} &=0.
\endaligned\right.$$}

  We have
$$C_{12}=\frac{1}{9} H (H^3+8 K^3), \quad
  {\frak C}_{12}=\frac{1}{27} K (K^3-H^3).$$
Hence,
$$j=-27 \frac{H^3 (H^3+8 K^3)^3}{K^3 (K^3-H^3)^3}.\tag 5.5$$
Put
$$t=\frac{H}{K}.$$
Then
$$j=-27 \frac{t^3 (t^3+8)^3}{(1-t^3)^3}.\tag 5.6$$

  Note that
$$\varphi=\frac{1}{9}(H-K), \quad \psi=\frac{1}{3}(H+2K).$$
Hence,
$$C_{18}=-\frac{1}{27}(H^6-20 H^3 K^3-8 K^6).\tag 5.7$$
Now, we get
$$E: \quad y^2=x^3-3H(H^3+8K^3)x-2(H^6-20 H^3 K^3-8 K^6).$$
Put
$$X=\frac{x}{K^2}, \quad Y=\frac{y}{K^3}.$$
We have
$$E_{2, t}: \quad Y^2=X^3-3t(t^3+8)X-2(t^6-20 t^3-8).\tag 5.8$$

  For the elliptic curve
$$E_{2, t}/{\Bbb Q}(t): \quad y^2=x^3-3t(t^3+8)x-2(t^6-20 t^3-8),\tag 5.9$$
we find that
$$P_2=(3t^2, 4(t^3-1)) \in E_{2, t}.\tag 5.10$$
Moreover,
$$[2] P_2=(3 t^2, -4 (t^3-1)), \quad -P_2=(3 t^2, -4 (t^3-1)).\tag 5.11$$
Hence, $[2] P_2=-P_2$, i.e., $[3] P_2=O$. So, $P_2$ is a
$3$-division point in $E_{2, t}$.

  The discriminant of the elliptic curve (5.9) is given by
$$\Delta=2^{12} \cdot 3^3 \cdot (t-1)^3 (t^2+t+1)^3.\tag 5.12$$
The bad fibers occur at
$$t=1, \omega, \overline{\omega}, \infty.$$
They are of types $(I_3, I_3, I_3, I_3)$ by \cite{Her}, p.336.
Correspondingly,
$$\left\{\aligned
  &t=1, \text{i.e.}, \varphi=0, \text{i.e.}, z_1 z_2 z_3=0,\\
  &t=\omega, \text{i.e.}, \psi-3 \overline{\omega} \varphi=0,
   \text{i.e.}, (z_1+z_2+\overline{\omega} z_3)(\omega z_1+z_2+\omega z_3)
   (\overline{\omega} z_1+z_2+z_3)=0,\\
  &t=\overline{\omega}, \text{i.e.}, \psi-3 \omega \varphi=0,
   \text{i.e.}, (z_1+\omega z_2+z_3)(\omega z_1+z_2+z_3)
   (\overline{\omega} z_1+\overline{\omega} z_2+z_3)=0,\\
  &t=\infty, \text{i.e.}, \psi-3 \varphi=0,
   \text{i.e.}, (z_1+z_2+z_3)(z_1+\omega z_2+\overline{\omega} z_3)
   (z_1+\overline{\omega} z_2+\omega z_3)=0,
\endaligned\right.$$
where $\omega=\exp(2 \pi i/3)$. They are the $12$ lines of the
Hessian arrangement on ${\Bbb P}^2$.

  Therefore, for the elliptic curve (5.3), we have
$$P_1=(3 C_6, 108 C_9) \in E,\tag 5.13$$
$$P_2=(3 H^2, 4 (H^3-K^3))=(3(\psi+6 \varphi)^2, 108 \varphi(\psi^2+3 \psi \varphi+9 \varphi^2)) \in E.\tag 5.14$$

  Note that
$$t=\frac{H}{K}=\frac{\psi+6 \varphi}{\psi-3 \varphi}.$$
This gives the following identity:
$$E_{1, t}: \quad (t-1)(z_1^3+z_2^3+z_3^3)-3(t+2) z_1 z_2 z_3=0,\tag 5.15$$
which is a family of elliptic curves over ${\Bbb P}^1$. In the
singular fiber: $t=1$, $E_{1, t}$ degenerates into $z_1 z_2
z_3=0$. In the singular fiber: $t=\infty$, $E_{1, t}$ degenerates
into $(z_1+z_2+z_3)(z_1+\omega z_2+\omega^2 z_3)(z_1+\omega^2
z_2+\omega z_3)=0.$ In the singular fiber: $t=\overline{\omega}$,
$E_{1, t}$ degenerates into $(z_1+\omega z_2+z_3)(\omega
z_1+z_2+z_3)(\overline{\omega} z_1+\overline{\omega} z_2+z_3)=0$.
In the singular fiber: $t=\omega$, $E_{1, t}$ degenerates into
$(z_1+z_2+\overline{\omega} z_3)(\omega z_1+z_2+\omega
z_3)(\overline{\omega} z_1+z_2+z_3)=0$. They are the $12$ lines of
the Hessian arrangement on ${\Bbb P}^2$. Now, $E$ decomposes into
the pair $(E_{1, t}, E_{2, t})$.

  By the linear transformation $t=9r+1$, $E_{1, t}$ becomes one of
the Beauville's families of elliptic curves, which are modular
elliptic surfaces for some conjugate of the given group (see
\cite{B} and \cite{Ve2}). The equation in ${\Bbb P}^2$ is given by
$$III: \quad z_1^3+z_2^3+z_3^3=\left(\frac{1}{r}+3\right) z_1 z_2 z_3.$$
The group $\Gamma=\Gamma(3)$. The number of components of singular
fibers is $(3, 3, 3, 3)$. The Picard-Fuchs equation for this
family of elliptic curves is given by
$$III: \quad r(27 r^2+9r+1) f^{\prime \prime}+(9r+1)^2 f^{\prime}+3(9r+1) f=0.$$
The solution of this equation in terms of modular forms is
$$III: \quad r=\eta(3 \tau)^3 \eta \left(\frac{1}{3} \tau\right)^{-3}, \quad
             f=\eta \left(\frac{1}{3} \tau\right)^3 \eta(\tau)^{-1}.$$
In our present case, the Picard-Fuchs equation becomes
$$(t^3-1) \frac{d^2 f}{d t^2}+3 t^2 \frac{df}{dt}+tf=0.\tag 5.16$$
In fact, for $\overline{C}={\Bbb P}^1$, $C={\Bbb P}^1-\{ 1,
\infty, \text{roots of} \quad t^2+t+1=0 \}$ and the family $f: X
\to C$ is the modular family for the group $\Gamma_1(3)$. The
Picard-Fuchs equation (for the holomorphic $1$-forms on $X_t$) is
(see \cite{Pe})
$$(t-1)(t^2+t+1) \frac{d^2 f}{dt^2}+3t^2 \frac{d f}{dt}+t f=0.$$
The solutions are given by
$$t=9 \frac{\eta(3 \tau)^3}{\eta \left(\frac{1}{3} \tau\right)^3}+1, \quad
  f=\frac{\eta \left(\frac{1}{3} \tau\right)^3}{\eta(\tau)}.\tag 5.17$$
$\Gamma(3)$ is a genus zero principal congruence subgroup of
$PSL(2, {\Bbb Z})$, the index of $\Gamma(3)$ in $PSL(2, {\Bbb Z})$
is $\mu=12$. As a rational function of the Hauptmodul $\rho$, the
$j$-function can be expressed as (see \cite{MS})
$$j(\rho)=\frac{\rho^3 (\rho+6)^3 (\rho^2-6 \rho+36)^3}{(\rho-3)^3 (\rho^2+3 \rho+9)^3},$$
which is also the $J$-invariant of the semi-stable arithmetic
elliptic surface over ${\Bbb P}^1$. In the elliptic curve $E_{1,
t}$, set
$$\rho=\frac{z_1^3+z_2^3+z_3^3}{z_1 z_2 z_3}=\frac{3(t+2)}{t-1}.\tag 5.18$$
We find that
$$j(\rho)=27 \frac{t^3 (t^3+8)^3}{(t^3-1)^3},\tag 5.19$$
which is the $j$-invariant of the elliptic curve $E_{2, t}$.

  Let us recall the following well-known results (see \cite{Shi1}):
let $k$ be a field of characteristic $\neq 3$ containing $3$ cubic
roots of unity. Consider the elliptic curve
$$E_3: \quad x^3+y^3+z^3-3 \mu xyz=0$$
defined over $K_3=k(\mu)$, $\mu$ being a variable over $k$. Then
$E_3$ has exactly $9$ $K_3$-rational points (i.e. base points of
the pencil) and they are of order $3$.

{\smc Proposition 5.2}. (Deuring Normal Form) (see \cite{Sil}).
{\it Let $E/K$ be an elliptic curve over a field with
$\text{char}(K) \neq 3$. Then $E$ has a Weierstrass equation over
$\overline{K}$ of the form
$$E_{\alpha}: \quad y^2+\alpha xy+y=x^3, \quad \alpha \in
  \overline{K}, \alpha^3 \neq 27.$$
This equation has discriminant and $j$-invariant
$$\Delta=\alpha^3-27, \quad j=\alpha^3
(\alpha^3-24)^3/(\alpha^3-27).$$}

  We will give some basic facts about families of elliptic curves with
points of order $3$: the Hessian family (see \cite{Hus}).

{\it Definition} 5.3. The Hessian family of elliptic curves
$E_{\alpha}: y^2+\alpha xy+y=x^3$ is defined for any field of
characteristic different from $3$ with $j$-invariant
$j(\alpha)=j(E_{\alpha})$ of $E_{\alpha}$ given by
$$j(\alpha)=\frac{\alpha^3 (\alpha^3-24)^3}{\alpha^3-27}.$$

  The curve $E_{\alpha}$ is nonsingular for $\alpha^3 \neq 27$, that is,
$\alpha$ is not in $3 \mu_3$, where $\mu_3$ is the group of third
roots of unity. Over the line $k$ minus $3$ points, $k-3 \mu_3$,
the family $E_{\alpha}$ consists of elliptic curves with a
constant section $(0, 0)$ of order $3$ where $2 (0, 0)=(0, -1)$.
The Hessian family $E_{\alpha}$ has three singular fibres which
are nodal cubics at the points of $3 \mu_3=\{ 3, 3 \omega, 3
\omega^2 \}$, where $\omega^2+\omega+1=0$.

  There are other versions of the Hessian family which in homogeneous
coordinates take the form
$$H_{\mu}: \quad u^3+v^3+w^3=3 \mu uvw,$$
or in affine coordinates with $w=-1$, it has the form
$$u^3+v^3=1-3 \mu uv.$$
If we set $y=-v^3$ and $x=-uv$, we obtain $x^3/y-y=1+3 \mu x$, or
$$E_{3 \mu}: \quad y^2+3 \mu xy+y=x^3.$$
This change of variable defines what is called a $3$-isogeny of
$H_{\mu}$ onto $E_{3 \mu}$. There are nine cross-sections of the
family $H_{\mu}$ given by
$$\matrix
  (0, -1, 1), & (0, -\omega, 1), & (0, -\omega^2, 1),\\
  (1, 0, -1), & (\omega, 0, -\omega^2), & (\omega^2, 0, -\omega),\\
  (-1, 1, 0), & (-1, \omega^2, 0), & (-1, \omega, 0).
\endmatrix$$
The family $H_{\mu}$ is nonsingular over the line minus $\mu_3$,
and choosing, for example, $0=(-1, 1, 0)$, one can show that these
nine points form the subgroup of $3$-division points of the family
$H_{\mu}$.

  For $E_{3t}: y^2+3txy+y=x^3$, we have
$$j_1=\frac{27 t^3 (9 t^3-8)^3}{t^3-1}.$$
For $E_{2, t}: y^2=x^3-3t(t^3+8)x-2(t^6-20 t^3-8)$, we have
$$j_2=\frac{27 t^3 (t^3+8)^3}{(t^3-1)^3}.$$
Let $r=t^3-1$, then
$$\left\{\aligned
  j_1 &=j_1(r)=\frac{27}{r} (r+1)(9r+1)^3,\\
  j_2 &=j_2(r)=\frac{27}{r^3}(r+1)(r+9)^3.
\endaligned\right.$$
Note that
$$j_1 \left(\frac{1}{r}\right)=j_2(r).$$
This gives a one-to-one correspondence between two families of
elliptic curves $E_{3t}$ and $E_{2, t}$.

  We will make use of certain elliptic modular surfaces associated
to torsion free congruence subgroups of $SL(2, {\Bbb Z})$. Each of
these has four singular fibers, all of which are of type $I_b$ for
various $b$. This is actually a complete list of rational elliptic
surfaces with exactly four semi-stable singular fibers (see
\cite{B} and \cite{Her} for more details). Combine with the
results in \cite{B} and \cite{Her}, we find that two families of
elliptic curves $E_{1, t}$ and $E_{2, t}$ are isogenous. According
to \cite{Shi1}, Theorem 5.2, the points of $E_{1, t}$ and $E_{2,
t}$ over ${\Bbb Q}(t)$ are of finite order because the group of
sections of any elliptic modular surface is finite. In fact,
$E_{1, t}$ is defined over ${\Bbb Q}(t)$, $E_{2, t}$ is defined
over ${\Bbb Q}[t]$.

  Isogenous elliptic curves are related by the modular polynomials
(see \cite{La}). For elliptic curves $E$ and $F$ with
$j$-invariants $j_E$ and $j_F$ respectively, there is an isogeny
$\lambda: E \to F$ with kernel isomorphic to the cyclic group of
order $n$ if and only if
$$\Phi_n(j_E, j_F)=0,$$
where $\Phi_n(x, y)$ is the modular polynomial of order $n$. The
roots of $\Phi_n(X, j(\tau))=0$ are $j(\alpha \tau)$ where
$\alpha=\left(\matrix a & b\\ 0 & d \endmatrix\right)$, $a, b, d
\in {\Bbb Z}$ with $(a, b, d)=1$, $ad=n$ and $0 \leq b<d$. Such
$\alpha$ are called primitive.

  According to \cite{Shi2}, suppose that the minimal Weierstrass equation
of $E$ over $K=k(t)$ is given by
$$y^2+a_1(t) xy+a_3(t) y=x^3+a_2(t) x^2+a_4(t) x+a_6(t), \quad a_i(t) \in k[t].$$
Then the associated elliptic surface $f: S \to {\Bbb P}^1$ is a
rational surface if and only if
$$\deg a_i(t) \leq i \quad \text{(all $i$) and} \quad
  \text{Sing}(f) \neq \emptyset.$$
In case $\text{char}(k) \neq 2, 3$, we can take the equation in a
more familiar form:
$$y^2=x^3+p(t) x+q(t), \quad p(t), q(t) \in k[t].$$
Then the above condition is equivalent to $\deg p(t) \leq 4$,
$\deg q(t) \leq 6$ and $\Delta=4 p(t)^3+27 q(t)^2$ is not a
constant, i.e. $\Delta \notin k$.

  For $E_{2, t}$, we have $p(t)=-3t(t^3+8)$, $q(t)=-2(t^6-20
t^3-8)$ and $\Delta=-2^8 \cdot 3^3 (t^3-1)^3$. Hence, $E_{2, t}$
is a rational elliptic surface.

  We need the following theorems due to Shioda (see \cite{Shi1}).

{\smc Proposition 5.4} (see \cite{Shi1}, Proposition 4.2). {\it
The elliptic modular surface $B_{\Gamma}$ has $t_1$ singular
fibres of types $I_b$ $(b \geq 1)$, $t_2$ singular fibres of types
$I_b^{*}$ $(b \geq 1)$ and $s$ singular fibres of type $IV^{*}$,
where $t_1$, $t_2$ and $s$ are respectively the number of cusps of
the first kind, the number of cusps of the second kind and the
number of elliptic points for $\Gamma$.}

{\smc Theorem 5.5} (see \cite{Shi1}, Theorem 5.2). {\it Let
$B_{\Gamma}$ be the elliptic modular surface attached to $\Gamma$.

\roster

\item If $\Gamma$ has torsion $($i.e. $s>0)$, then the group of
sections ${\Cal S}(B_{\Gamma})$ is either trivial or a cyclic
group of order $3$.

\item If $\Gamma$ has a cusp of the second kind $($i.e. $t_2>0)$,
then the group of sections ${\Cal S}(B_{\Gamma})$ is either
trivial or isomorphic to one of the groups
$${\Bbb Z}/(2), \quad {\Bbb Z}/(4) \quad \text{or} \quad
  {\Bbb Z}/(2) \times {\Bbb Z}/(2).$$

\item If $\Gamma$ is torsion-free and all cusps are of the first
kind, then the group of sections ${\Cal S}(B_{\Gamma})$ is
isomorphic to a subgroup of ${\Bbb Z}/(m) \times {\Bbb Z}/(m)$,
where $m$ denotes the least common multiple of $b_i$'s $(1 \leq i
\leq t_1)$. Here we suppose that the singular fibres of
$B_{\Gamma}$ are of types $I_{b_i}$ $(1 \leq i \leq t_1)$.

\endroster}

  Combine with the above results, we find that the Mordell-Weil
groups for $E_{1, t}$ and $E_{2, t}$ are finite groups. More
precisely, the Mordell-Weil groups of $E_{1, t}$ and $E_{2, t}$
over ${\Bbb Q}(t)$ are isomorphic to the subgroups of ${\Bbb Z}/3
{\Bbb Z} \times {\Bbb Z}/3 {\Bbb Z}$.

  Let $Y \to {\Bbb P}^1$, $Y^{\prime} \to {\Bbb P}^1$ be two
rational elliptic surfaces. Under some mild technical conditions
on the singular fibers, the fibre product $Y \times_{{\Bbb P}^1}
Y^{\prime}$ admits a nice resolution $W$ which is a projective
simply connected threefold with trivial canonical bundle. It can
be showed that some of them are rigid, i.e., admit no nontrivial
deformation (see \cite{Sc}).

  According to \cite{Sc}, let $r: Y \to {\Bbb P}^1$ and $r^{\prime}: Y^{\prime}
\to {\Bbb P}^1$ denote two relatively minimal, rational, elliptic
surfaces with sections. Write $S$ (resp. $S^{\prime}$) for the
images of the singular fibers of $Y$ (resp. $Y^{\prime}$) in
${\Bbb P}^1$. The fiber product $p: W=Y \times_{{\Bbb P}^1}
Y^{\prime} \to {\Bbb P}^1$ is nonsingular except at points in the
fibers over $S^{\prime \prime}=S \cap S^{\prime}$. In order that
the singularities of $W$ be no worse than ordinary double points,
we shall assume that the singular fibers of $r$ and $r^{\prime}$
above the points in $S^{\prime \prime}$ are either irreducible
nodal rational curves or cycles of smooth rational curves. In the
notation of Kodaira such fibers are of type $I_b$, where $b>0$
denotes the number of irreducible components. Such fibers are also
called semi-stable.

  The dualizing sheaf of $W$ is most readily computed by regarding
$W$ as the hypersurface in $Y \times Y^{\prime}$ obtained by
pulling back the diagonal in ${\Bbb P}^1 \times {\Bbb P}^1$ via
the map $r \times r^{\prime}$. For a relatively minimal, regular
elliptic surface, $r: Y \to {\Bbb P}^1$, the canonical sheaf is
given by
$$\omega_Y \simeq r^{*} {\Cal O}_{{\Bbb P}^1}(p_g(Y)-1) \otimes
  {\Cal O}_Y(\sum_{i} (m_i-1) F_i),$$
where the second term is the contribution of the multiple fibers.
An easy computation using the adjunction formula reveals that the
dualizing sheaf $\omega_W$ is trivial exactly when
$p_g(Y)=p_g(Y^{\prime})=0$ and there are no multiple fibers. Of
course these conditions are fulfilled when both $Y$ and
$Y^{\prime}$ are rational and have sections. However it follows
from Castelnuovo's rationality criterion and the vanishing of the
Tate-Shafarevich group for rational elliptic surfaces with
sections that this is the only case when the conditions for
triviality of $\omega_W$ are fulfilled. In this case any small
resolution has trivial canonical bundle.

$$\text{Table 1. Certain rational elliptic modular surfaces (see \cite{Sc})}$$
$$\matrix
  \text{Level} & \text{Congruence subgroup} & \text{Number of components in singular fibers}\\
   3 & \Gamma(3) & 3, 3, 3, 3\\
   4 & \Gamma_1(4) \cap \Gamma(2) & 4, 4, 2, 2\\
   5 & \Gamma_1(5) & 5, 5, 1, 1\\
   6 & \Gamma_1(6) & 6, 3, 2, 1\\
   8 & \Gamma_0(8) \cap \Gamma_1(4) & 8, 2, 1, 1\\
   9 & \Gamma_0(9) \cap \Gamma_1(3) & 9, 1, 1, 1
\endmatrix$$

{\smc Proposition 5.6} (see \cite{Sc}). {\it Let $r: Y \to {\Bbb
P}^1$ and $r^{\prime}: Y^{\prime} \to {\Bbb P}^1$ denote
relatively minimal rational elliptic surfaces with section.
Suppose that the fibers of $r$ and $r^{\prime}$ above all points
in $S^{\prime \prime}$ are semi-stable. Then the fiber product $Y
\times_{{\Bbb P}^1} Y^{\prime}$ has only ordinary double point
singularities. The projective variety, $\widetilde{W}$, obtained
by blowing up the nodes has no infinitesimal deformations, exactly
when all fibers of $r$ and $r^{\prime}$ are semi-stable and one is
in one of the following four cases, each of which actually occurs:
\roster

\item $Y$ and $Y^{\prime}$ appear in Table 1 and are isogenous. In
particular $S=S^{\prime}=S^{\prime \prime}$ has $4$ elements.

\item $Y$ and $Y^{\prime}$ appear in Table 1, each has at least
one $I_1$ fiber, and the map $r^{\prime}$ has been modified by an
automorphism of ${\Bbb P}^1$ so that $\sharp(S^{\prime
\prime})=3$. Furthermore the singular fibers of $r$ and
$r^{\prime}$ which do not lie above points of $S^{\prime \prime}$
are of type $I_1$.

\item $Y$ and $Y^{\prime}$ are not isogenous and
$S=S^{\prime}=S^{\prime \prime}$ has $5$ elements.

\item $Y$ appears in Table 1, $S=S^{\prime \prime}$,
$\sharp(S^{\prime})=5$, and the singular fiber of $r^{\prime}$
which does not map to $S^{\prime \prime}$ has type $I_1$.
\endroster}

  Therefore, we find that the fiber product $E=E_{1, t} \times_{{\Bbb P}^1}
E_{2, t}$ is a rigid Calabi-Yau threefold.

  The fibre product can be considered as a kind of topological
surgery. We will study the structure of rational points under this
topological surgery. From the point of view of Diophantine
geometry, fibre products can be used to generate some new rational
points in Calabi-Yau threefolds. For example, for the self-fiber
product
$$E_{1, t} \times_{{\Bbb P}^1} E_{1, t}: \quad (x^3+y^3+z^3) rst=(r^3+s^3+t^3) xyz,$$
there are some trivial new rational points, such as $x=r$, $y=s$,
$z=t$ and so on. It will be interesting to find some nontrivial
new rational points in the fiber products. We will prove that
there are infinitely many nontrivial new rational points in the
fiber product $E=E_{1, t} \times_{{\Bbb P}^1} E_{2, t}$.

  In actually dealing with specific elliptic curves over ${\Bbb Q}$
we find the following theorem due to Lutz and Nagell very useful
and it gives a manageable method for determining
$E_{\text{tor}}({\Bbb Q})$.

{\smc Theorem 5.7} (Lutz-Nagell) (see \cite{Sil}). {\it Let
$E/{\Bbb Q}$ be an elliptic curve with Weierstrass equation
$$y^2=x^3+Ax+B, \quad A, B \in {\Bbb Z}.$$
Suppose $P \in E({\Bbb Q})$ is a non-zero torsion point. Then
\roster

\item $x(P), y(P) \in {\Bbb Z}$.

\item Either $[2] P=O$, or else $y(P)^2$ divides $4 A^3+27 B^2$.

\endroster}

  For the elliptic curve (5.3), we know that
$P=(3 C_6, 108 C_9) \in E$. Here, $A=-27 C_{12}$ and $B=54
C_{18}$. Note that $A, B \in {\Bbb Z}$ when $(z_1, z_2, z_3) \in
{\Bbb Z}^3$. In fact, $E$ can be defined over ${\Bbb Q}[z_1, z_2,
z_3]$. Suppose that $P$ is a torsion point in $E$, then $P$ is a
torsion point in $E$ with $(z_1, z_2, z_3) \in {\Bbb Z}^3$. We
find that
$$[2] P=\left(-6 C_6+\frac{(C_6^2-C_{12})^2}{64 C_9^2},
        -108 C_9+\frac{9 C_6 (C_6^2-C_{12})}{8 C_9}
        -\frac{(C_6^2-C_{12})^3}{512 C_9^3}\right).\tag 5.20$$
If $[2] P=O$, then $C_9=0$ for every $(z_1, z_2, z_3) \in {\Bbb
Z}^3$, which is impossible! Thus, $[2] P \neq O$. By Lutz-Nagell
theorem, we have that $108^2 C_9^2$ divides $2^8 \cdot 3^{12}
{\frak C}_{12}^3$, i.e., $C_9^2$ divides $2^4 \cdot 3^6 {\frak
C}_{12}^3$ for every $(z_1, z_2, z_3) \in {\Bbb Z}^3$, which is
also impossible! Therefore, $P$ is of infinite order.

  Now, we obtain the following theorem:

{\smc Theorem 5.8} (Main Theorem 5). {\it The variety associated
with Hessian polyhedra
$$E: \quad y^2=x^3-27 C_{12} x+54 C_{18}\tag 5.21$$
can be expressed as the fiber product of two isogenous,
semi-stable, rational elliptic modular surfaces
$$E_{1, t}: \quad (t-1)(z_1^3+z_2^3+z_3^3)-3(t+2) z_1 z_2 z_3=0$$
and
$$E_{2, t}: \quad y^2=x^3-3t(t^3+8)x-2(t^6-20 t^3-8).$$
Hence, $E$ is a rigid Calabi-Yau threefold. There are only
finitely many rational points in $E_{1, t}$ and $E_{2, t}$.
However, there exist infinitely many nontrivial rational points in
$E$. More precisely, $P=(3 C_6, 108 C_9) \in E$ is of infinite
order.}

  The developments in physics stimulated the interest of mathematicians
in Calabi-Yau varieties. There are some important conjectures for
these special varieties. One of these is the modularity conjecture
for Calabi-Yau threefolds defined over ${\Bbb Q}$ (see \cite{Me}
and \cite{Yu}).

  Let us recall some basic facts about Calabi-Yau varieties (see \cite{Yu}).

{\it Definition} 5.9. Let $X$ be a smooth complex projective
variety of dimension $d$. $X$ is called a Calabi-Yau variety if

\roster

\item $H^i(X, {\Cal O}_X)=0$ for every $i$, $0<i<d$.

\item The canonical bundle of $X$ is trivial.

\endroster

  We introduce the Hodge cohomology groups and the Hodge numbers
$$H^{i, j}(X):=H^j(X, \Omega_X^i), \quad h^{i, j}(X):=\dim H^{i, j}(X).$$
By the duality induced from complex conjugation, $h^{i,
j}(X)=h^{j, i}(X)$, and the Serre duality on the Hodge cohomology
groups asserts that $h^{i, j}(X)=h^{d-j, d-i}(X)$. Also there is
the Hodge decomposition
$$H^k(X, {\Bbb C})=\bigoplus_{i+j=k} H^{i, j}(X).$$
The $k$-th Betti number $b_k(X):=\dim H^k(X, {\Bbb C})$ and the
Euler characteristic is
$$\chi(X)=\sum_{k=0}^{2d} (-1)^i b_k(X).$$
The condition (1) is that $h^{i, 0}(X)=0$ for every $i$, $0<i<d$,
and the condition (2) implies that the geometric genus $p_g:=h^{d,
0}(X)=1$.

  The dimension three objects of our interest here are Calabi-Yau threefolds.
Note that Calabi-Yau threefolds are K\"{a}hler manifolds, so that
$h^{1, 1}>0$. The Hodge diamond of a Calabi-Yau threefold is given
by
$$\matrix
   &   &   & 1        &   &   &    \\
   &   & 0 &          & 0 &   &    \\
   & 0 &   & h^{1, 1} &   & 0 &    \\
 1 &   & h^{2, 1} &   & h^{1, 2} &  & 1\\
   & 0 &   & h^{2, 2} &   & 0 &    \\
   &   & 0 &          & 0 &   &    \\
   &   &   & 1        &   &   &
\endmatrix$$
The Betti numbers are
$$b_0=b_6=1, \quad b_1=b_5=0, \quad b_2=b_4=h^{1, 1}, \quad b_3=2(1+h^{2, 1}).$$
The Euler characteristic is $\chi=2(h^{1, 1}-h^{2, 1})$.

{\it Definition} 5.10. A Calabi-Yau threefold $X$ is said to be
rigid if $h^{2, 1}=0$

  The name comes from the fact that this condition implies $H^1(X, T_X)=0$,
and hence that $X$ has no infinitesimal complex deformations. In
this case, the middle cohomology
$$H^3(X)=H^{3, 0}(X) \oplus H^{0, 3}(X) \cong {\Bbb C}^2$$
is two-dimensional.

  Now we will formulate the modularity conjecture for rigid Calabi-Yau
threefolds over ${\Bbb Q}$ (see \cite{Yu}). Let $X$ be a rigid
Calabi-Yau threefold over ${\Bbb Q}$. We know that an integral
model for $X$ always exists. Let $p$ be a good rational prime.
Then the characteristic polynomial of the endomorphism
$\text{Frob}_p^{*}$ on $H_{\text{\'{e}t}}^3(\overline{X}, {\Bbb
Q}_{\ell})$ is
$$P_3^p(T):=\det(1-\text{Frob}_p^{*} T|H_{\text{\'{e}t}}^3(\overline{X}, {\Bbb Q}_{\ell})).$$
Then it is given by
$$P_3^p(T)=1-t_3(p) T+p^3 T^2 \in {\Bbb Z}[T]$$
with $t_3(p) \in {\Bbb Z}$, $|t_3(p)| \leq 2 p^{3/2}$.
Furthermore, by the Lefschetz fixed point formula, $t_3(p)$ can be
expressed in terms of the number of rational points on $X$ over
${\Bbb F}_p$. Let
$t_i(p):=\text{tr}(\text{Frob}_p^{*}|H_{\text{\'{e}t}}^{i}(\overline{X},
{\Bbb Q}_{\ell}))$ denote the trace of the endomorphism
$\text{Frob}_p^{*}$ on $H_{\text{\'{e}t}}^{i}(\overline{X}, {\Bbb
Q}_{\ell})$. Then the Lefschetz fixed point formula asserts that
$$\sharp X({\Bbb F}_p)=\sum_{i=0}^{6} (-1)^i t_i(p).$$
For a rigid Calabi-Yau threefold, this formula gives rise to the
identity:
$$t_3(p)=1+p^3+(1+p) t_2(p)-\sharp X({\Bbb F}_p)$$
where $t_2(p) \leq p h^{1, 1}$ and the equality holds when all
cycles generating $H^{1, 1}(X)$ are defined over ${\Bbb Q}$. The
$L$-series of $X$ is given by putting together all the local data:
$$L(X_{{\Bbb Q}}, s):=L(H_{\text{\'{e}t}}^3(\overline{X}, {\Bbb Q}_{\ell}), s)
 =(\ast) \prod_{p: \text{good}} \frac{1}{1-t_3(p) p^{-s}+p^{3-2s}}$$
where $(\ast)$ is the factor corresponding to bad primes.

  Now, we give the modularity conjecture for rigid Calabi-Yau threefolds over ${\Bbb Q}$.

{\smc Conjecture 5.11} (see \cite{Yu}). {\it Let $X$ be a rigid
Calabi-Yau threefold defined over ${\Bbb Q}$. Then there is a cusp
form $f$ of weight $4$ on some $\Gamma_0(N)$ where $N$ is
divisible only by bad primes, such that the two $2$-dimensional
$\ell$-adic Galois representations $\rho_{X, \ell}$ and $\rho_{f,
\ell}$ associated to $X$ and $f$, respectively, have the
isomorphic semi-simplifications $(\rho_{X, \ell} \sim \rho_{f,
\ell})$. Consequently,
$$L(X, s)=L(f, s).$$
More precisely, if we write $f(q)=\sum_{m=1}^{\infty} a_f(m) q^m$
with $q=e^{2 \pi i z}$, then
$$t_3(p)=a_f(p)$$
for all primes not dividing $N$.}

  The fiber products of rational elliptic surfaces give some examples
of rigid Calabi-Yau threefolds and their modularity. It is known
that (see \cite{B}, \cite{Sc} and \cite{Seb}) there are exactly
six torsion free genus zero congruence subgroups of $PSL(2, {\Bbb
Z})$, and the list is complete. They are $\Gamma(3)$, $\Gamma_1(4)
\cap \Gamma(2)$, $\Gamma_1(5)$, $\Gamma_1(6)$, $\Gamma_0(8) \cap
\Gamma_1(4)$ and $\Gamma_0(9) \cap \Gamma_1(3)$ in the Table 1.
The modularity of the rigid Calabi-Yau threefolds defined over
${\Bbb Q}$ corresponding to these subgroups are given as follows:

{\smc Theorem 5.12}. (see \cite{SY}, \cite{Ve1}) {\it Let $\Gamma$
be one of the six groups in the Table 1. Then the rigid Calabi-Yau
threefolds defined over ${\Bbb Q}$ associated to $\Gamma$ are all
modular. More precisely, let
$\widetilde{S_{\Gamma}^2}:=\widetilde{S_{\Gamma}
\times_{C_{\Gamma}} S_{\Gamma}}$ be a rigid Calabi-Yau threefold
defined over ${\Bbb Q}$ associated to $\Gamma$. Then there is a
cusp form of weight $4$ on $\Gamma$ such that
$L(\widetilde{S_{\Gamma}^2}, s)=L(f, s)$.}

  In particular, when $\Gamma=\Gamma(3)$, the singular fibers are of
types $(I_3, I_3, I_3, I_3)$. The Hodge number $h^{1, 1}(X)=36$.
The cusp form of weight $4$ is given by $\eta(q^3)^8$. Here
$\eta(q)$ is Dedekind eta-function $\eta(q)=q^{1/24}
\prod_{n=1}^{\infty} (1-q^n)$ with $q=e^{2 \pi i z}$, $z \in {\Bbb
H}$. The equation for the rigid Calabi-Yau threefold is
$$(x^3+y^3+z^3) rst=(r^3+s^3+t^3) xyz.$$

  Note that in the above construction, we only use the self-fiber
products of semi-stable, rational elliptic modular surfaces. It
will be interesting to study the fiber products of two isogenous,
semi-stable, rational elliptic modular surfaces, such as our fiber
product $E_{1, t} \times_{{\Bbb P}^1} E_{2, t}$. A natural
question is: which modular form corresponds to this fiber product?

  According to Vafa (see \cite{V}), sometimes the physics of the singularities
are unconventional. For example when a $4$-cycle (say a ${\Bbb C}
{\Bbb P}^2$) shrinks in a Calabi-Yau threefold, it gives rise to
very interesting unconventional new physical theories which were
not anticipated! This is thus a great source of insight into new
physics. In particular what types of singularities occur as well
as what are the ways to resolve them will be of extreme importance
for unravelling aspects of this new physics. It is tempting to
speculate that these singularities may also lead to new invariants
for four manifolds.

  For the rigid Calabi-Yau threefold $E$ in (5.21) where
$C_{12}=C_{12}(z_1, z_2, z_3)$ and $C_{18}=C_{18}(z_1, z_2, z_3)$
with $(z_1, z_2, z_3) \in {\Bbb C} {\Bbb P}^2$, we construct
infinitely many subvarieties of $E$:
$$Y_m=\{ (z_1, z_2, z_3, x([m] P), y([m] P)): (z_1, z_2, z_3) \in {\Bbb C} {\Bbb P}^2,
         P=(3 C_6, 108 C_9) \in E, m \in {\Bbb Z} \}.\tag 5.22$$
It will be interesting to study the singularity when $Y_m$ shrinks
in $E$. This singularity maybe involves vanishing cycles and
Picard-Lefschetz theory.

  Note that the Calabi-Yau threefold $E$ in (5.21) can be considered
as an elliptic curve defined over ${\Bbb Q}(T_1, T_2)$:
$$E/{\Bbb Q}(T_1, T_2): \quad y^2=x^3-27 C_{12} x+54 C_{18},\tag 5.23$$
where $T_1=z_1/z_3$, $T_2=z_2/z_3$ and
$$\left\{\aligned
  C_{12}=C_{12}(T_1, T_2) &=(T_1^3+T_2^3+1)[(T_1^3+T_2^3+1)^3+216 T_1^3 T_2^3],\\
  {\frak C}_{12}={\frak C}_{12}(T_1, T_2) &=T_1 T_2 [27 T_1^3 T_2^3-(T_1^3+T_2^3+1)^3],\\
  C_{18}=C_{18}(T_1, T_2) &=(T_1^3+T_2^3+1)^6-540 T_1^3 T_2^3 (T_1^3+T_2^3+1)^3-5832 T_1^6 T_2^6.
\endaligned\right.\tag 5.24$$

  Let $k$ be a number field and $A$ an abelian variety over
$K=k(T_1, \cdots, T_n)$ where the $T_i$ are indeterminates over
$k$. By a theorem of N\'{e}ron (see \cite{Se}), $A(K)$ is a
finitely generated group. In our case, $E({\Bbb Q}(T_1, T_2))$ is
finitely generated by N\'{e}ron's theorem. The $j$-invariant
$$j=j(T_1, T_2)=-\frac{C_{12}(T_1, T_2)^3}{{\frak C}_{12}(T_1, T_2)^3} \notin {\Bbb Q}.\tag 5.25$$
In fact, $E$ is an elliptic threefold (see \cite{Mi} for more
details).

\vskip 2.0 cm

{\smc Department of Mathematics, Peking University}

{\smc Beijing 100871, P. R. China}

{\it E-mail address}: yanglei\@math.pku.edu.cn
\vskip 1.5 cm
\Refs

\item{[Ar1]} {\smc E. Artin}, \"{U}ber eine neue Art von L-Reihen,
            Hamb. Abh. {\bf 3} (1923), 89-108, In: {\it The
            Collected Papers of Emil Artin}, Edited by S. Lang and
            J. Tate, 105-124, Addison-Wesley Publishing Company,
            1965.

\item{[Ar2]} {\smc E. Artin}, Zur Theorie der L-Reihen mit
            allgemeinen Gruppencharakteren, Hamb. Abh. {\bf 8} (1930),
            292-306, In: {\it The Collected Papers of Emil Artin},
            Edited by S. Lang and J. Tate, 165-179, Addison-Wesley
            Publishing Company, 1965.

\item{[B]} {\smc A. Beauville}, Les familles stables de courbes
            elliptiques sur ${\Bbb P}^1$ admettant quatre fibres
            singuli\`{e}res, C. R. Acad. Sci. Paris S\'{e}r. I
            Math. {\bf 294} (1982), 657-660.

\item{[Br]} {\smc R. Brauer}, On simple groups of order $5 \cdot
            3^a \cdot 2^b$, In: {\it Collected Papers}, Vol. II,
            Edited by P. Fong and W. J. Wong, 421-470, The MIT
            Press, 1970.

\item{[BE]} {\smc E. Brieskorn}, Singularities in the work of
           Friedrich Hirzebruch, Surveys in Differential Geometry, Vol. VII (2000), 17-60.

\item{[BrE]} {\smc M. Brou\'{e} and M. Enguehard}, Polyn\^{o}mes
            des poids de certains codes et fonctions th\^{e}ta de
            certains r\'{e}seaux, Ann. Sci. \'{E}cole Norm. Sup.
            (4) {\bf 5} (1972), 157-181.

\item{[Bu1]} {\smc H. Burkhardt}, Untersuchungen aus dem Gebiete
            der hyperelliptischen Modulfunctionen, II, Math. Ann.
            {\bf 38} (1891), 161-224.

\item{[Bu2]} {\smc H. Burkhardt}, Untersuchungen aus dem Gebiete
            der hyperelliptischen Modulfunctionen, III, Math. Ann.
            {\bf 41} (1893), 313-343.

\item{[CMS]} {\smc Y. Choie, B. Mourrain and P. Sol\'{e}},
             Rankin-Cohen brackets and invariant theory, J.
             Algebraic Combinatorics {\bf 13} (2001), 5-13.

\item{[CC]} {\smc J. H. Conway, R. T. Curtis, S. P. Norton, R.
             A. Parker and R. A. Wilson}, {\it Atlas of Finite Groups,
             Maximal Subgroups and Ordinary Characters for Simple Groups},
             Clarendon Press, Oxford, 1985.

\item{[CS]} {\smc J. H. Conway and N. J. A. Sloane}, {\it Sphere
            Packings, Lattices and Groups}, Second Edition,
            Springer-Verlag, 1993.

\item{[Co]} {\smc H. S. M. Coxeter}, {\it Regular Complex
            Polytopes}, Cambridge University Press, 1974.

\item{[Eb]} {\smc W. Ebeling}, {\it Lattices and Codes, A Course
             Partially Based on Lectures by F. Hirzebruch}, 2nd Edition,
             Advanced Lectures in Mathematics, Vieweg, 2002.

\item{[F]} {\smc R. Fricke}, Lehrbuch der Algebra, Bd. II,
           Braunschweig, 1926.

\item{[G]} {\smc S. Gelbart}, An elementary introduction to the
           Langlands program, Bull. Amer. Math. Soc. {\bf 10}
           (1984), 177-219.

\item{[Her]} {\smc S. Herfurtner}, Elliptic surfaces with four
           singular fibres, Math. Ann. {\bf 291} (1991), 319-342.

\item{[H]} {\smc D. Hilbert}, Mathematical problems, Bull. Amer.
            Math. Soc. {\bf 8} (1902), 437-479.

\item{[Hi]} {\smc F. Hirzebruch}, {\it Gesammelte Abhandlungen},
            Bd. II, 1963-1987, Springer-Verlag, 1987.

\item{[Hi1]} {\smc F. Hirzebruch}, Hilbert's modular group of the
             field ${\Bbb Q}(\sqrt{5})$ and the cubic diagonal surface
             of Clebsch and Klein, Russ. Math. Surv. {\bf 31}: 5(1976),
             96-110, In: {\it Gesammelte Abhandlungen}, Bd. II, 394-408,
             Springer-Verlag, 1987.

\item{[Hi2]} {\smc F. Hirzebruch}, The ring of Hilbert modular
            forms for real quadratic fields of small discriminant, In:
            {\it Modular Functions of One Variable VI, Proceedings,
            Bonn 1976}, Edited by J. P. Serre and D. B. Zagier,
            287-323, Lecture Notes in Math. {\bf 627}, Springer-Verlag,
            1977, In: {\it Gesammelte Abhandlungen}, Bd. II, 501-536,
            Springer-Verlag, 1987.

\item{[Hi3]} {\smc F. Hirzebruch}, The icosahedron, Raymond and
             Beverly Sackler Distinguished Lectures in Mathematics,
             Tel Aviv University, 1981, In: {\it Gesammelte Abhandlungen},
             Bd. II, 656-661, Springer-Verlag, 1987.

\item{[Ho1]} {\smc R. P. Holzapfel}, {\it Geometry and Arithmetic
            Around Euler Partial Differential Equations}, Dt. Verl.
            d. Wiss., Berlin/Reidel, Dordrecht 1986.

\item{[Ho2]} {\smc R. P. Holzapfel}, {\it The Ball and Some
            Hilbert Problems}, Lectures in Mathematics, ETH Z\"{u}rich,
            Birkh\"{a}user, 1995.

\item{[Ho3]} {\smc R. P. Holzapfel}, On the Nebentypus of Picard
             modular forms, Math. Nachr. {\bf 139} (1988), 115-137.

\item{[Ho4]} {\smc R. P. Holzapfel}, Chern numbers of algebraic
             surfaces-Hirzebruch's examples are Picard modular
             surfaces, Math. Nachr. {\bf 126} (1986), 255-273.

\item{[Hu]} {\smc B. Hunt}, {\it The Geometry of some special
            Arithmetic Quotients}, Lecture Notes in Math. {\bf
            1637}, Springer-Verlag, 1996.

\item{[Hus]} {\smc D. Husem\"{o}ller}, {\it Elliptic Curves}, GTM
             {\bf 111}, Springer-Verlag, 1987.

\item{[Kl]} {\smc F. Klein}, {\it Lectures on the Icosahedron and
             the Solution of Equations of the Fifth Degree},
             Translated by G. G. Morrice, second and revised edition,
             Dover Publications, Inc., 1956.

\item{[La]} {\smc S. Lang}, {\it Elliptic Functions}, GTM {\bf
            112}, Springer-Verlag, 1987.

\item{[L1]} {\smc R. P. Langlands}, Some contemporary problems
            with origins in the Jugendtraum, In: {\it Mathematical
            developments arising from Hilbert problems} (Proc.
            Sympos. Pure Math., Vol. XXVIII, Northern Illinois
            Univ., De Kalb, Ill., 1974), 401-418, Amer. Math.
            Soc., Providence, R. I., 1976.

\item{[L2]} {\smc R. P. Langlands}, {\it Base Change for $GL(2)$},
            Ann. of Math. Studies {\bf 96}, Princeton University
            Press, 1980.

\item{[Ma]} {\smc Yu. I. Manin}, {\it Cubic Forms, Algebra,
            Geometry, Arithmetic}, Translated from Russian by
            M. Hazewinkel, Second Edition, North-Holland, 1986.

\item{[MT]} {\smc Yu. I. Manin and M. A. Tsfasman}, Rational
            varieties: algebra, geometry and arithmetic, Russian
            Math. Surveys {\bf 41}:2 (1986), 51-116.

\item{[MP]} {\smc Yu. I. Manin and A. A. Panchishkin}, {\it
            Introduction to Modern Number Theory}, Second Edition,
            Springer-Verlag, 2005.

\item{[Mas]} {\smc H. Maschke}, Aufstellung des vollen
            Formensystems einer quatern\"{a}ren Gruppe von $51840$ linearen
            Substitutionen, Math. Ann. {\bf 33} (1889), 317-344.

\item{[MS]} {\smc J. McKay and A. Sebbar}, Arithmetic semistable
            elliptic surfaces, In: {\it Proceedings on Moonshine
            and Related Topics}, Edited by J. McKay and A. Sebbar,
            CRM Proceedings and Lecture Notes, Vol. {\bf 30},
            119-130, AMS Providence, 2001.

\item{[Me]} {\smc C. Meyer}, {\it Modular Calabi-Yau Threefolds},
            Fields Institute Monographs, Vol. {\bf 22}, AMS Providence, 2005.

\item{[Mi]} {\smc R. Miranda}, Smooth models for elliptic
            threefolds, In: {\it The Birational Geometry of Degenerations},
            Edited by R. Friedman and D. R. Morrison, Progress in
            Math. {\bf 29}, 85-133, Birkh\"{a}user, 1983.

\item{[Pe]} {\smc C. Peters}, Monodromy and Picard-Fuchs equations
            for families of K3-surfaces and elliptic curves, Ann.
            Sci. \'{E}cole Norm. Sup. (4) {\bf 19} (1986), 583-607.

\item{[Pi]} {\smc E. Picard}, Sur les formes quadratiques
            ternaires ind\'{e}finies \`{a} ind\'{e}termin\'{e}es conjugu\'{e}es
            et sur les fonctions hyperfuchsiennes correspondantes,
            Acta Math. {\bf 5} (1884), 121-184.

\item{[SY]} {\smc M. -H. Saito and N. Yui}, The modularity
            conjecture for rigid Calabi-Yau threefolds over ${\Bbb Q}$, J.
            Math. Kyoto Univ. {\bf 41} (2001), 403-419.

\item{[Sc]} {\smc C. Schoen}, On fiber products of rational
            elliptic surfaces with section, Math. Z. {\bf 197} (1988), 177-199.

\item{[Seb]} {\smc A. Sebbar}, Classification of torsion-free
             genus zero congruence groups, Proc. Amer. Math. Soc.
             {\bf 129}(2001), 2517-2527.

\item{[Seg]} {\smc B. Segre}, {\it The Nonsingular Cubic Surfaces,
             A New Method of Investigation with Special Reference to
             Questions of Reality}, Oxford, 1942.

\item{[Se]} {\smc J.-P. Serre}, {\it Lectures on the Mordell-Weil
             Theorem}, Translated from the French and edited by
             Martin Brown from notes by Michel Waldschmidt,
             Aspects of Mathematics, E15, Friedr. Vieweg $\&$ Sohn,
             Braunschweig, 1989.

\item{[Shi1]} {\smc T. Shioda}, On elliptic modular surfaces, J.
            Math. Soc. Japan {\bf 24} (1972), 20-59.

\item{[Shi2]} {\smc T. Shioda}, On the Mordell-Weil lattices,
            Comment. Math. Univ. St. Pauli {\bf 39} (1990), 211-240.

\item{[Shi3]} {\smc T. Shioda}, Theory of Mordell-Weil lattices,
            Proc. ICM Kyoto, Vol. I (1990), 473-489.

\item{[Sil]} {\smc J. H. Silverman}, {\it The Arithmetic of
            Elliptic Curves}, GTM {\bf 106}, Springer-Verlag, 1986.

\item{[Ta]} {\smc J. Tate}, Problem 9: The general reciprocity
            law, In: {\it Mathematical developments arising from
            Hilbert problems} (Proc. Sympos. Pure Math., Vol. XXVIII,
            Northern Illinois Univ., De Kalb, Ill., 1974), 311-322,
            Amer. Math. Soc., Providence, R. I., 1976.

\item{[V]} {\smc C. Vafa}, Geometric physics, Proc. ICM,
           Vol. I (Berlin, 1998), Doc. Math. 1998, Extra Vol. I, 537-556.

\item{[Ve1]} {\smc H. A. Verrill}, The $L$-series of certain rigid
            Calabi-Yau threefolds, J. Number Theory {\bf 81}
            (2000), 310-334.

\item{[Ve2]} {\smc H. A. Verrill}, Picard-Fuchs equations of some
            families of elliptic curves, In: {\it Proceedings on Moonshine
            and Related Topics}, Edited by J. McKay and A. Sebbar,
            CRM Proceedings and Lecture Notes, Vol. {\bf 30},
            253-268, AMS Providence, 2001.

\item{[Wei1]} {\smc A. Weil}, Arithmetic on algebraic varieties,
             In: {\it Collected Papers}, Vol. I, 116-126, Springer-Verlag, 1979.

\item{[Wei2]} {\smc A. Weil}, Abstract versus classical algebraic
             geometry, Proc. ICM Amsterdam, Vol. III (1954), 550-558,
             In: {\it Collected Papers}, Vol. II, 180-188, Springer-Verlag, 1979.

\item{[Wei3]} {\smc A. Weil}, Two lectures on number theory, past
             and present, In: {\it Collected Papers}, Vol. III, 279-302, Springer-Verlag, 1979.

\item{[Y1]} {\smc L. Yang}, Product formulas on a unitary group in
            three variables, math.NT/0104259.

\item{[Y2]} {\smc L. Yang}, Geometry and arithmetic associated to
           Appell hypergeometric partial differential equations,
           math.NT/0309415.

\item{[Y3]} {\smc L. Yang}, Hessian polyhedra, invariant theory
           and Appell hypergeometric partial differential equations,
           math.NT/0412065.

\item{[Yu]} {\smc N. Yui}, Update on the modularity of Calabi-Yau
            varieties, In: {\it Calabi-Yau Varieties and Mirror
            Symmetry}, Edited by N. Yui and J. D. Lewis, Fields
            Institute Communications Vol. {\bf 38} (2003),
            307-362, AMS Providence, 2003.

\endRefs
\end{document}